\documentclass[11pt]{article}

\usepackage{latexsym}
\usepackage{amssymb}
\usepackage{amsthm}
\usepackage{amscd}
\usepackage{amsmath}
\usepackage{mathrsfs}
\usepackage{graphicx}
\usepackage{hyperref}
\usepackage{ytableau}
\usepackage{stmaryrd}

\usepackage[colorinlistoftodos]{todonotes}
\usetikzlibrary{chains,scopes,decorations.markings}

\theoremstyle{definition}
\newtheorem* {theorem*}{Theorem}
\newtheorem* {proposition*}{Proposition}
\newtheorem* {conjecture*}{Conjecture}
\newtheorem{theorem}{Theorem}[section]

\theoremstyle{definition}

\newtheorem* {example*}{Example}

\newtheorem{lemma}[theorem]{Lemma}
\theoremstyle{definition}
\newtheorem{definition}[theorem]{Definition}
\theoremstyle{definition}

\newtheorem{proposition}[theorem]{Proposition}
\newtheorem{corollary}[theorem]{Corollary}

\newtheorem*{remark*}{Remark}
\theoremstyle{definition}
\newtheorem{remark}[theorem]{Remark}
\theoremstyle{definition}
\newtheorem {example}[theorem]{Example}
\theoremstyle{definition}

\theoremstyle{definition}

\theoremstyle{definition}

\theoremstyle{definition}

\newcommand{\cC}{\mathcal{C}}

\newcommand{\NN}{\mathbb{Z}_{\geq0}}
\newcommand{\PP}{\mathbb{Z}_{>0}}
\newcommand{\ZZ}{\mathbb{Z}}

\newcommand{\barr}{\begin{array}}
\newcommand{\earr}{\end{array}}
\newcommand{\ba}{\begin{aligned}}
\newcommand{\ea}{\end{aligned}}
\newcommand{\be}{\begin{equation}}
\newcommand{\ee}{\end{equation}}

\newcommand{\quand}{\quad\text{and}\quad}

\newcommand{\quord}{\quad\text{or}\quad}

\newcommand{\ben}{\begin{enumerate}}
\newcommand{\een}{\end{enumerate}}

\newcommand{\bei}{\begin{itemize}}
\newcommand{\eei}{\end{itemize}}

\newcommand{\e}{\varepsilon}

\definecolor{darkred}{rgb}{0.7,0,0} 
\newcommand{\defn}[1]{{\color{darkred}\emph{#1}}} 

\newcommand{\SD}{\mathsf{SD}}

\newcommand{\SetShTab}{\mathsf{SetShTab}}

\newcommand{\GP}{G\hspace{-0.25mm}P}
\newcommand{\GQ}{G\hspace{-0.25mm}Q}
\newcommand{\GS}{G\hspace{-0.25mm}S}

\numberwithin{equation}{section}
\allowdisplaybreaks[1]

\newcommand{\res}{\mathsf{res}}
\newcommand{\sfA}{\mathsf{A}}
\newcommand{\sfB}{\mathsf{B}}
\newcommand{\sfC}{\mathsf{C}}
\newcommand{\sfD}{\mathsf{D}}

\newcommand{\mpmap}{\mathsf{mperm}}

\newcommand{\FM}{\mathsf{FSeq}}
\newcommand{\MP}{\mathsf{MPerm}}

\newcommand{\TBC}{T^{\sfB}}
\newcommand{\TD}{T^{\sfD}}

\newcommand{\cK}[1]{\Theta_{#1}}
\newcommand{\cKA}[1]{\cK{#1}^{\sfA}}
\newcommand{\cKB}[1]{\cK{#1}^{\sfB}}
\newcommand{\cKC}[1]{\cK{#1}^{\sfC}}
\newcommand{\cKD}[1]{\cK{#1}^{\sfD}}
\newcommand{\cKX}[1]{\cK{#1}^X}

\newcommand{\cH}{\mathcal{H}}

\newcommand{\UU}{\mathbb{U}}
\newcommand{\bVV}{\widehat{\mathbb{V}}}
\newcommand{\fkS}{\mathfrak{S}}

\newcommand{\cF}{G}
\newcommand{\cFX}{\cF^X}
\newcommand{\cFB}{\cF^{\sfB}}
\newcommand{\cFC}{\cF^{\sfC}}
\newcommand{\cFD}{\cF^{\sfD}}

\newcommand{\VV}{\mathbb{V}}

\renewcommand{\AA}{\mathbb{A}_\beta}

\newcommand{\inv}{\mathsf{inv}}

\newcommand{\QQ}{\mathbb{Q}}

\newcommand{\bfQB}{\widetilde{\mathfrak{R}}^\sfB}
\newcommand{\bfu}{\mathfrak{u}}
\newcommand{\bft}{\mathfrak{t}}
\newcommand{\bfv}{\mathfrak{v}}
\newcommand{\bftB}{\bft^{\sfB}}
\newcommand{\bftC}{\bft^{\sfC}}
\newcommand{\bftD}{\bft^{\sfD}}
\newcommand{\bfuB}{\bfu^{\sfB}}
\newcommand{\bfuC}{\bfu^{\sfC}}
\newcommand{\bfuD}{\bfu^{\sfD}}
\newcommand{\bfvX}{\bfv^X}

\newcommand{\bfvB}{\bfv^{\sfB}}

\newcommand{\bftX}{\bft^X}
\newcommand{\bfuX}{\bfu^X}
\newcommand{\bfM}{\mathfrak{M}}
\newcommand{\bfMX}{\bfM^X}
\newcommand{\bfMB}{\bfM^{\sfB}}

\newcommand{\bfo}{\mathfrak{o}}
\newcommand{\bfoX}{\bfo^{X}}
\newcommand{\bfoB}{\bfo^{\sfB}}
\newcommand{\bfoC}{\bfo^{\sfC}}
\newcommand{\bfoD}{\bfo^{\sfD}}
\newcommand{\bfn}{\mathfrak{n}}
\newcommand{\bfnX}{\bfn^X}
\newcommand{\bfnB}{\bfn^{\sfB}}
\newcommand{\bfnC}{\bfn^{\sfC}}
\newcommand{\bfnD}{\bfn^{\sfD}}
\newcommand{\bfR}{\mathfrak{R}}
\newcommand{\bfRX}{\bfR^X}
\newcommand{\bfRB}{\bfR^{\sfB}}

\newcommand{\fku}{\mathbf{u}}
\newcommand{\fkt}{\mathbf{t}}
\newcommand{\fkv}{\mathbf{v}}
\newcommand{\fktB}{\fkt^{\sfB}}

\newcommand{\fktD}{\fkt^{\sfD}}
\newcommand{\fkuB}{\fku^{\sfB}}

\newcommand{\fkuD}{\fku^{\sfD}}
\newcommand{\fkvX}{\fkv^X}
\newcommand{\fkvA}{\fkv^{\sfA}}
\newcommand{\fkvB}{\fkv^{\sfB}}

\newcommand{\fktX}{\fkt^X}
\newcommand{\fkuX}{\fku^X}
\newcommand{\fkM}{\mathbf{M}}
\newcommand{\fkMX}{\fkM^X}
\newcommand{\fkMA}{\fkM^{\sfA}}
\newcommand{\fkMB}{\fkM^{\sfB}}

\newcommand{\fktA}{\fkt^{\sfA}}
\newcommand{\fkuA}{\fku^{\sfA}}
\newcommand{\fko}{\mathbf{o}}
\newcommand{\fkoX}{\fko^{X}}
\newcommand{\fkoB}{\fko^{\sfB}}
\newcommand{\fkoC}{\fko^{\sfC}}
\newcommand{\fkoD}{\fko^{\sfD}}

\newcommand{\Des}{\mathrm{Des}}
\newcommand{\Peak}{\mathrm{Peak}}

\newcommand{\ellA}{\ell^{\sfA}}
\newcommand{\ellB}{\ell^{\sfB}}
\newcommand{\ellC}{\ell^{\sfC}}
\newcommand{\ellD}{\ell^{\sfD}}

\newcommand{\LD}{\mathsf{LD}}

\newcommand{\fkGX}{\fkG^X}
\newcommand{\fkGB}{\fkG^{\sfB}}
\newcommand{\fkGC}{\fkG^{\sfC}}
\newcommand{\fkGD}{\fkG^{\sfD}}

\newcommand{\last}{\mathsf{end}}
\newcommand{\cO}{\mathcal{O}} 
\newcommand{\cR}{\mathcal{R}}
\newcommand{\cS}{\mathcal{S}}
\newcommand{\cL}{\mathcal{L}}

\newcommand{\raff}{r^{\mathrm{aff}}}

\newcommand{\cU}{\mathcal{U}}
\newcommand{\CC}{\mathbb{C}}
\newcommand{\RR}{\mathbb{R}}

\newcommand{\spanning}{\text{-}\mathrm{span}}

\newcommand{\fkh}{\mathfrak{h}}
\newcommand{\fkg}{\mathfrak{g}}

\newcommand{\GL}{\mathsf{GL}}
\newcommand{\SL}{\mathsf{SL}}
\newcommand{\SO}{\mathsf{SO}}
\newcommand{\Sp}{\mathsf{Sp}}

\newcommand{\fkG}{\mathfrak{G}}

\newcommand{\bfx}{\mathbf{x}}
\newcommand{\bfy}{\mathbf{y}}
\newcommand{\bfz}{\mathbf{z}}
\renewcommand{\(}{\left(}
\renewcommand{\)}{\right)}

\newcommand{\WA}{W^{\sfA}}
\newcommand{\WB}{W^{\sfB}}
\newcommand{\WC}{W^{\sfC}}
\newcommand{\WD}{W^{\sfD}}

\newcommand{\ILambda}{I\hspace{-0.2mm}\Lambda}

\newcommand{\QSym}{\mathbf{QSym}}
\newcommand{\Sym}{\mathbf{Sym}}
\newcommand{\GPP}{\mathbf{GP}^+}
\newcommand{\GQP}{\mathbf{GQ}^+}
\renewcommand{\SS}{\mathbf{SSym}}
\newcommand{\bfG}{\mathbf{G}}

\newcommand{\GPB}{\mathbf{GP}}
\newcommand{\GQB}{\mathbf{GQ}}

\usepackage{fullpage}

\begin{document}
\title{Classical double Grothendieck transitions}
\author{
    Eric Marberg\\
    Department of Mathematics \\
    Hong Kong University of Science and Technology \\
    {\tt emarberg@ust.hk}
}

\date{}

\maketitle

\begin{abstract}
Kirillov and Naruse have constructed double Grothendieck polynomials to represent the 
equivariant $K$-theory classes of Schubert varieties in the complete flag manifolds of types $B$, $C$, and $D$.
We derive a recursive formula for these polynomials,
extending certain $K$-theoretic transition equations known in type $A$ to all classical types.
As an application, we obtain an identity that expands the $K$-Stanley symmetric functions in types $B$, $C$, and $D$
into positive linear combinations of $K$-theoretic Schur $P$- and $Q$-functions.
We also resolve several positivity conjectures related to the  skew generalizations of the latter  functions.
\end{abstract}


\section{Introduction}

This article studies certain formal power series introduced by Kirillov and Naruse \cite{KN}
in connection with the equivariant $K$-theory rings of  flag varieties in classical type.
We show that the \defn{$K$-Stanley symmetric functions} for types $\sfB$, $\sfC$, and $\sfD$
defined in \cite{KN}
 have finite, positive  expansions into the \defn{$K$-theoretic Schur $P$- and $Q$-functions}
constructed by Ikeda and Naruse in \cite{IN}.
Before explaining these statements precisely, we review some analogous properties in type $\sfA$ that serve as motivation.

\subsection{Symmetric Grothendieck functions}

Write $\ZZ$ for the set of integers, and
let $\beta$ and $z_1,z_2,z_3,\dots$ be commuting variables.
We define $\Sym_\beta$ to be the subring of formal power series in $\ZZ[\beta]\llbracket \bfz\rrbracket = \ZZ[\beta]\llbracket z_1,z_2,z_3,\dots\rrbracket$ that are invariant under all permutations of the $z$-variables.
Let $S_\infty$ denote the group of permutations of the positive integers with finite support
and write $\ell : S_\infty \to \NN$ for its Coxeter length function, which assigns to each permutation its number of inversions.

Each element $w \in S_\infty$ has an associated \defn{symmetric Grothendieck function} $G_w(\bfz)\in \Sym_\beta$,
which is a certain symmetric formal power series of unbounded degree.
 These functions, which were first defined by Fomin and Kirillov in \cite{FK1994},
 arise geometrically as the stable limits of the \defn{Grothendieck polynomials} of Lascoux and Sch\"utzenberger \cite{LS1983},
 and are sometimes called  \defn{stable Grothendieck polynomials}.
They also
 generalize the  \defn{Stanley symmetric functions} from \cite{Stanley}, which are recovered by setting $\beta=0$.
 For the precise definition, see Section~\ref{k-stanley-sect}.

There is a natural bialgebra structure on the module
 \be
 \bfG:=\ZZ[\beta]\spanning  \{ G_w : w \in S_\infty\}
 \ee
of finite $\ZZ[\beta]$-linear combinations of symmetric Grothendieck functions.
 To describe this, given 
  $u,v \in S_\infty$ with
$n = \max( \{0\}\sqcup \{ i \in \PP : u(i)\neq i\})$, write $u\times v \in S_\infty$ for the element defined by 
 \[
 (u\times v)( i) = \begin{cases} u(i) &\text{if }i \in \{1,2,\dots,n\} \\ 
 v(i-n) +n &\text{if }i \in \{n+1,n+2,n+3,\dots\}.
 \end{cases}
 \]
 In addition, let $\circ$ denote the \defn{Demazure product} for $S_\infty$, which is
 the unique associative binary operation
with
$u\circ v = uv$ if $\ell(uv)=\ell(u)+\ell(v)$ and $s\circ s = s$ if $\ell(s)=1$.

\begin{proposition}[\cite{Buch2002}] 
\label{dem-prop}
The module $ \bfG$ is a bialgebra
for the $\ZZ[\beta]$-linear (co)product operations
\[  G_u(\bfz)  G_v(\bfz)  = G_{u\times v} (\bfz)
\quand \Delta(G_w(\bfz) ) = \sum_{\substack{u,v\in S_\infty \\ u\circ v = w}}  \beta^{\ell(u)+\ell(v)-\ell(w)}G_u(\bfz) \otimes G_v(\bfz),\]
which extend the usual (co)product on the Hopf algebra of bounded degree symmetric functions.
\end{proposition}

\begin{proof}
This result is implicit in \cite{Buch2002}. The 
product formula follows as an exercise from  \cite[Eq.~(7)]{BKSTY},
and the 
coproduct formula can be obtained by specializing \cite[Thm.~4.15]{LLS}.
\end{proof}

A permutation $w \in S_\infty$ is \defn{Grassmannian} if 
 it has at most one \defn{descent} $i \in \PP$ with $w(i)>w(i+1)$; in this case the \defn{shape} of $w$
is the partition $\lambda = (w(i) - i, \dots, w(2) -2, w(1)-1)$.
All Grassmannian permutations with a given shape $\lambda$ have the same symmetric Grothendieck function,
which we denote by $G_\lambda(\bfz)$. 
The subset of symmetric Grothendieck functions obtained in this way is significant for the following reasons:
\bei
\item Buch \cite{Buch2002} showed that  $G_\lambda(\bfz)$ has an explicit tableau generating function formula. Specifically,
  $G_\lambda(\bfz)$  is
the weight generating function for all \defn{semistandard set-valued tableaux} of shape $\lambda$.
One can deduce from this formula that the $G_\lambda(\bfz)$'s are linearly independent over $\ZZ[\beta]$.

\item Lascoux \cite{LascouxTransitions} proved that for each $w \in S_\infty$,
 there is a finite list of (not necessarily distinct) Grassmannian permutations $u_1,u_2,\dots,u_k \in S_\infty$
such that $G_w(\bfz) =  \sum_{i=1}^k \beta^{\ell(u_i)-\ell(w)} G_{u_i}(\bfz)$.
A later bijective proof of this result appears in \cite{BKSTY}.

\item It follows that the set
$  \{ G_\lambda(\bfz) : \lambda\text{ any partition}\}$ is a  $\ZZ[\beta]$-basis for the bialgebra $\bfG$
whose structure constants for both the product and coproduct are all elements of $\NN[\beta]$.
More strongly, Buch's work \cite{Buch2002} provides explicit combinatorial formulas for these coefficients.
\eei
Thus, in type $\sfA$, the set of symmetric Grothendieck functions indexed by Grassmannian permutations
provide a basis for a bialgebra with certain strong positivity properties.
Our goal now is to explain how one can extend these results to the other classical types $\sfB$, $\sfC$, and $\sfD$.

\subsection{K-Stanley symmetric functions}\label{intro2-sect}

A \defn{signed permutation} 
is a bijection
$ w:\ZZ\to\ZZ$ that has $w(-i) = -w(i)$ for all $i \in \ZZ$
and $w(i)=i$ for all but finitely many $i \in \ZZ$.
Define $\WB_\infty=\WC_\infty$ to be the group of all signed permutations 
and let $\WD_\infty$ be the subgroup of signed permutations $w$ with an even number of sign changes,
that is, for which the set $\{ i \in \PP : w(i) < 0 \}$ has even size.

Motivated by the problem of constructing equivariant $K$-theory representatives for Schubert varieties in
 flag varieties of all classical types,
Kirillov and Naruse have introduced three analogues of $G_w(\bfz)$
indexed by each of the classical Weyl groups $W^\sfB_\infty$, $W^\sfC_\infty$, and $W^\sfD_\infty$.
These power series,
which we write  as $G^\sfB_w(\bfz)$, $G^\sfC_w(\bfz)$, and $G^\sfD_w(\bfz)$,
are again elements of $\Sym_\beta$ and will be referred to as \defn{$K$-Stanley symmetric functions}.
This terminology is motivated by the fact
that setting $\beta=0$ reduces $\cFX_w(\bfz)$ to the \defn{Stanley symmetric functions} for classical types introduced in \cite{BilleyHaiman} and further studied in \cite{BL,FK1996,TKLamThesis,TKLam}.
For the precise definitions, see Section~\ref{k-stanley-sect}.

Fix a type $X \in \{\sfB,\sfC,\sfD\}$.
Then we may consider the submodule of symmetric functions
\be
\bfG^X := \ZZ[\beta]\spanning \{ G^X_w(\bfz) : w \in W^X_\infty \}.
\ee
One consequence of our main results is the following analogue of Proposition~\ref{dem-prop} for $\bfG^X$.

\begin{theorem}\label{conj1}
The modules $\bfG^\sfB$, $\bfG^\sfC$, and $\bfG^\sfD$
are sub-bialgebras of $\bfG$
with
$\bfG^\sfC\subset \bfG^\sfB=\bfG^\sfD$.

\end{theorem}

These bialgebras have a more precise characterization.
Let $\SS_\beta$ be the subring of power series $f(\bfz)\in\Sym_\beta$
that have the following \defn{$K$-supersymmetry property} introduced in \cite{IN}:
\be\label{ksym-eq}
f(t, \tfrac{-t}{1+\beta t}, z_3,z_4,z_5,\dots) = f(0, 0, z_3,z_4,z_5,\dots).
\ee
Ikeda and Naruse \cite{IN} have identified two distinguished pseudobases for $\SS_\beta$,
given by the \defn{$K$-theoretic Schur $P$- and $Q$-functions} $\GP_\lambda(\bfz)$
and $\GQ_\lambda(\bfz)$.
These are the generating functions for certain \defn{semistandard set-valued shifted tableaux}
of strict partition shape $\lambda=(\lambda_1>\lambda_2>\dots >0)$; see Section~\ref{ss-sect} for a detailed review. 
The sets of finite linear combinations
\be
\ba \GPB&:=\ZZ[\beta]\spanning\left\{\GP_\lambda(\bfz)  : \lambda\text{ any strict partition}\right\}\quand
\\
\GQB&:=\ZZ[\beta]\spanning\left\{\GQ_\lambda(\bfz) : \lambda\text{ any strict partition}\right\}
\ea
\ee
are both proper submodules of $\SS_\beta$.

The modules $\GPB$ and $\GQB$ have been studied previously in \cite{ChiuMarberg,CTY,IN,LM2,MarScr}.
It is known  that both are closed under multiplication and thus form rings \cite{LM2}, with 
 strict containment $\GQB\subsetneq \GPB\subsetneq \bfG$ \cite[\S3]{MarScr}.
Our results will show  that these subrings are actually sub-bialgebras.

A signed permutation $w $ is \defn{Grassmannian} if $w(i) < w(i+1)$ for all $i \in \PP$.
If $w$ has this property and $n\in \NN$ is maximal with $w(n) \leq 0$ then we can define two strict partitions by
\be\label{gr-eq0}
\ba 
\lambda_\sfB(w) =\lambda_\sfC(w) &:= (-w(1),-w(2),\dots,-w(n)),
\\
 \lambda_\sfD(w) &:=(-w(1)-1,-w(2)-1,\dots,-w(n)-1).
 \ea\ee
If $X \in \{\sfB,\sfC,\sfD\}$ and $w \in W^X_\infty$ is Grassmannian
with   $\lambda = \lambda_X(w)$,
 then it is known \cite[\S4.2]{KN} that
\be\label{gr-eq}
\cFX_w(\bfz) = \begin{cases} \GP_{\lambda}(\bfz)&\text{when }X \in \{\sfB,\sfD\} \\ 
 \GQ_{\lambda}(\bfz)&\text{when $X =\sfC$.}
 \end{cases}
 \ee
The following theorem confirms predictions of Kirillov--Naruse \cite[Rem.~3]{KN} and Arroyo--Hamaker--Hawkes--Pan \cite[Conj.~7.2]{AHHP}.
Here, we write $\ell^X$  for the Coxeter length function $W^X_\infty \to \NN$ (see \S\ref{infinite-sect} for an explicit formula).

\begin{theorem*}[{See Theorem~\ref{k-stanley-positivity-thm}}]
\label{conj2}
Fix a type $X \in \{\sfB,\sfC,\sfD\}$. Then for each $w \in W^X_\infty$,
 there is a finite list of (not necessarily distinct) Grassmannian elements $u_1,u_2,\dots,u_k \in W^X_\infty$
such that \[\textstyle\cFX_w(\bfz) =   \sum_{i=1}^k \beta^{\ell^X(u_i)-\ell^X(w)} \cFX_{u_i}(\bfz)
 \]
Consequently, it holds that  $\bfG^\sfB=\bfG^\sfD =\GPB$ and $\bfG^\sfC =\GQB$.
\end{theorem*}

 We prove this result  in a more precise form in Section~\ref{k-positivity-sect}.
 Our methods provide an explicit algorithm for determining the $\GP$- and $\GQ$-expansions
 of $\cFX_w(\bfz)$ (see Corollary~\ref{bcd-thm3}), but it is an open problem to find a positive combinatorial formula
 for the associated coefficients, in the style of \cite{BKSTY}. The recent paper \cite{AHHP} provides a conjectural formula of this kind for type $\sfC$.
 
Almost all of Theorem~\ref{conj1} follows as a corollary to the preceding discussion.
It only remains to explain why each subalgebra $\bfG^X\subset \bfG$ is also a sub-coalgebra.
We establish this property in the following way.

The shifted tableau generating functions for $\GP_\lambda(\bfz)$ and $\GQ_\lambda(\bfz)$ reviewed in Section~\ref{ss-sect}
have natural skew analogues $\GP_{\lambda/\mu}(\bfz)$ and $\GQ_{\lambda/\mu}(\bfz)$ indexed by pairs of strict partition with $\mu\subset\lambda$ in the sense that $\mu_i \leq \lambda_i$ for all $i$. 
These functions naturally arise when considering the coproduct of $\bfG$, as it is known \cite[Prop.~5.6]{LM} that 
\be\label{coprod-eq}
\ba
\Delta(\GP_\lambda(\bfz))&\in \NN[\beta]\spanning\{ \GP_{\lambda/\kappa}(\bfz) \otimes \GP_\mu(\bfz) : \kappa \subseteq \mu\subseteq\lambda\}\quand
\\
\Delta(\GQ_\lambda(\bfz))&\in \NN[\beta]\spanning\{ \GQ_{\lambda/\kappa}(\bfz) \otimes \GQ_\mu(\bfz) : \kappa \subseteq \mu\subseteq\lambda\}.
\ea
\ee 
The following result, which is proved in Section~\ref{ss-sect},
extends \cite[Thm.~1.5]{AHHP} from type $\sfC$ to the remaining classical
types $\sfB$ and $\sfD$. This statement implies that $\bfG^X\subset \bfG$ is a sub-coalgebra.

\begin{theorem*}[{See Theorem~\ref{fc-prop}}]\label{conj3}
For any strict partitions $\mu\subseteq\lambda$
there are signed permutations 
$u,v\in W^\sfB_\infty=W^\sfC_\infty$ and $w \in W^\sfD_\infty$ 
such that 
$
\cFB_{u}(\bfz)=\cFD_{w}(\bfz) = \GP_{\lambda/\mu}(\bfz)
$ and
$\cFC_{v}(\bfz) = \GQ_{\lambda/\mu}(\bfz)$.
\end{theorem*}

We highlight some other consequences of our main results.
Combining the preceding theorems with \eqref{gr-eq}
gives the following corollary, which was predicted as \cite[Conj.~5.14]{LM}: 

\begin{corollary}\label{coprod-cor}
There are coefficients $f^{\nu}_{\lambda\mu},g^{\nu}_{\lambda\mu} \in \NN$
indexed by triples of strict partitions
with
\[
\GP_{\lambda/\mu}  =\sum_\nu f^{\nu}_{\lambda\mu}\cdot \beta^{|\nu|+|\mu|-|\lambda|}\cdot \GP_\nu \in \GPB
\quand
\GQ_{\lambda/\mu}  =\sum_\nu g^{\nu}_{\lambda\mu}  \cdot \beta^{|\nu|+|\mu|-|\lambda|}\cdot \GQ_\nu \in \GQB.
\]
When $\mu\subseteq \lambda$ are fixed, the integers $f^{\nu}_{\lambda\mu}$ and $g^{\nu}_{\lambda\mu}$ are nonzero for only finitely many choices of $\nu$.
\end{corollary}

A conjectural combinatorial formula for the coefficients $g^\nu_{\lambda\mu}$ is given in \cite[Conj.~6.1]{AHHP}.
 
We can summarize this discussion with the following stronger statement:

\begin{corollary}\label{conj4}
The modules $\GPB$ and $\GQB$ are sub-bialgebras of $\bfG$, for which 
the respective sets $\{\GP_\lambda : \lambda\text{ strict}\}$ and $\{\GQ_\lambda : \lambda\text{ strict}\}$
are positive bases in the sense that 
\[
\ba \GP_\lambda\GP_\mu &=\sum_\nu A_{\lambda\mu}^\nu \cdot \beta^{|\nu|-|\lambda|-|\mu|}\cdot \GP_\nu
\\
 \GQ_\lambda\GQ_\mu &=\sum_\nu B_{\lambda\mu}^\nu  \cdot \beta^{|\nu|-|\lambda|-|\mu|}\cdot \GQ_\nu
 \ea
 \quand
 \ba
 \Delta(\GP_\nu) &= \sum_\lambda\sum_{\mu} a_{\lambda\mu}^\nu  \cdot \beta^{|\lambda|+|\mu|-|\nu|}\cdot \GP_\lambda\otimes \GP_\mu
 \\
 \Delta(\GQ_\nu) &=\sum_\lambda \sum_{\mu} b_{\lambda\mu}^\nu  \cdot \beta^{|\lambda|+|\mu|-|\nu|}\cdot\GQ_\lambda\otimes \GQ_\mu
 \ea
 \]
 for certain nonnegative integers $A_{\lambda\mu}^\nu,B_{\lambda\mu}^\nu,a_{\lambda\mu}^\nu,b_{\lambda\mu}^\nu \in \NN$
 indexed by triples of strict partitions.
\end{corollary}

\begin{proof}
The positivity of the multiplicative structure constants $A_{\lambda\mu}^\nu$ and $B_{\lambda\mu}^\nu$
is established in \cite[Thm.~1.6]{LM2}.
The claim that $\GQB\subsetneq\GPB\subsetneq\bfG$ is shown in \cite[\S3]{MarScr}.
The positivity of the comultiplicative structure constants is new, and follows from \eqref{coprod-eq}
and Corollary~\ref{coprod-cor}. 
\end{proof}

Combinatorial formulas are known for  the numbers $A_{\lambda\mu}^\nu$ from \cite{CTY,PechenikYong}, but
it is an open problem to find similar expressions for 
the other coefficients in Corollary~\ref{conj4}; see the conjectures in \cite[\S6]{AHHP}.
The nonnegativity of the numbers $ a_{\lambda\mu}^\nu$ and $ b_{\lambda\mu}^\nu$ confirms \cite[Conjs.~4.35 and 4.36]{MarbergHopf}.

\subsection{Double Grothendieck transitions}

Kirillov and Naruse's definition of the $K$-Stanley symmetric functions $\cFX_w(\bfz)$  is a preliminary step in their construction
of certain \defn{double Grothendieck polynomials} $\fkGX_w(\bfx;\bfy;\bfz) $ for each classical type $X \in \{\sfB,\sfC,\sfD\}$.
The objects, which belong to the ring 
\[\SS_\beta[\bfx,\bfy] := \SS_\beta[x_1,x_2,\dots,y_1,y_2,\dots],\] 
 generalize the \defn{(double) Schubert polynomials} for classical groups introduced
in \cite{BilleyHaiman,IkMiNa}.
We review the precise definition of $\fkGX_w(\bfx;\bfy;\bfz)$ in Section~\ref{knp-sect}, and just emphasize here that $\fkGX_w(0;0;\bfz) = \cFX_w(\bfz)$.

The results in \cite[\S6]{KN} show that certain specializations of
$\fkGX_w(\bfx;\bfy;\bfz) $ represent the (opposite) Schubert classes in the equivariant $K$-theory ring of 
the complete flag variety $G/B$  in type $X\in \{\sfB,\sfC,\sfD\}$.
These representatives are distinguished by several desirable stability properties.
In particular, the structure constants that expand the product
of $\fkGX_v(\bfx;\bfy;\bfz) $ and $\fkGX_w(\bfx;\bfy;\bfz) $ are the same as the ones for the Schubert basis of the equivariant $K$-theory ring $K_T(G/B)$. 

These rings have been studied previously, e.g., in  \cite{AGM,GK,GR,KK,LP,MNS,Willems}.
Most notably, Lenart and Postnikov \cite{LP} have proved an
equivariant Chevalley formula giving a multiplication rule in $K_T(G/B)$
for all types.
 By translating this rule to a polynomial identity for $\fkGX_w(\bfx;\bfy;\bfz)$,
 we are able to prove a \defn{transition formula} of the following kind:

\begin{theorem*}[{See Theorem~\ref{bcd-thm2}}]
Fix a type $X \in \{\sfB,\sfC,\sfD\}$ and a signed permutation $w \in W^X_\infty$.
Let $a \in [1,\infty)$ be maximal with $w(a) > w(a+1)$, define $b \in [a+1,\infty)$ to be maximal
with $w(a)>w(b)$, and set $v = w(a,b)(-a,-b)$ and $c= w(b)$. Then there is an explicit $ \ZZ[\beta][\bfy]$-linear operator
\[
\bfRX_a : \ZZ[\beta][\bfy]\spanning\left\{ \fkGX_u(\bfx; \bfy;\bfz)  : u \in W^X_\infty\right\} \to  
 \ZZ[\beta][\bfy]\spanning\left\{ \fkGX_u(\bfx; \bfy;\bfz)  : u \in W^X_\infty\right\}
\]
such that 
\[\fkG^X_w(\bfx; \bfy; \bfz) =   (1+ \beta y_{c})(1+ \beta x_a) \bfRX_a \fkGX_v(\bfx; \bfy;\bfz) -\fkGX_v(\bfx;\bfy;\bfz) 
\] where we set $y_{c} := \frac{-y_{|c|}}{1+ \beta y_{|c|}}$ if $c<0$.
\end{theorem*}

The form of this statement is based on an analogous result in type $\sfA$ stated in unpublished work of Lascoux \cite{Lascoux}
and proved combinatorially by Weigandt \cite[Thm.~3.1]{Weigandt}. (We provide an alternate proof in Section~\ref{a-transitions-sect}.)
While not obvious from our presentation, setting $\bfx=\bfy=0$ in the preceding theorem
will provide an $\NN[\beta]$-linear recurrence for $\cFX_w(\bfz)$ which can be used to show 
the main results in the previous section. This approach is similar to (but more complicated than)
than the methods used for the cohomological case in \cite{Billey,IkMiNa}.

Efficient methods of computing $\fkGX_w(\bfx;\bfy;\bfz)$ are somewhat limited.
For Grassmannian signed permutations, there are tableau generating functions in \cite{HIMN,KN}.
Tamvakis \cite{Tamvakis} has found some related formulas for certain \defn{skew} signed permutations. 
 Outside of these cases, it is not particularly easy to calculate or represent the power series $\fkGX_w(\bfx;\bfy;\bfz)$.
The result above (Theorem~\ref{bcd-thm2}) provides a recursive solution to this problem,
by giving an explicit finite algorithm to expand
\be\fkGX_w(\bfx;\bfy;\bfz)  \in \ZZ[\beta][\bfy]\spanning\Bigl\{ \fkGX_u(\bfx;\bfy;\bfz) : u \in W^X_{\infty}\text{ is Grassmannian}\Bigr\}.
\ee

The rest of this paper is organized as follows. Section~\ref{prelim-sect} contains some standard preliminaries on root systems and Weyl groups.
Section~\ref{ecf-sect} reviews Lenart and Postnikov's Chevalley formula for equivariant $K$-theory, and its specializations to classical type. Section~\ref{te-sect}, finally, derives transition formulas for double Grothendieck polynomials
and $K$-Stanley symmetric functions in classical type in order to prove the main results summarized in this introduction.

\subsection*{Acknowledgments}

This article is based on work supported by the National Science Foundation under grant DMS-1929284 while the author was in residence at the Institute for Computational and Experimental Research in Mathematics in Providence, RI, during the Categorification and Computation in Algebraic Combinatorics semester program.
The author was also partially supported by Hong Kong RGC grants 16304122 and 16304625.

We thank David Anderson, Joshua Arroyo, Stephen Griffeth, Zachary Hamaker, Cristian Lenart, Jianping Pan, Brendan Pawlowski, and Travis Scrimshaw for helpful discussions. 
 
\section{Preliminaries}\label{prelim-sect}

This section reviews some preliminaries on root systems and Weyl groups in classical type.
This material, which serves primarily to set up our notation, follows the conventions in \cite{Billey,IkMiNa}.

\subsection{General setup}\label{genset-sect}

Let $G$ be a connected, simply connected, simple complex Lie group. 
Although the setup that we are about to present is completely general,
our focus elsewhere will be exclusively on the cases
when $G$ is one of the classical groups $\SL_n(\CC)$, $\SO_{2n}(\CC)$, $\Sp_{2n}(\CC)$, or $\SO_{2n+1}(\CC)$.

Fix a Borel subgroup $B\subset G$ and a maximal torus $T\subset B$.
Let $\fkh$ be the corresponding Cartan subalgebra of the Lie algebra $\fkg$ of $G$. 
Let $\Phi\subset \fkh^*$ be the corresponding irreducible root system,
and write $V= \RR\spanning\{\alpha \in \Phi\}  \subset \fkh^*$ for the real span of the roots.
Let $(\cdot,\cdot) : V\times V \to \RR$ be the nondegenerate bilinear form 
induced by the Killing form on $\fkg$. 

The choice of $B$ determines a partition of $\Phi = \Phi^+\sqcup -\Phi^+$ into positive and negative roots.
The set of \defn{positive roots} $\Phi^+$ then contains a unique set of \defn{simple roots} 
$\Pi = \{\alpha_1,\alpha_2,\dots,\alpha_n\}$ which form an $\RR$-basis for $V$.
The \defn{fundamental weights} 
$\omega_1,\omega_2,\dots,\omega_n \in V$
are the elements with 
\be(\omega_i,\alpha_j^\vee) =\delta_{ij} := |\{i\}\cap\{j\}|, 
\quad\text{where $\alpha^\vee := \tfrac{2}{(\alpha, \alpha)}\alpha$ for any $\alpha \in \Phi$}.
\ee
These vectors are a $\ZZ$-basis for
the \defn{weight lattice} 
\be\Lambda = \{ \lambda \in V : (v,\alpha^\vee) \in \ZZ\text{ for all }\alpha \in \Phi\}.\ee
Finally,
the \defn{Weyl group} $W \subset \GL(V)$ of $G$ is the (finite) group generated by the reflections
\be
\barr{cccl} r_\alpha :& V & \to & V \\
& v & \mapsto & v - (v,\alpha^\vee)\alpha
\earr
\quad\text{for all }\alpha \in \Phi^+.
\ee
The smaller set $S = \{r_\alpha : \alpha \in \Pi\}$ also generates $W$, and the 
pair $(W,S)$ is a Coxeter system
with length function $
\ell ( w )= | \{\alpha \in \Phi^+ : w(\alpha) \notin \Phi^+\}|.
$
We write $w_0 \in W$ for the unique
\defn{longest element} in the Weyl group, which has length $\ell(w_0) = |\Phi^+|$.

\subsection{Notation for classical groups}

We fix some notation to describe the primary cases of interest when $G$ is a classical group.
Throughout, $n$ denotes a positive integer and $[n] =\{1,2,\dots,n\}$. We write $\binom{[n]}{2}$ to denote the set of 2-element subsets $\{i,j\}\subset[n]$ with $i\neq j$.
Write $\e_1,\e_2,\dots,\e_n$ for the standard basis of $\RR^n$
and let $(\cdot,\cdot)$ be the  usual bilinear form on $\RR^n$ with  $(\e_i,\e_j) = \delta_{ij}$.

\subsubsection{Type A}\label{a-sect1}
In type $\sfA_{n-1}$ when $G= \SL_{n}(\CC)$ 
we choose $T\subset B\subset G$ and $\fkh\subset \fkg$
such that the positive roots $\Phi^+$ and simple roots $\Pi$ are
\be 
\ba
\Phi^+_{\sfA_{n-1}}&= \Bigl\{ \e_j-\e_i : 1\leq i<j\leq n\Bigr\}
\quand\\
\Pi_{\sfA_{n-1}} &= \Bigl\{ \alpha_i = \e_{i+1}-\e_{i} : i \in [n-1]\Bigr\},
\ea
\ee
viewed as subsets of the quotient vector space $V = \RR^n / \RR\spanning\{\e_1+\e_2+\dots+\e_n\}$.
The bilinear form on $V$ is defined by setting
$(\alpha_i, \alpha_i)=2$, $(\alpha_i, \alpha_{i\pm 1}) = -1$, and $(\alpha_i, \alpha_j)=0$ when $|i-j|>1$,

so that $\alpha^\vee = \alpha$ for all roots $\alpha \in \Phi$.
The fundamental weights   are 
\be \omega_i = \e_{i+1}+\e_{i+2}+\dots +\e_n\quad\text{for }1\leq i \leq n-1.
\ee
Notice that 
$-\omega_1 = \e_1 \in V$ while $\omega_{i}-\omega_{i+1} = \e_{i+1}$ for $i \in [n-2]$ and $\omega_{n-1}=\e_n$.

Letting $s_i = r_{\alpha_i}$ act on $\ZZ$ as the transposition $(i,i+1)$ determines an isomorphism $W \cong S_n$, where we define the symmetric group $S_n$  to be the set of permutations of $\ZZ$ that fix all points outside $[n]$.
The longest element of  $S_n=\langle s_1,s_2,\dots,s_{n-1}\rangle$ is the reverse permutation
 $ w_0 = n\cdots 321$.

\subsubsection{Type B}\label{b-sect1}

For the other three classical groups the ambient vector space is $V=\RR^n$ equipped with the standard bilinear form $(\cdot,\cdot)$ satisfying $(\e_i,\e_j) = \delta_{ij}$.
In type $\sfB_n$ when $G=\SO_{2n}(\CC)$ 
we choose $T\subset B\subset G$ and $\fkh\subset \fkg$
such that the positive roots $\Phi^+$ and simple roots $\Pi$ are the sets
\be 
\ba
\Phi^+_{\sfB_{n}}&= \Bigl\{\e_1,\e_2,\dots,\e_n\Bigr\} \sqcup \Bigl\{ \e_j\pm\e_i : 1\leq i<j\leq n\Bigr\}
\quand \\
\Pi_{\sfB_{n}} &= \Bigl\{ \alpha_0 = \e_1\Bigl\}\sqcup \Bigl\{ \alpha_i = \e_{i+1}-\e_{i} : i \in [n-1]\Bigl\}.
\ea
\ee
The fundamental weights are 
\be
\ba
\omega_0 &= \tfrac{1}{2}(\e_1+\e_2+\dots+\e_n) \text{ and }
\\
\omega_i &= \e_{i+1}+\e_{i+2}+\dots+\e_n\text{ for }1\leq i \leq n-1.
\ea
\ee

Letting $t_0 = r_{\alpha_0}$ act on $\ZZ$ as the transposition $(-1,1)$
while letting $t_i = r_{\alpha_i}$ for $i>0$ act as $(i,i+1)(-i,-i-1)$
determines an 
isomorphism $W \cong \WB_n$, where $\WB_n$ is the group of bijections $w:\ZZ\to \ZZ$
that fix all numbers $i>n$ and satisfy $w(-i) =-w(i)$ for all $i \in \ZZ$.

The one-line representation of  $w\in \WB_n$  is the sequence $w(1)w(2)\cdots w(n)$. When writing elements in this way, we use the notation $\overline a$ in place of $-a$. For example, we have 
$\WB_2 = \{12,\ 21,\ \overline 12,\ 2\overline 1,\ 1\overline2,\ \overline21,\ \overline1\overline2,\ \overline2\overline1\}.$
The longest element of $\WB_n=\langle t_0,t_1,\dots,t_{n-1}\rangle$
 is  $ w_0 = \overline1\overline2\overline3\cdots\overline n$.

\subsubsection{Type C} \label{c-sect1}

In type $\sfC_n$ when $G=\Sp_{2n}(\CC)$ 
we choose $T\subset B\subset G$ and $\fkh\subset \fkg$
such that $\Phi^+$ and $\Pi$ are 
\be 
\ba
\Phi^+_{\sfC_{n}}&= \Bigl\{2\e_1,2\e_2,\dots,2\e_n\Bigr\} \sqcup \Bigl\{ \e_j\pm\e_i : 1\leq i<j\leq n\Bigr\}
\quand \\
\Pi_{\sfC_{n}} &= \Bigl\{ \alpha_0 = 2\e_1\Bigl\}\sqcup \Bigl\{ \alpha_i = \e_{i+1}-\e_{i} : i \in [n-1]\Bigl\}.
\ea
\ee
The fundamental weights are 
\be
\ba
\omega_i &= \e_{i+1}+\e_{i+2}+\dots+\e_n\text{ for }0\leq i \leq n-1
\ea
\ee
The Weyl group $W $ has the same description as in type $\sfB_n$. For convenience, let 
$\WC_n :=\WB_n.$
 
 \subsubsection{Type D}\label{d-sect1}
 
 In type $\sfD_n$ when $G=\SO_{2n+1}(\CC)$ 
we choose $T\subset B\subset G$ and $\fkh\subset \fkg$
such that $\Phi^+$ and $\Pi$ are  
\be 
\ba
\Phi^+_{\sfD_{n}}&= \Bigl\{ \e_j\pm\e_i : 1\leq i<j\leq n\Bigr\}\quand
\\
\Pi_{\sfD_{n}} &= \Bigl\{ \alpha_{-1} =  \e_2+\e_1\Bigl\}\sqcup \Bigl\{ \alpha_i = \e_{i+1}-\e_{i} : i\in[n-1]\Bigl\}.
\ea
\ee
The fundamental weights are 
\be
\ba
\omega_{-1} &=\tfrac{1}{2}(\e_1+\e_2+\dots+\e_n),\\
\omega_1 &= \tfrac{1}{2}(-\e_1+\e_2+\dots+\e_n), \text{ and }
\\
\omega_i &= \e_{i+1}+\e_{i+2}+\dots+\e_n\text{ for }2\leq i \leq n-1.
\ea
\ee

Letting $t_{-1} = r_{\alpha_{-1}}$ act on $\ZZ$ as the permutation $(1,-2)(2,-1)$
while letting $t_i = r_{\alpha_i}$ for $i>0$ act as $(i,i+1)(-i,-i-1)$
determines an 
isomorphism $W \cong \WD_n$, where $\WD_n$ is the subgroup of elements $w \in \WB_n$ such that
the number of integers $i \in [n]$ with $w(i)<0$ is even.
The longest element of the Coxeter group $\WD_n=\langle t_{-1},t_1,t_2,\dots,t_{n-1}\rangle$ is  
\be
 w_0 =\begin{cases}
\overline1\overline2\overline3\cdots\overline n&\text{when $n$ is even}, \\
 1\overline2\overline3\cdots\overline n&\text{when $n$ is odd}.
\end{cases}
 \ee

\subsection{Infinite rank Weyl groups}\label{infinite-sect}
Our conventions make it automatically true that $S_n\subset S_{n+1}$,  $\WB_n\subset \WB_{n+1}$, and $ \WD_n\subset \WD_{n+1}$ for all positive integers $n$.
The infinite unions 
\be\textstyle
S_\infty := \bigcup_{n\geq 1} S_n
\quand
\WB_\infty:= \bigcup_{n\geq 1} \WB_n
\quand
\WD_\infty := \bigcup_{n\geq 1} \WD_n
\ee
are also Coxeter groups,
which respectively consist
of all finitely supported permutations of the positive integers, all finitely supported signed permutations of $\ZZ$,
and all finitely supported even signed permutations.
For convenience we also set
$
\WA_\infty = \langle t_1,t_2,t_3,\dots\rangle$ and $
\WC_\infty := \WB_\infty.
$

\begin{remark}
Notice that  $S_\infty = \langle s_1,s_2,s_3,\dots\rangle$ is technically distinct from the subgroup 
\[
\WA_\infty = \langle t_1,t_2,t_3,\dots\rangle\subsetneq \WD_\infty = \langle t_{-1},t_1,t_2,t_3,\dots\rangle \subsetneq \WB_\infty = \WC_\infty= \langle t_0,t_1,t_2,t_3,\dots\rangle.
\] However, the map sending $s_i \mapsto t_i$ determines an isomorphism $S_\infty \to \WA_\infty$ of Coxeter groups.
We sometimes identify the elements of the two groups via this bijection.
\end{remark}

The Coxeter length functions of $S_\infty$, $\WB_\infty = \WC_\infty$, and $\WD_\infty$ 
restrict to the corresponding length functions of   $S_n$, $\WB_n=\WC_n$, and $\WD_n$ for each $n\in\PP$.
These maps have the formulas \cite[\S3]{Billey}
\be
\ellA(w)  := \inv(w),\quad
\ellB(w)= \ellC(w) := \tfrac{\inv(w) +\ell_0(w)}{2},
\quand
\ellD(w) :=  \tfrac{\inv(w) -\ell_0(w)}{2} 
\ee
where for any finitely supported permutation $w$ of $\ZZ$ we define
\be
\ba
\inv(w) &:= |\{ (i,j) \in \ZZ\times \ZZ : i<j\text{ and }w(i)>w(j)\}|\quand
\\
\ell_0(w) &:= |\{ i \in \PP : w(i) < 0\}|.
\ea
\ee

The \defn{reflections} in a Coxeter group are the elements conjugate to the simple generators.
In $S_\infty$ these elements are the transpositions $(i,j)$ for $1\leq i < j$.

To denote the reflections in $\WB_\infty=\WC_\infty$ and $\WD_\infty$,
we define 
\be\label{tij-eq}
t_{ij} := (i,j)(-j,-i) \in   \WD_\infty  
\quand
t_{0k} := (-k,k) \in  \WB_\infty=  \WC_\infty
\ee
for all integers $i,j,k \in \ZZ$ with $k \neq 0$.
In the notation of \cite{Billey}, we have \be t_{0k} = s_{kk}\quand t_{-i,j} = s_{ij}
\quad\text{for }i,j,k\in [n]\text{ with }i\neq j
.\ee  
Notice that
$
t_{jj}=t_{-j,j} = 1
$ and that $t_{ij}=t_{ji}=t_{-i,-j}=t_{-j,-i} $.
The distinct reflections in $\WB_\infty=\WC_\infty$ consist of the elements $t_{ij}$ with $|i|<j$.
The reflections in $\WD_\infty$ consist of the same set but with $t_{0j}$ excluded for all $j \in \PP$.

We conclude this section by reviewing 
the properties  \cite[Lemmas~1 and 2]{Billey}
that characterize when multiplication by a reflection
increases the length of a classical Weyl group element by exactly one.

\begin{lemma}[{\cite[Lem.~1]{Billey}}] \label{len-lem1}
Fix $X \in \{\sfA,\sfB,\sfC,\sfD\}$ and write $\ell=\ell^X$. 
Suppose $t=t_{ij}=(i,j)(-i,-j)$
where $0<i<j$.
Then $w \in W^X_\infty$ has $\ell(wt)=\ell(w)+1$ if and only if 
\[w(i) < w(j)
\quand
\text{no $e\in\ZZ$ exists with $i<e<j$ and $w(i)<w(e)<w(j)$.}\]
\end{lemma}

\begin{lemma}[{\cite[Lem.~2]{Billey}}]  \label{len-lem2}
Fix $X\in \{\sfB,\sfC\}$ and write $\ell=\ell^X$. Suppose $t=t_{0j} = (-j,j)$ where $0<j$.
Then $w \in W^X_\infty$ has $\ell(wt)=\ell(w)+1$ if and only if 
\[
0<w(j)
\quand
\text{no $e \in \ZZ$ exists with $0<e < j$ and $-w(j) < w(e) < w(j)$.}
\]
\end{lemma}

\begin{lemma}[{\cite[Lem.~2]{Billey}}]  \label{len-lem3}
Fix $X \in \{\sfB,\sfC,\sfD\}$ and write $\ell=\ell^X$.  Let $t=t_{-i,j}=(-i,j)(i,-j)$ where $0<i <j$.
Then $w \in W^X_\infty$ has $\ell(wt)=\ell(w)+1$ if and only if the following conditions hold:
\ben
\item[(a)] we have $-w(i) < w(j)$; 
\item[(b)] when $X=\sfB$, either $w(i)<0$ or $w(j) < 0$;
\item[(c)] no $e \in \ZZ$ exists with $0<e < i$ and $-w(j) < w(e) < w(i)$; and
\item[(d)] no $e \in \ZZ$ exists with $0<e < j$ and $-w(i) < w(e) < w(j)$.
\een
\end{lemma}

Denote the \defn{descent set} of  $w \in \WB_\infty$ by
$\Des(w) =  \{ i \in \PP : w(i) > w(i+1)\}$. 
Then, following \cite[\S5]{Billey}, we define the \defn{least descent} of $w \in \WB_\infty$ to be 
the nonnegative integer
\be\label{LD-eq}
\LD(w) = \max\( \{0\}\sqcup \Des(w)\) \in \NN.\ee
The following technical consequence of Lemmas~\ref{len-lem1}, \ref{len-lem2}, and \ref{len-lem3} will be used in  Section~\ref{k-positivity-sect}.

\begin{lemma}\label{tech-lem}
Fix $X \in \{\sfB,\sfC,\sfD\}$ and $n \in \PP$.
Suppose $w \in W^X_n$. Choose integers  $i\in\ZZ$ and $k\in[n]$ with $i<k$ such that $t:=t_{ik} \in W^X_\infty$ and $\ell^X(wt) = \ell^X(w)+1$.
Then the following holds:
\ben
\item[(a)] $wt \in W^X_{n+1}$;
\item[(b)] $wt(k) < w(k)$; and
\item[(c)] $\Des(wt) \subseteq \Des(w)\cup \{k-1\}$.
\een
Moreover, if we also have $wt \notin W^X_n$ then these additional properties hold:
\ben
\item[(d)]  $i=-n-1$; and
\item[(e)] $\ell^X(wt\cdot t_{jk}) \neq \ell^X(wt)+1$ for all integers $i<j<k$ with $t_{jk} \in W^X_\infty$.
\een
\end{lemma}

\begin{proof}
If $-n\leq i < k$ then $t \in W^X_{n}$ so $wt \in W^X_{n}$, 
and if $i < -n-1$ then $\ell(wt) \neq \ell(w)+1$ by Lemma~\ref{len-lem3}. 
Only in the remaining case when $i=-n-1$ and $t \in W^X_{n+1}$ can we have $wt \in W^X_{n+1}$.
This establishes part (a),
and also shows that if $wt \notin W^X_n$ then $i=-n-1$.
In this case  $wt(k) = -n-1$ and $-n \leq wt(j) \leq n$ for all $i<j<k$, so 
part (e) is evident from Lemmas~\ref{len-lem1}, \ref{len-lem2}, and \ref{len-lem3}.

Parts (b) is also clear from the same lemmas.
It remains only to check part (c).
If $0<i<k$, then Lemma~\ref{len-lem1} 
implies that
$\Des(wt)$ is formed from $\Des(w)$ by simultaneously removing $i-1$ when both $i>1$ and $w(k) > w(i-1) > w(i)$,
 removing $k$ when $w(k) > w(k+1) > w(i)$, and adding $k-1$ in the event that $i=k-1$.
The desired containment $\Des(wt) \subseteq \Des(w)\cup \{k-1\}$ is therefore clear in this case.
When $i<0<k$ one can check that $\Des(wt) \subseteq \Des(w)$; see \cite[Lem.~3]{Billey}.
\end{proof}

\section{Equivariant Chevalley formulas}\label{ecf-sect}

This section reviews a multiplication rule for the equivariant $K$-theory ring of $G/B$
due to Lenart and Postnikov \cite{LP}. We first present this rule in its original version,
then reformulated for the basis of opposite Schubert classes, and then
in a more compact form specialized to classical type.
The main results here are Theorems~\ref{lp-thm1} and \ref{lp-thm2}
and Corollaries~\ref{lp-cor1} and \ref{lp-cor2}.

\subsection{Equivariant K-theory rings}

We return to the general setup from Section~\ref{genset-sect}.
Below, we review from \cite[\S3]{LP} the standard presentation of the equivariant $K$-theory ring $K_T(G/B)$
of the complete flag variety $G/B$.

Let $\ZZ[\Lambda]$ and $\ZZ[X]$ be the free $\ZZ$-modules with respective bases given by the 
formal symbols $e^\lambda$ and $x^\lambda$ for $\lambda \in\Lambda$.
These abelian groups are rings
with multiplication rules
\[e^\lambda  e^\mu = e^{\lambda+\mu}
\quand
x^\lambda  x^\mu = x^{\lambda+\mu}.\]
We view $\ZZ[\Lambda]$ as the group algebra of the weight lattice 
and $\ZZ[X]$ as the representation ring of the 
group 
of
characters of $T$.
Both
$\ZZ[\Lambda] = \ZZ\left[e^{\pm \omega_1},\dots,e^{\pm \omega_n}\right]$
and
$\ZZ[X] = \ZZ\left[x^{\pm \omega_1},\dots,x^{\pm \omega_n}\right]$
are isomorphic to Laurent polynomial rings in $n$ variables.

The \defn{equivariant $K$-theory ring} $K_T(G/B)$ 
is the Grothendieck ring of coherent $T$-equivariant sheaves on $G/B$. 
This ring is naturally a $\ZZ[X]$-module, which turns out to be free \cite{KK}.
Some distinguished $\ZZ[X]$-bases will be described below.

The ring $K_T(G/B)$ is isomorphic to the following quotient of 
$\ZZ[X]\otimes\ZZ[\Lambda] $.
The Weyl group $W$ acts on  $\ZZ[\Lambda]$ by $w(e^\lambda) := e^{w(\lambda)}$. 
Let $\ZZ[\Lambda]^W$ be the subalgebra of $W$-invariant elements. 
Write $\iota : \ZZ[\Lambda] \to \ZZ[X]$ for the isomorphism sending 
$e^\lambda \mapsto  x^{\lambda}$
and consider the ideal  
\be I := \left\langle \iota (f)\otimes 1-1\otimes f : f \in \ZZ[\Lambda]^W\right\rangle
\subset \ZZ[\Lambda] \otimes \ZZ[X] .\ee
For each $\lambda \in \Lambda$ let $\CC_{-\lambda}$ be the one-dimensional complex $B$-representation whose character restricts to $-\lambda$ on $T$, 
define $\cL_\lambda := G\times _\sfB \CC_{-\lambda} $ to be the line bundle over $G/B$ associated with the weight $\lambda \in \Lambda$, and write $[\cL_\lambda]$ for its class in $K_T(G/B)$.
Then there is a unique $\ZZ[X]$-algebra isomorphism \be K_T(G/B) \xrightarrow{\sim} (\ZZ[X] \otimes \ZZ[\Lambda])/I\ee
sending
$
[\cL_\lambda] \mapsto e^{-\lambda}.$
We identify $K_T(G/B)$ with its image under this map, so that from now on
\be
[\cL_\lambda] = e^{-\lambda}\in K_T(G/B)\quad\text{for all }\lambda \in\Lambda.\ee

The flag variety $G/B$ is a smooth projective variety that
 decomposes into a disjoint union of \defn{Schubert cells} 
 $BwB/B$ indexed by elements of the Weyl group $w \in W$. 
 The \defn{Schubert variety} 
 \be X_w := \overline{BwB/B}\quad\text{for $w \in W$}\ee is the closure of the corresponding Schubert cell.
In this notation, the subvariety $X_w$ has dimension $\ell(w)$ and codimension $\ell(w_0) - \ell(w)$, and we have $X_u \subseteq X_v$ 
if and only if $u \leq v$ in the \defn{(strong) Bruhat order} on $W$ 
(see \cite[\S2]{LP}).

For each $w \in W$ let
$\cO_{X_w}$ be the structure sheaf of $X_w$,
and write $[\cO_{X_w}]$ for the corresponding class in $K_T(G/B)$.
These \defn{Schubert classes}
are a basis for $K_T(G/B)$ as a free $\ZZ[X]$-module \cite[\S4]{KK}. Since $X_{w_0} = G/B$ has codimension zero, one has 
\be [\cO_{X_{w_0}}] = 1 \in K_T(G/B).\ee
Lenart and Postnikov \cite[Lem.~8.1]{LP} (see also \cite[Eq.~(0.3)]{GR}) derive a simple formula for $[\cO_{X_w}]$ in the codimension one case:
if $\alpha_k \in \Pi$ is a simple root and $s_k= r_{\alpha_k} \in W$, then
 \be \label{codim1-eq}
[\cO_{X_{w_0 s_k}}] = 1 - x^{-w_0(\omega _k)} [\cL_{-\omega_k}]= 1 - x^{-w_0(\omega _k)} e^{\omega_k} \in K_T(G/B)
\ee
where $\omega_k\in\Lambda$ is the unique fundamental weight satisfying $( \omega_k, \alpha_k^\vee)=1$.

\begin{remark}\label{notation-rem}
There are some potentially confusing differences in notation in the literature \cite{AGM,GK,GR,KK,LP,MNS,Willems} studying $K_T(G/B)$.
The following table attempts to clarify this situation:
\begin{center}
\begin{tabular}{l | l| l | l }
symbol here and in \cite{LP} & symbol in \cite{GR} & symbol in \cite{AGM,GK}  & symbol in \cite{MNS} 
 \\ \hline  && & \\[-8pt]
 $\cL_\lambda$ & $X^\lambda$ & $\cL(\lambda)$ & $\cL_{-\lambda}$
 \\&&& \\[-8pt]
$e^\lambda =[\cL_{-\lambda}]$ & $[X^{-\lambda}]$ & $[\cL(-\lambda)]$  & $[\cL_\lambda]$ 
\\&&& \\[-8pt]
$x^\lambda$ & $e^{-\lambda}$ & $e^{\lambda}$  & $e^\lambda$
\end{tabular}
\end{center}
Finally, there is a ring involution $\ast : K_T(G/B) \to K_T(G/B)$ sending $x^\lambda e^\mu \mapsto x^{-\lambda} e^{-\mu}$,
and most results in \cite{LP} are stated in terms of the modified Schubert classes $[\cO_w] := \ast [\cO_{X_w}]$.
\end{remark}

\subsection{Chevalley formula of Lenart--Postnikov}

Everywhere in this section we fix an arbitrary weight $\lambda \in \Lambda$. 
Here, we review the \defn{Chevalley formula}
proved by Lenart and Postnikov to expand  $[\cL_\lambda]\cdot[\cO_{X_w}]$ in the Schubert basis of $K_T(G/B)$.

\begin{definition}[{See \cite[\S6]{LP}}]
\label{cS-def}
Define $\cS_\lambda$ to be the set of pairs $(\alpha,k) \in \Phi^+\times \ZZ$ 
 such that 
 \be(\lambda,\alpha^\vee) \neq 0
\quand
 \begin{cases}
0\leq k < (\lambda,\alpha^\vee) & \text{if }(\lambda,\alpha^\vee) >0, \\
0 > k \geq (\lambda,\alpha^\vee) &\text{if }(\lambda,\alpha^\vee) <0.
\end{cases}
\ee
Now, let $<$ be the total order on $\cS_\lambda$
with $(\alpha,k) < (\alpha',k')$ if and only if 
we have $h((\alpha,k)) < h((\alpha',k'))$ in the lexicographic order on $\RR^{n+1}$,
where
 $h : \cS_\lambda \to \RR^{n+1}$ is the map sending
\be
(\alpha,k) \mapsto \tfrac{1}{(\lambda,\alpha^\vee)} \Bigl[\barr{ccccc}
k &
(\omega_1,\alpha^\vee) &
(\omega_2,\alpha^\vee) &
\cdots &
(\omega_n,\alpha^\vee)
\earr\Bigr]^\top.
\ee
Finally, suppose that with respect to this order
$\cS_\lambda = \{ (\alpha_1,k_1) < (\alpha_2,k_2) < \cdots < (\alpha_m,k_m)\}.
$
\end{definition}

Given $\beta = (\alpha,k) \in \cS_\lambda $, 
let $\raff_\beta : V\to V$ denote the \defn{affine reflection} 
\be
\raff_\beta(v) = r_\alpha(v) + k\alpha = v - (v,\alpha^\vee)\alpha + k\alpha.
\ee
For each set $J=\{j_1<\dots<j_l\}\subseteq [m]$ let 
\be n(J)  = |\{ j \in J : k_j < 0\}|.
\ee
Finally, define
$
 r_J^{[p]} =  r_{\alpha_{j_1}}\cdots r_{\alpha_{j_p}}$ for  $p \in [l]$
along with
$ r_J =   r_{\alpha_{j_1}}\cdots r_{\alpha_{j_l}}$
and
$\raff_J = \raff_{\beta_{j_1}} \cdots \raff_{\beta_{j_l}}.$

\begin{theorem}[{Lenart--Postnikov \cite[Thm.~6.1 and Prop.~6.7]{LP}}]
\label{lp-thm1}
If $u \in W$ then
\[
[\cL_\lambda] \cdot [\cO_{X_{w_0u}}] = \sum_{w \in W}\sum_{\mu \in \Lambda}  b_{u,w}^{\lambda,\mu} \cdot x^{-w_0(\mu)} \cdot  [\cO_{X_{w_0w}}]
\]
where $b_{u,w}^{\lambda,\mu} = \sum_{J} (-1)^{n(J)} $
is the sum over all subsets  $J \subseteq [m]$ with the following properties:
\be
\mu =   u \cdot  \raff_J(\lambda),
\quad
w=u\cdot r_J,
\quand \ell(u\cdot r_{J}^{[p]}) = \ell(u)+p\text{ for all $1\leq p \leq |J|$.}
\ee

\end{theorem}

\begin{remark}
To see the equivalence between Theorem~\ref{lp-thm1} and the results in \cite[\S6]{LP},
notice that 
\ben
\item[(1)] our set $\cS_\lambda$ is related to the set $\cR_\lambda$ in \cite[Prop.~6.7]{LP}
by 
$\cS_{\lambda} = \{ (\alpha, -k) : (\alpha,k) \in \cR_\lambda\}$;

\item[(2)] our coefficients $b_{u,w}^{\lambda,\mu}$ are related to the numbers $c_{u,w}^{\lambda,\mu}$ in 
\cite[Thm.~6.1]{LP} by $b_{u,w}^{\lambda,\mu} =c_{w_0u,w_0w}^{\lambda,w_0\mu}$.
\een
\end{remark}
Lenart and Postnikov provide several examples illustrating Theorem~\ref{lp-thm1}
in \cite[\S15 and \S16]{LP}.
We will not reproduce these here since our focus will be on the reformulation of Theorem~\ref{lp-thm1} for \defn{opposite Schubert classes} explained in the next section.

\subsection{Reformulation for opposite Schubert classes}\label{oppo-sect}

The \defn{opposite Schubert variety} of $v \in W$ is defined to be  
\be\label{opp-def} X^v = w_0 X_{w_0v}\subseteq G/B\ee
where $w_0 \in W$ continues to denote the longest element of the Weyl group.
Equivalently, this is the closure of the orbit of the opposite Borel subgroup $B_-=w_0Bw_0$ acting on $G/B$.
Each variety $X^v$ has dimension $\ell(w_0)-\ell(v)$ and codimension $\ell(v)$.

The classes 
$
[\cO_{X^w}] \in K_T(G/B)
$
of the structure sheaves of the opposite Schubert varieties provide another $\ZZ[X]$-basis
for the equivariant $K$-theory ring of $G/B$. 
To make our notation less cumbersome, we define
\be
\cK{w} := [\cO_{X^w}] \in K_T(G/B)\quad\text{for each }w \in W.
\ee
Because $X^1 = w_0 X_{w_0}=w_0G/B=G/B$, we have
\be
\cK{1} = 1\in K_T(G/B).
\ee
Many other ``opposite'' formulas
of interest can be derived
using the following result from \cite{MNS}, which generalizes earlier observations like \cite[Rem.~3.11]{GK} and \cite[Lem.~7.5]{GK}.

\begin{lemma}[{See \cite[\S5.2 and \S9.2.2]{MNS}}] \label{wL-lem}
The ring involution $w_0^L: K_T(G/B) \to K_T(G/B)$ sending
$x^\lambda e^\mu \mapsto x^{w_0(\lambda)} e^\mu$ for $\lambda,\mu \in \Lambda$
acts on Schubert classes as $
w_0^L([\cO_{w_0v}]) =\cK{v}$ for all $v \in W$.
\end{lemma}

Applying $w_0^L$ to both sides of \eqref{codim1-eq} gives the following:

\begin{corollary}\label{cKs-cor}
If $\alpha_k \in \Pi$ is a simple root and $s_k= r_{\alpha_k} \in S$, then 
\[
\cK{s_k} = 1 - x^{-\omega _k} [\cL_{-\omega_k}]= 1 - x^{-\omega _k} e^{\omega_k}
 \in K_T(G/B)
\]
where $\omega_k\in\Lambda$ is again the unique fundamental weight satisfying $( \omega_k, \alpha_k^\vee)=1$.
\end{corollary}

 When $c$ is a permutation or signed permutation and $G$ is the classical group of type $\sfA_{n-1}$, $\sfB_n$, $\sfC_n$, or $\sfD_n$, respectively, we write $\cKA{v}$, $\cKB{v}$, $\cKC{v}$, and $\cKD{v}$ in place of $\cK{v}$.

 \begin{example}
 \label{cKs-cor-eq}
  We can use Corollary~\ref{cKs-cor} to compute $\Theta^X_v$ when $\ell(v)=1$.
\ben
\item[(a)] In type $\sfA_{n-1}$ when $G=\SL_n(\CC)$, 
 for each $k \in [n-1]$ we have
\[
\ba 1-\cKA{s_k} 
= x^{-\omega_k}  e^{\omega_k}
&=x^{-\e_{k+1}  -\dots -\e_n}  e^{\e_{k+1}  + \dots +\e_n}
\\&= x^{\e_1 +\e_2+ \dots + \e_{k}}  e^{-\e_{1}-\e_{2} -\dots-\e_k} 
=\textstyle\prod_{i \in [k]} x^{\e_{i}}  e^{-\e_i}.
\ea
\]

\item[(b)]
In type $\sfB_n$ when $G=\SO_{2n}(\CC)$,
 for each $k \in [n-1]$ we have
\[
\ba
(1-\cKB{t_0})^2 &=  x^{-2\omega_0} e^{2\omega_0}&&=\textstyle \prod_{i =1}^n x^{-\e_i}  e^{\e_i},
\\
1-\cKB{t_k} &=  x^{-\omega_k} e^{\omega_k} &&= \textstyle \prod_{i=k+1}^n x^{-\e_i}  e^{\e_i}.
\ea
\]

\item[(c)]
In type $\sfC_n$ when $G=\Sp_{2n}(\CC)$, 
for each $k \in \{0,1,\dots,n-1\}$ we have
\[1-\cKC{t_k} =  x^{-\omega_k} e^{\omega_k} = \textstyle \prod_{i=k+1}^n x^{-\e_i} e^{\e_i}.
\]

\item[(d)] 
In type $\sfD_n$ when $G=\SO_{2n+1}(\CC)$, for each $k \in \{2,3,\dots,n-1\}$ we have
\[
\ba
(1-\cKD{t_{1}})(1-\cKD{t_{-1}}) &=  x^{-\omega_{1} - \omega_{-1}}  e^{\omega_{1} + \omega_{-1}}
&&=\textstyle \prod_{i =2}^n x^{-\e_i}  e^{-\e_i},
\\
1-\cKD{t_k} &=  x^{-\omega_k} e^{\omega_k}&&= \textstyle \prod_{i=k+1}^n x^{-\e_i}  e^{\e_i}.
\ea
\]
\een
\end{example} 

Using these examples, one derives the following proposition by a straightforward calculation.

\begin{proposition}
  \label{cL-prop}
The following identities hold in $K_T(G/B)$ when $G$ is of classical type:
\ben
\item[(a)] If $G=\SL_n(\CC)$ then $[\cL_{\e_1}] = x^{-\e_1} (1-\cKA{s_1})$
and 
$[\cL_{\e_k}] = x^{-\e_{k}} \tfrac{1-\cKA{s_k}}{1-\cKA{s_{k-1}}}$ for $2 \leq k<n$.

\item[(b)] If $G=\SO_{2n}(\CC)$ then $[\cL_{\e_1}] = x^{-\e_1} \tfrac{1-\cKB{t_1}}{(1-\cKB{t_0})^2}
$ and $[\cL_{\e_k}] = x^{-\e_{k}} \tfrac{1-\cKB{t_k}}{1-\cKB{t_{k-1}}}$ for $2\leq k<n$.

\item[(c)] If $G=\Sp_{2n}(\CC)$ then $[\cL_{\e_k}] = x^{-\e_{k}} \tfrac{1-\cKC{t_k}}{1-\cKC{t_{k-1}}}$ for $1\leq k<n$.

\item[(d)] If $G=\SO_{2n+1}(\CC)$ then $[\cL_{\e_k}] = x^{-\e_{k}} \tfrac{1-\cKD{t_k}}{1-\cKD{t_{k-1}}}$ for $3\leq k<n$, along with
\[  [\cL_{\e_1}] = x^{-\e_1}  \tfrac{1-\cKD{t_1}}{1-\cKD{t_{-1}}}
\quand
[\cL_{\e_2}] = x^{-\e_2}  \tfrac{1-\cKD{t_2}}{(1-\cKD{t_1})(1-\cKD{t_{-1}})}
.\] 
\een
\end{proposition}

 Similarly,
applying $w_0^L$ to Theorem~\ref{lp-thm1} gives this ``opposite'' Chevalley formula.

\begin{theorem}
\label{lp-thm2}
If $\lambda \in \Lambda$ and $u \in W$ then
$
[\cL_\lambda] \cdot \cK{u} = \sum_{w \in W}\sum_{\mu \in \Lambda}  b_{u,w}^{\lambda,\mu} \cdot x^{-\mu} \cdot  \cK{w}
$
where $b_{u,w}^{\lambda,\mu}\in\ZZ$ is defined as in Theorem~\ref{lp-thm1}.
\end{theorem}

\subsection{Classical specializations}

When $G$ is a classical group and $\lambda=\e_k$ is a standard basis vector,
we can express Theorem~\ref{lp-thm2} in a more explicit form, which we present in this section
along with some relevant examples.

First assume $G=\SL_n(\CC)$ and define $\ZZ$-linear operators 
 $K_T(G/B) \to K_T(G/B)$ by the formulas 
  \be
  \ba
  \fktA_{ij} &: x^\mu  \cdot \cKA{w} \mapsto \begin{cases}
x^{\mu} \cdot\cKA{w(i,j)}&\text{if }\ell(w(i,j))=\ell(w) + 1\\
0 & \text{otherwise}
\end{cases}
\\[-10pt]\\
  \fkuA_{ij} &: x^\mu  \cdot \cKA{w} \mapsto \begin{cases}
x^{w(i,j) w^{-1}(\mu)} \cdot\cKA{w(i,j)} &\text{if }\ell(w(i,j))=\ell(w) + 1\\
0 & \text{otherwise}
\end{cases}
\\[-10pt]\\
 \fkvA_{k}  &: x^\mu \cdot  \cKA{w} \mapsto x^{\mu - w(\e_k)} \cdot \cKA{w}
\ea
\ee
for all integers $1\leq i < j \leq n$ and $k\in[n]$. Then let
\be
\fkMA_{n,k}=(1 - \fktA_{k,k+1})  (1 - \fktA_{k,k+2})\cdots (1 - \fktA_{k,n}) (1 + \fkuA_{1,k}) \cdots (1 + \fkuA_{k-2,k})(1 + \fkuA_{k-1,k}) \fkvA_{k}.
\ee

\begin{corollary} 
\label{lp-cor1}
If $u \in S_n$ and $k \in [n]$ then
$
[\cL_{\e_k}] \cdot \cKA{u}
=
\fkMA_{n,k}\cKA{u}.
$
\end{corollary}

\begin{proof}
For $1\leq j<k<l\leq n$ define
$
\raff_{jk}  = \raff_{(\e_k-\e_j,0)}  
$ and
$\raff_{kl} = \raff_{(\e_l-\e_k,-1)}.$ 
Now choose integers
\[ 1\leq a_p <   \dots < a_1 < k < b_{q}< b_{q-1} <\dots <b_{1}\leq n \] and define  
$
w = u(a_1,k)(a_2,k)\cdots (a_p,k) (k,b_1)(k,b_2)\cdots (k,b_q) \in S_n.
$
Since $\raff_{kl}(\e_k)=\e_k$, we have 
\[
\mu := u\cdot \raff_{a_1,k}\cdot\raff_{a_2,k}\cdots \raff_{a_p,k} \cdot\raff_{k,b_1}\cdot\raff_{k,b_2}\cdots \raff_{k,b_q}(\e_k)
=u(a_1,k)(a_2,k)\cdots (a_p,k)(\e_k).
\]
It follows that 
\be\label{cKA-cor-eq}
 (-\fktA_{k,b_q})\cdots  (-\fktA_{k,b_2})(- \fktA_{k,b_1}) \fkuA_{a_p,k}\cdots  \fkuA_{a_2,k} \fkuA_{a_1,k} \fkvA_{k} \cKA{u}
 = \begin{cases}
 (-1)^q \cdot x^{-\mu}\cdot \cKA{w} &\text{if (*) holds} 
 \\ 
 0&\text{otherwise}
 \end{cases}
 \ee
 where (*) means that 
the successive multiplications in the product defining $w$ each increase the length by exactly one.
Since in type $\sfA_{n-1}$ when $\lambda=\e_k$, the set $\cS_{\lambda}$
from Definition~\ref{cS-def} is
\be
\cS_{\e_k}=
\left\{ 
\ba &(\e_{k}-\e_{k-1},0) <  (\e_{k}-\e_{k-2},0) < \dots <   (\e_{k}-\e_{1},0)
\\ &<
(\e_{n}-\e_{k},-1) < (\e_{n-1}-\e_{k},-1)< \dots <  (\e_{k+1}-\e_{k},-1)
\ea
\right\},
\ee
the corollary follows by comparing \eqref{cKA-cor-eq}
with the definition of $b_{u,w}^{\lambda,\mu}$ in Theorem~\ref{lp-thm2}.
\end{proof}

\begin{example}\label{A-M-ex}
 When $n=k=3$ we compute that
\[
\ba \fkMA_{3,3}\cKA{1}
&=
(1 + \fkuA_{1,3})(1 + \fkuA_{2,3}) \fkvA_{3}
\cKA{1}
\\&
=
(1 + \fkuA_{1,3})(1 + \fkuA_{2,3})\( x^{-\e_3} \cdot \cKA{1}\)
\\&
=
(1 + \fkuA_{1,3})\( x^{-\e_3} \cdot \cKA{1} +  x^{-\e_2} \cdot \cKA{132}\)
=
 x^{-\e_3} \cdot \cKA{1} +  x^{-\e_2} \cdot \cKA{132}+  x^{-\e_1} \cdot \cKA{231}
\ea
\]
and so
\[
[\cL_{\e_3}]\cdot \cKA{1} 
=x^{-\e_3}\cdot \cKA{1}
        + x^{-\e_2}\cdot \cKA{s_2}
        + x^{-\e_1}\cdot \cKA{s_1s_2}.
\]
Using Lemma~\ref{wL-lem} (recall the definition of $[\cO_w] = *[\cO_{X_w}]$ from Remark~\ref{notation-rem}), we then get
\[
e^{\e_3}\cdot [\cO_{{w_0}}] 
=x^{\e_1}\cdot [\cO_{{w_0}}] 
        + x^{\e_2}\cdot [\cO_{{s_2s_1}}] 
        + x^{\e_3}\cdot [\cO_{{s_1}}] .
\]
This recovers \cite[Ex.~15.6]{LP} once we note that the weights denoted $[a,b,c]$ in \cite[\S15]{LP} correspond to what we write as $c\e_1 + b\e_2 + a\e_3 \in \Lambda$.
When $n>k=3$ there is a longer, stable  expansion 
{\[
\ba{} [\cL_{\e_3}]\cdot\cKA{1} &=  \fkMA_{n,3}\cKA{1}
\\&=
(1 - \fktA_{3,4})  (1 - \fktA_{3,5})\cdots (1 - \fktA_{3,n})\( x^{-\e_3} \cdot \cKA{1} +  x^{-\e_2} \cdot \cKA{132}+  x^{-\e_1} \cdot \cKA{231}\)
\\
&=
(1 - \fktA_{3,4}) \( x^{-\e_3} \cdot \cKA{1} +  x^{-\e_2} \cdot \cKA{132}+  x^{-\e_1} \cdot \cKA{231}\)
\\
&=x^{-\e_3} \cdot \cKA{1} +  x^{-\e_2} \cdot \cKA{132}+  x^{-\e_1} \cdot \cKA{231}
+x^{-\e_3} \cdot \cKA{1243} +  x^{-\e_2} \cdot \cKA{1342}+  x^{-\e_1} \cdot \cKA{2341}
        .
        \ea
\]}
\end{example}

We streamline our notation for the other classical types.
Fix $X \in \{\sfB,\sfC,\sfD\}$ and 
define the following $\ZZ$-linear operators on $K_T(G/B)$
for integers $i \in \ZZ$ and $j,k\in[n]$ with $-n\leq i <j$:
  \be
  \ba
  \fktX_{ij} &: x^\mu  \cdot \cKX{w} \mapsto \begin{cases}
x^{\mu} \cdot\cKX{wt_{ij}}&\text{if $t_{ij}\in W^X_n$ and }\ell^X(wt_{ij})=\ell^X(w) + 1\\
0 & \text{otherwise}
\end{cases}
\\[-10pt]\\
  \fkuX_{ij} &: x^\mu  \cdot \cKX{w} \mapsto \begin{cases}
x^{wt_{ij} w^{-1}(\mu)} \cdot\cKX{wt_{ij}} &\text{if $t_{ij}\in W^X_n$ and }\ell^X(wt_{ij})=\ell^X(w) + 1\\
0 & \text{otherwise}
\end{cases}
\\[-10pt]\\
 \fkvX_{k}  &: x^\mu \cdot  \cKX{w} \mapsto x^{\mu - w(\e_k)} \cdot \cKX{w}
\\[-10pt]\\
\fkoX_{k}  &: x^\mu\cdot  \cKX{w} \mapsto \begin{cases}
 \cKX{w(-k,k)} &\text{if $X=\sfB$ and }\ell^\sfB(w(-k,k))=\ell^\sfB(w) + 1\\
0 & \text{otherwise},
\end{cases}
\ea
\ee
 Note that $\fktX_{-j,j} = \fkuX_{-j,j} =0$ as $t_{-j,j}=1$, and that $\fkoC_k=\fkoD_k=\fktD_{0,j} = \fkuD_{0,j} =0$.
 Then let 
{\small \be
\fkMX_k=
(1 - \fktX_{k,k+1})(1 - \fktX_{k,k+2})  \cdots (1 - \fktX_{k,n}) (1 + \fkoX_{k})(1 + \fkuX_{-n,k})  \cdots (1 + \fkuX_{k-2,k})(1 + \fkuX_{k-1,k}) \fkvX_{k} .
\ee}

\begin{corollary}\label{lp-cor2} 
Fix a type $X \in \{\sfB,\sfC,\sfD\}$. If $u \in W^X_n$ and $k \in [n]$ then
$
[\cL_{\e_k}]  \cdot\cKX{u}
=
\fkMX_k  \cKX{u}
$.
\end{corollary}

\begin{proof}
One just needs to verify that our operator formula is equivalent to Theorem~\ref{lp-thm2}
in types $\sfB_n$, $\sfC_n$, and $\sfD_n$ when $\lambda=\e_k$.
A straightforward calculation shows that $\cS_{\e_k}$ is given in type $\sfB_n$ by
\be
\cS_{\e_k}=
\left\{ 
\ba &(\e_{k}-\e_{k-1},0) <  (\e_{k}-\e_{k-2},0) < \dots <   (\e_{k}-\e_{1},0)
\\ &<(\e_k,0)
\\ &<(\e_{k}+\e_{1},0)< (\e_{k}+\e_{2},0) < \dots < (\e_{k}+\e_{k-1},0)
\\ &<(\e_{k+1}+\e_{k},0)< (\e_{k+2}+\e_{k},0) < \dots < (\e_{n}+\e_{k},0)
\\ &<(\e_k,1)
\\ &<
(\e_{n}-\e_{k},-1) < (\e_{n-1}-\e_{k},-1)< \dots <  (\e_{k+1}-\e_{k},-1)
\ea
\right\},
\ee
in type $\sfC_n$ by
\be
\cS_{\e_k}=
\left\{ 
\ba &(\e_{k}-\e_{k-1},0) <  (\e_{k}-\e_{k-2},0) < \dots <   (\e_{k}-\e_{1},0)
\\ &<(2\e_k,0)
\\ &<(\e_{k}+\e_{1},0)< (\e_{k}+\e_{2},0) < \dots < (\e_{k}+\e_{k-1},0)
\\ &<(\e_{k+1}+\e_{k},0)< (\e_{k+2}+\e_{k},0) < \dots < (\e_{n}+\e_{k},0)
\\ &<
(\e_{n}-\e_{k},-1) < (\e_{n-1}-\e_{k},-1)< \dots <  (\e_{k+1}-\e_{k},-1)
\ea
\right\},
\ee
and in 
type $\sfD_n$ by
\be
\cS_{\e_k}=
\left\{ 
\ba &(\e_{k}-\e_{k-1},0) <  (\e_{k}-\e_{k-2},0) < \dots <   (\e_{k}-\e_{1},0)
\\ &<(\e_{k}+\e_{1},0)< (\e_{k}+\e_{2},0) < \dots < (\e_{k}+\e_{k-1},0)
\\ &<(\e_{k+1}+\e_{k},0)< (\e_{k+2}+\e_{k},0) < \dots < (\e_{n}+\e_{k},0)
\\ &<
(\e_{n}-\e_{k},-1) < (\e_{n-1}-\e_{k},-1)< \dots <  (\e_{k+1}-\e_{k},-1)
\ea
\right\}.
\ee
Using these formulas and noting that $\raff_{(\e_l-\e_k,-1)}(\e_k) = \e_k $ for $k<l$ and $ \raff_{(\e_k,1)}(\e_k) = 0$, we deduce  the result by
an argument is similar to the proof of Corollary~\ref{lp-cor1}.
\end{proof}

\begin{example} \label{B-M-ex}
When $X=\sfB$, $n=2$, and $k=1$ we can use Corollary~\ref{lp-cor2} to calculate 
\[\ba
{} [\cL_{\e_1}]\cdot  \cKB{12} &=\fkMB_1\cKB{12}
\\&
=(1 - \fktB_{1,2})(1 - \fkoB_{1})(1 + \fkuB_{-2,1}) (1 + \fkuB_{0,1}) \fkvB_{1} \cKB{12}  
\\&
=(1 - \fktB_{1,2})(1 - \fkoB_{1})(1 + \fkuB_{-2,1}) (1 + \fkuB_{0,1}) \(x^{-\e_1}\cdot \cKB{12}  \)
\\&
=(1 - \fktB_{1,2})(1 - \fkoB_{1})(1 + \fkuB_{-2,1})   \(x^{-\e_1}\cdot \cKB{12}  + x^{\e_1}\cdot \cKB{\overline{1}2} \)
\\&
=(1 - \fktB_{1,2})(1 - \fkoB_{1})   \(x^{-\e_1}\cdot \cKB{12}  + x^{\e_1}\cdot \cKB{\overline{1}2}   + x^{\e_2}\cdot \cKB{\overline{2}1}  \)
\\&
=(1 - \fktB_{1,2})    \(x^{-\e_1}\cdot \cKB{12}  - (1-x^{\e_1})\cdot \cKB{\overline{1}2}    + x^{\e_2}\cdot \cKB{\overline{2}1} \)
\\&
=
x^{-\e_1}\cdot \(\cKB{12}-\cKB{21}\)  -(1-x^{\e_1})\cdot \(\cKB{\overline{1}2}-\cKB{2\overline{1}}\)    + x^{\e_2}\cdot \(\cKB{\overline{2}1}-\cKB{1\overline{2}}\).
\ea
\]
Similarly, one can check that 
$
[\cL_{\e_1}]\cdot  \cKB{\overline{1}2} =\fkMB_1\cKB{\overline{1}2}
=
x^{\e_1}\cdot\( \cKB{\overline{1}2} -\cKB{\overline{2}1}\) + 
x^{-\e_2}\cdot \(\cKB{\overline{2}1} - \cKB{\overline{1}2}\).
$
\end{example}

 \section{Transition equations}\label{te-sect}
 
In this section we translate the multiplication formulas for $K_T(G/B)$ discussed above
into certain polynomials identities for the \defn{double Grothendieck polynomials} of Lascoux--Sch\"utzenberger \cite{LS1983}
and their analogues in classical type due to Kirillov--Naruse \cite{KN}.
The main new results here are Theorems~\ref{bcd-thm1}, \ref{bcd-thm2}, and \ref{k-stanley-positivity-thm}
and their corollaries in Section~\ref{k-positivity-sect}.
 
\subsection{Double Grothendieck polynomials}\label{a-transitions-sect}

Let $\beta$ for a formal parameter and let $x_1,x_2,x_3,\dots$ and $y_1,y_2,y_3,\dots$
be two sequences of commuting variables.
Let $\ZZ[\beta][\bfy]=\ZZ[\beta,y_1,y_2,y_3,\dots]$
and
$\ZZ[\beta][\bfx;\bfy] = \ZZ[\beta,x_1,x_2,x_3,\dots,y_1,y_2,y_3,\dots]$
denote the usual polynomial rings in these variables.

 We 
denote
 typical elements of $\ZZ[\beta][\bfy]$ and $\ZZ[\beta][\bfx;\bfy]$ as $f(\bfy)$
 and
 $f(\bfx;\bfy)$. The group $S_\infty$ acts on $\ZZ[\beta][\bfx;\bfy]$ by permuting the $x$-variables, 
 and for
each $i \in \PP$ the formula 
\be \pi_i f(\bfx;\bfy) := \tfrac{(1+\beta x_{i+1}) f(\bfx;\bfy) - (1+\beta x_i) s_i f(\bfx;\bfy)}{x_i-x_{i+1}}\ee
defines a $\ZZ[\beta][\bfy]$-linear map $\ZZ[\beta][\bfx;\bfy] \to \ZZ[\beta][\bfx;\bfy]$.
These operators satisfy the Coxeter braid relations for $S_\infty$ along with the identity $\pi_i \pi_i = -\beta \pi_i$.

\begin{definition}[\cite{FK1994,LS1983}]
\label{fkG-def}
The \defn{double Grothendieck polynomials} $\fkG_w(\bfx;\bfy)$ for $w \in S_\infty$ are the unique elements of $\ZZ[\beta][\bfx;\bfy]$
satisfying the following conditions:
\ben
\item[(a)] if $n \in \PP$ then $\displaystyle\fkG_{n\cdots 321}(\bfx;\bfy) = \prod_{\substack{i,j \in [n-1] \\ i+j \leq n}} x_i \oplus y_j$ where $x\oplus y = x+y+\beta xy$; and
\item[(b)] if $i \in \PP$ is such that $w(i) > w(i+1)$ then $ \fkG_{w s_i}(\bfx;\bfy) =\pi_i \fkG_w(\bfx;\bfy)$.
\een
We also form
$
\fkG_w(\bfx) := \fkG_w(\bfx;0) \in \ZZ[\beta][x_1,x_2,x_3,\dots]
$
by setting all $y$-variables to zero.

\end{definition}

Taking $n=1$ in condition (a) implies that $\fkG_1(\bfx;\bfy)=1$.
Condition (a) implies that $\fkG_w(\bfx;\bfy)$ is always homogeneous of degree $\ell(w)$
if we define $\deg\beta=-1$ and $\deg x_i = \deg y_i= 1$.
Setting $\beta=0$
transforms $\fkG_w(\bfx;\bfy)$ to the \defn{double Schubert polynomial} $\fkS_w(\bfx;\bfy)$; see \cite[\S2.3]{Manivel}.

\begin{example}\label{sk-ex}
A variety of more explicit formulas for $\fkG_w(\bfx;\bfy)$ are known; see, e.g., \cite{BucScr,FK1994,KnutsonMiller1,KMY,Lascoux,Weigandt}.
It follows from \cite[Thm.~1.1]{Weigandt} that we actually have 
\be\label{nn-eq}
\fkG_w(\bfx;\bfy) \in \NN[\beta][x_1,x_2,\dots,x_{n-1},y_1,y_2,\dots,y_{n-1}]\quad\text{for all }w \in S_n.
\ee
In addition, it is essentially immediate from \cite[Cor.~5.4]{KnutsonMiller1} that
\be
\textstyle
1+\beta \fkG_{s_k}(\bfx;\bfy) = \prod_{i\in[k] } (1+\beta x_i)(1+\beta y_i).
\ee
This means that we can write
\be\label{x1A-eq}
1+\beta x_1= \tfrac{1+\beta \fkG_{s_1}(\bfx;\bfy) }{1+\beta y_1}
\quand
1+\beta x_k =\tfrac{1}{1+\beta y_k}\cdot \tfrac{1+\beta \fkG_{s_k}(\bfx;\bfy) }{1+\beta \fkG_{s_{k-1}}(\bfx;\bfy)}\text{ for }k>1.
\ee
\end{example}

The polynomials $\fkG_w(\bfx;\bfy)$ for $w \in S_n$ are related to the Schubert basis
of the equivariant $K$-theory ring in type $\sfA_{n-1}$ by this identity \cite[Thm.~3]{FL} (see also \cite[Prop.~3.13]{Hudson}):
\be\label{G-rep-eq}
\cKA{w} = \fkG(x_i \mapsto 1-e^{-\e_i}; y_i \mapsto 1-x^{\e_i})|_{\beta=-1} \in K_T(\SL_n(\CC)/B).
\ee
This property implies 
the following lemma, which is implicit in the literature. Consider the rings
 \be
 \ba 
 \AA^{(n)} &:=  \ZZ[\beta]\left[y_1, \tfrac{1}{1+\beta y_1}, y_2, \tfrac{1}{1+\beta y_2},\dots,y_{n-1}, \tfrac{1}{1+\beta y_{n-1}} \right],
\\
  \AA &:=  \ZZ[\beta]\left[y_1, \tfrac{1}{1+\beta y_1}, y_2, \tfrac{1}{1+\beta y_2},y_3, \tfrac{1}{1+\beta y_3},\dots \right] = \textstyle\bigcup_{n\geq 1} \AA^{(n)}
 .
 \ea
 \ee
 As with $\ZZ[\beta][\bfy]$, we denote a typical element of $\AA$ as $f(\bfy)$.
 
 An element of $\AA$ is \defn{homogeneous}
 if it is a $\ZZ$-linear combination of terms of the same degree where  $\deg \beta=-1$, $\deg \frac{1}{1+\beta y_i} = 0$, and $\deg x_i=\deg y_i = 1$. If $f(\bfy) \in \ZZ[\beta][\bfy]$ is homogeneous of degree $d$, then it can be recovered from $g(\bfy) := f(\bfy)|_{\beta=-1}$
 by the formula $f(\bfy) = (-\beta)^{-d} g(y_i\mapsto -\beta y_i)$.
 
\begin{lemma}\label{structure-lem1}
Fix permutations $u,v \in S_\infty$. Then there are unique homogeneous coefficients
 $c_{u,v}^w(\bfy) \in \AA$ for   $w\in S_\infty$
with
$\fkG_u(\bfx;\bfy)\cdot\fkG_v(\bfx;\bfy) = \sum_{w\in S_\infty}c_{u,v}^{w}(\bfy) \cdot \fkG_w(\bfx;\bfy)$,
and it holds that
\ben
\item[(a)] $\left\{ c_{u,v}^{w}(\bfy) : w \in S_n\right\} \subset \AA^{(n)}$, and
\item[(b)] $\cKA{u}  \cKA{v} = \sum_{w\in S_n}c_{u,v}^{w}(y_i \mapsto 1-x^{\e_i})|_{\beta=-1}   \cKA{w}  
\in K_T(\SL_n(\CC)/B)$
for all $n\in \PP$ with $u,v\in S_n$.
 \een
\end{lemma}

\begin{proof}
It is known that the set of all double Grothendieck polynomials is a homogeneous $\AA$-basis for a ring (see, e.g., \cite[\S8.6]{LLS}), so the coefficients $c_{u,v}^w(\bfy)\in\AA$ exist and are unique.
If $u,v\in S_n$ then \eqref{nn-eq} implies that $\{ c_{u,v}^{w}(\bfy) : w \in S_\infty\} \subset \AA^{(n)}$,
and the remaining identity follows from \eqref{G-rep-eq}.
\end{proof}

We have already defined an action of $S_\infty$ on $\ZZ[\beta][\bfx;\bfy]$ which permutes the $x$-variables.
Let $\star $ denote the analogous action of $S_\infty$ on $\AA$ that permutes the $y$-variables, so that 
\be\label{ast-eq} w \star f(\bfy) = f(y_{w(1)},y_{w(2)},y_{w(3)},\dots)\quad\text{for }w \in S_\infty.\ee
We now define several $\ZZ[\beta]$-linear operators on the $\AA$-module
\be
\label{AAx-eq}
\VV := \AA\spanning\{\fkG_w(\bfx;\bfy):  w \in S_\infty\} \subseteq \AA[\bfx] := \AA[x_1,x_2,x_3,\dots].
\ee
These operators will be denoted $\bft_{ij}$, $\bfu_{ij}$, $\bfv_k$, and $\bfM_k$.
Fix positive integers $ i < j$  and $k$. 
First,
let $\bft_{ij}$ be the $\AA$-linear operator sending
\be
\bft_{ij}: 
\fkG_w(\bfx; \bfy)\mapsto
\begin{cases}
\fkG_{w(i,j)}(\bfx; \bfy) &\text{if }\ell(w) +1 = \ell(w(i,j))
\\
0&\text{otherwise} 
\end{cases}
\ee
for each $w \in S_\infty$. Next, let $\bfu_{ij}$ be the $\ZZ[\beta]$-linear operator on $\VV$ that sends
\be
\bfu_{ij}: f(\bfy)\cdot \fkG_w(\bfx; \bfy)\mapsto
\(w(i,j)w^{-1} \star f(\bfy)\) \cdot \bft_{ij}\fkG_{w}(\bfx; \bfy)
\ee
for each $f(\bfy) \in \AA$ and $w \in S_\infty$. 
Finally, let $\bfv_{k} :\VV\to\VV$ be the $\AA$-linear map with
\be
\bfv_{k}:  \fkG_w(\bfx; \bfy)\mapsto
\tfrac{1}{1+\beta y_{w(k)}}\cdot \fkG_w(\bfx; \bfy).
\ee
Notice that $\bfu_{ij}$ is not $\AA$-linear.
Finally, define $\bfM_k$ to be the composition
\be\label{fkMa-eq}
\bfM_k=\prod_{k<l<\infty}^{[<]} (1 +\beta \bft_{kl})  \cdot \prod_{0< j < k}^{[<]} (1 -\beta \bfu_{jk}) \cdot  \bfv_{k}
\ee
where $\displaystyle\prod^{[<]}$ indicates a product that is expanded from left to right as indices increase. 
Although the first product has an infinite number of terms, we have
$
(1 +\beta \bft_{kl}) \fkG_w(\bfx;\bfy) = \fkG_w(\bfx;\bfy)
$
whenever $w \in S_n$ and $l>\max\{k,n\} +1$. Therefore $\bfM_k$ is a well-defined
$\ZZ[\beta]$-linear map $\VV\to\VV$.

\begin{theorem}\label{len-cor1} 
If $u \in S_\infty$ and $k \in \PP$ then
$
(1+\beta x_k)   \fkG_u(\bfx;\bfy)
=
\bfM_k \fkG_u(\bfx;\bfy).
$
\end{theorem}

\begin{proof}
In view of \eqref{x1A-eq}, it suffices to show that 
\be\label{len-eq}
\ba
(1+\beta \fkG_{s_1}(\bfx;\bfy)) \cdot  \fkG_u(\bfx;\bfy) &= 
(1+\beta y_1)\cdot\bfM_1 \fkG_u(\bfx;\bfy) &&\text{and }
\\
(1+\beta \fkG_{s_k}(\bfx;\bfy)) \cdot  \fkG_u(\bfx;\bfy) &= 
(1+\beta y_k)\cdot(1+\beta \fkG_{s_{k-1}}(\bfx;\bfy))\cdot\bfM_k \fkG_u(\bfx;\bfy) &&\text{for $k\geq 2$}.
\ea
\ee
This involves checking, for each $w\in S_\infty$,
 the equality of certain homogeneous expressions in $\AA$ involving the coefficients $c_{s_ku}^w(\bfy)$ and $c_{s_{k-1}v}^w(\bfy)$
 where $v$ ranges over the indices of the terms in the expansion of $\bfM_k \fkG_u(\bfx;\bfy)$.
By Lemma~\ref{structure-lem1}, we know that these equations are true 
in the representation ring $\ZZ[X]$ of the type $\sfA_{n-1}$ torus (for all sufficiently large $n$)
if we substitute $\beta\mapsto -1$ and $y_i \mapsto 1-x^{\e_i}$,
since Proposition~\ref{cL-prop}(a) and Theorem~\ref{lp-cor1} imply that the identities
\[
\ba
(1- \cKA{s_1}) \cdot \cKA{u}
&=
x^{\e_1}\cdot \fkMA_{n,k}\cKA{u}
&&\text{and} \\
(1- \cKA{s_k}) \cdot \cKA{u}
&=
x^{\e_k} \cdot (1- \cKA{s_{k-1}})\cdot\fkMA_{n,k}\cKA{u}&&\text{for $k\geq 2$}
\ea
\]
hold in $K_T(\SL_n(\CC)/B)$ for all $n\gg0$. Recall that $x^{\e_1+\e_2+\dots+\e_n} = 1 $ in $ X[T]$. As setting $\beta$ to $-1$ loses no information by homogeneity,
we deduce that our original equations hold in $\AA^{(n)}$ for all $n\gg 0$
modulo the \defn{coinvariant ideal} $ \ILambda_n \subseteq \AA^{(n)}$
generated by the symmetric polynomials in $y_1,y_2,\dots,y_n$ without constant term.
Since $\bigcap_{n\geq N} \ILambda_n =0$ for any given $N\in\PP$, we conclude that our equations are actually valid in $\AA$,
so \eqref{len-eq} holds and the theorem follows.
\end{proof} 

\begin{example} Taking $k=3$ and $u=1$ in Theorem~\ref{len-cor1} gives  the identity
\[\ba 1+\beta x_3 &= (1+\beta x_3) \fkG_1 
\\&=(1 + \beta y_{3})^{-1} \cdot \fkG_{1}
        -\beta\cdot (1 + \beta y_{2})^{-1} \cdot \fkG_{1 3 2}
        + \beta^2 \cdot(1 + \beta y_{1})^{-1} \cdot \fkG_{2 3 1}
\\&\quad        + \beta\cdot (1 + \beta y_{3})^{-1} \cdot \fkG_{1 2 4 3}
        -\beta^2 \cdot(1 + \beta y_{2})^{-1} \cdot \fkG_{1 3 4 2}
        + \beta^3 \cdot(1 + \beta y_{1})^{-1} \cdot \fkG_{2 3 4 1}.
\ea
\]
Compare with the second calculation in Example~\ref{A-M-ex}.
\end{example}

\begin{example} This instance of Theorem~\ref{len-cor1} recovers \cite[Ex.~3.9]{Lenart}
when $\beta=-1$ and $\bfy=0$:
\[\ba
(1+\beta x_3) \fkG_{13452}(\bfx;\bfy) &=
(1 + \beta y_{4})^{-1} \cdot \fkG_{1 3 4 5 2}(\bfx;\bfy)
   + \beta \cdot (1 + \beta y_{4})^{-1} \cdot \fkG_{1 3 5 4 2}(\bfx;\bfy)
\\&\quad        -\beta \cdot(1 + \beta y_{3})^{-1} \cdot \fkG_{1 4 3 5 2}(\bfx;\bfy)
 -\beta^2 \cdot(1 + \beta y_{3})^{-1} \cdot \fkG_{1 4 5 3 2}(\bfx;\bfy)
\\&\quad       + \beta^2 \cdot(1 + \beta y_{1})^{-1} \cdot \fkG_{3 4 1 5 2}(\bfx;\bfy)
 + \beta^3 \cdot(1 + \beta y_{1})^{-1} \cdot \fkG_{3 4 5 1 2}(\bfx;\bfy)
\\&\quad       + \beta^3 \cdot(1 + \beta y_{1})^{-1} \cdot \fkG_{3 4 2 5 1}(\bfx;\bfy)
 + \beta^4 \cdot(1 + \beta y_{1})^{-1} \cdot \fkG_{3 4 5 2 1}(\bfx;\bfy).
\ea\]
\end{example}

Theorem~\ref{len-cor1} is the natural equivariant generalization of \cite[Thm.~3.1]{Lenart}, 
but does not seem to have appeared in the literature before.
We can rephrase this result to more closely resemble the notation in \cite{Lenart}.
For permutations $u,v \in S_\infty$ and positive integers $ a<b$, we write 
$
u \xrightarrow{(a,b)} v
$
 to indicate that $v = u(a,b)$ and $\ell(v) = \ell(u)+1$.
The following is equivalent to Theorem~\ref{len-cor1}:

\begin{corollary}
If $u \in S_n$ and $k \in [n]$ then
\[
(1+\beta x_k)\fkG_u(\bfx;\bfy)
 = \sum_\gamma 
\mathsf{sgn}(\gamma) 
\cdot  (1+\beta  y_\gamma)^{-1} \cdot \beta^{\mathsf{len}(\gamma)} \cdot \fkG_{\last(\gamma)}(\bfx;\bfy)
\]
where the sum is over all chains $\gamma$
of the form
\[
u = u_0 \xrightarrow{(a_1,k)}
\cdots  \xrightarrow{(a_p,k)} u_p  \xrightarrow{(k,b_{1})} 
\cdots 
 \xrightarrow{(k,b_{q})}u_{p+q} = \last(\gamma)
\]
with $ 1\leq a_p <   \dots < a_1 <a_0= k < b_{q}< b_{q-1} <\dots <b_{1} $, and where we set
\[
\mathsf{sgn}(\gamma) = (-1)^p,
\quad
y_\gamma = y_{u(a_p)},
\quand
\mathsf{len}(\gamma) = p+q.
\]
\end{corollary}

We note another quick corollary.

\begin{corollary} It holds that $\VV=\AA[\bfx]$.
\end{corollary}

\begin{proof}
Theorem~\ref{len-cor1} shows that $x_k \fkG_u(\bfx;\bfy) = \beta^{-1} (\bfM_k-1)\fkG_u(\bfx;\bfy)$, which is a finite $\AA$-linear combination of
double Grothendieck polynomials. Hence $x_{i_1}x_{i_2}\cdots  x_{i_k} = x_{i_1}(x_{i_2}(\cdots(  x_{i_k} \fkG_1(\bfx;\bfy))\cdots))$ is also such a linear combination,
so $\AA[\bfx]\subseteq \VV$. The reverse containment holds by definition.
\end{proof}

Finally, we explain how to recover the transition equation \cite[Thm.~3.1]{Weigandt},
which reformulates an unpublished result of Lascoux \cite{Lascoux}.
Let $ \UU:=\ZZ[\beta]\spanning\{\fkG_w(\bfx;\bfy) : w\in S_\infty\}\subsetneq \VV$.
Then define
\[
\bfR_{k} = (1+\beta \bft_{k-1,k})\cdots (1+\beta \bft_{2,k})(1+ \beta\bft_{1,k}) =\prod_{0< j < k}^{[>]} (1+ \beta\bft_{jk}) : \VV\to\VV
\]
for each $k\in\PP$,
where $\displaystyle\prod^{[>]}$ indicates a product  expanded left to right as indices decrease.

\begin{lemma}\label{weig-lem}
Let $k \in \PP$, $w \in S_\infty$,
and $F \in \UU$.  
Then
$\displaystyle(1+\beta x_k) \bfR_k F= \(\prod_{k<l<\infty}^{[<]}(1 +\beta \bft_{kl})\cdot \bfv_{k}\) F$.
\end{lemma}

\begin{proof}
Since $\bfR_k$ preserves $ \UU$ 
and since
 $\bfM_k$ is a $ \ZZ[\beta]$-linear operator, Theorem~\ref{len-cor1} implies that 
\[\textstyle
(1+\beta x_k) \bfR_k F = \bfM_k \bfR_k F
=
 \(\prod_{k<l<\infty}^{[<]}(1 +\beta \bft_{kl})\cdot \prod_{0<j<k}^{[<]}(1 -\beta \bfu_{jk}) \cdot \bfv_{k}\cdot  \prod_{0<j<k}^{[>]}(1 +\beta \bft_{jk})\) F.
\]
Let 
 $t=(j,k) \in S_\infty$ where $1\leq j<k$.
 Since $1+\beta\bft_{jk}$ maps $\UU\to\UU$,  it suffices to show     
that \[
   ( (1-\beta \bfu_{jk})\cdot \bfv_{k} \cdot (1+\beta \bft_{jk}))F =  
 \bfv_k F,\]
 as then
 $
\( \prod_{0<j<k}^{[<]}(1 -\beta \bfu_{jk}) \cdot \bfv_{k}\cdot  \prod_{0<j<k}^{[>]}(1 +\beta \bft_{jk})\) F = \bfv_{k}  F
 $
 follows by induction.
 As $\bft_{jk}$, $\bfu_{jk}$, and $\bfv_k$ are $ \ZZ[\beta]$-linear,
we may assume  that
$F= \fkG_w(\bfx;\bfy)$ for some $w \in S_\infty$, and then the needed identity is 
 straightforward from the relevant definitions.
\end{proof}

For any integer $n \in \ZZ$ we let $[n,\infty) = \{ i \in \ZZ : i \geq n\}$.

\begin{theorem}[{\cite[Thm.~2.3]{Weigandt}}]\label{a-thm2}
Fix $1 \neq w \in S_\infty$. 
Let $a \in [1,\infty)$ be maximal with $w(a) > w(a+1)$, define $b \in [a+1,\infty)$ to be maximal
with $w(a)>w(b)$, and set $v = w(a,b)$ and $c=v(a)$. Then
\[
\fkG_w(\bfx; \bfy) = \beta^{-1} \Bigl((1+\beta y_{c})(1+\beta x_a) \bfR_{a}\fkG_v(\bfx; \bfy)-\fkG_v(\bfx; \bfy)\Bigr) 
.
\]
\end{theorem}

\begin{proof} Let $c = v(a) = w(b)$.
The maximality of $a$ and $b$ implies that if $a<i<\infty$ then $\ell(v(a,i)) = \ell(v)+1$ if and only if $i=b$,
and that $\ell(w(a,i))\neq \ell(w)+1$ for all $a<i<b$.
Hence, by Lemma~\ref{weig-lem}  
\[\textstyle
(1+\beta x_a) \bfR_a \fkG_v(\bfx; \bfy) = \(\prod_{a<i<\infty}^{[<]}(1 +\beta \bft_{ai})\cdot \bfv_a\) \fkG_v(\bfx; \bfy) = \tfrac{1}{1+\beta y_c} \( \fkG_v(\bfx;\bfy) + \beta \fkG_w(\bfx;\bfy)\).
\]
Substituting this identity into the right hand side of the desired formula gives $\fkG_w(\bfx;\bfy)$.
\end{proof}

\begin{remark}
In \cite[Thm.~2.3]{Weigandt}, the integer $b$ is defined to be the largest number with $(a,b) \in D(w)$ 
where $D(w) := \{ (i,j) \in \PP\times \PP : w(i) > j\text{ and }w^{-1}(j) > i\}$
is the \defn{Rothe diagram} of $ w \in S_\infty$.
One can check that this is same as requiring $b \in [a+1,\infty)$ to be maximal
with $w(a)>w(b)$.
\end{remark}

 Theorems~\ref{len-cor1} and \ref{a-thm2} generalize several known identities for (single) Grothendieck polynomials
 and (double) Schubert polynomials from, e.g., \cite{KV,LascouxTransitions,LS1985,Lenart}; see the discussion in \cite[\S2.5]{Weigandt}.

\subsection{K-Stanley symmetric functions}\label{k-stanley-sect}

Before we can describe the analogues of $\fkG_w(\bfx;\bfy)$ in the other classical types,
we must review Kirillov and Naruse's definition of \defn{$K$-Stanley symmetric functions} for signed permutations.

For this purpose, we add to the parameters $\beta$, $x_1,x_2,x_3,\dots$, and $y_1,y_2,y_3,\dots$ a third sequence of mutually commuting variables
$z_1,z_2,z_3,\dots$.
Write 
\[\ZZ[\beta]\llbracket \bfz\rrbracket = \ZZ[\beta]\llbracket z_1,z_2,z_3,\dots\rrbracket 
\quand
\ZZ[\beta]\llbracket \bfx;\bfy;\bfz\rrbracket = \ZZ[\beta]\llbracket x_1,x_2,,\dots,y_1,y_2,\dots z_1,z_2,\dots\rrbracket \]
for the associated formal power series rings in these variables.

Recall that the \defn{Demazure product} on a Coxeter group $W$ with length function $\ell :W\to\NN$ is the unique associative operation
$\circ : W\times W \to W$ that satisfies 
\be u\circ v = uv\text{ if }\ell(uv)=\ell(u)+\ell(v)\quand s\circ s = s\text{ if }\ell(s)=1.\ee
 Fix a type $X \in \{\sfB,\sfC,\sfD\}$, and for $w \in W^X_\infty$ let $\cH^X(w)$ be the set of finite integer sequences $a=(a_1,a_2,\dots,a_k)$
with $a_i \in \{-1,0,1,2,\dots\}$ such that  
\[w = t_{a_1}\circ t_{a_2}\circ \cdots\circ t_{a_k}\]
when the Demazure product $\circ $ is computed relative to the Coxeter group $W= W^X_\infty$.
All entries of $a \in \cH^X(w)$ must be in the set $\{0,1,2,3,\dots\}$ when $X \in \{\sfB,\sfC\}$,
and in the set $\{-1,1,2,3,\dots\}$ when $X=\sfD$.
Note that $\cH^\sfB(w)=\cH^\sfC(w)$.
We refer to elements of $\cH^X(w)$ as \defn{Hecke words} for $w$.

\begin{definition}\label{cs-def}
Suppose $a=(a_1,a_2,\dots,a_k)  $ is any integer sequence.
Let $\ell(a) = k$ and define
\be o^\sfB(a) = |\{ i \in [k] : a_i = 0\}|, \quad o^\sfC(a)=0,\quand o^\sfD(a) = |\{ i \in [k] : a_i = \pm1\}|.
\ee
Next, let $\cC^X(a)$ be the set of integer sequences $b=(b_1,b_2,\dots,b_k)$ with the following properties:
\bei
\item[(a)] one has $1\leq b_1 \leq b_2 \leq \dots \leq b_k$;
\item[(b)] if $|a_{i-1}|\leq |a_i| \geq |a_{i+1}|$ for some $i \in \{2,3,\dots,k-1\}$ then $b_{i-1} < b_{i+1}$;
\item[(c)] if $X=\sfB$ and $a_i=a_{i+1}=0$ then $b_i < b_{i+1}$; and
\item[(d)] if $X=\sfD$ and $a_i=a_{i+1}\in\{-1,1\}$ then $b_i <b_{i+1}$.
\eei
Given $b \in \cC^X(a)$, let 
$\bfz_b = z_{b_1}z_{b_2}\cdots z_{b_k}$
and 
$ |b| = |\{ b_1,b_2,\dots,b_k\}|$, and define
\be
 \gamma(a,b) = |\{ i \in [k-1] : a_i=a_{i+1}\text{ and }b_i=b_{i+1}\}|.
 \ee
 We refer to the elements of $\cC^X(a)$ as \defn{compatible sequences} for $a$.
\end{definition}

The following is not the first definition of $K$-Stanley symmetric functions given in \cite{KN},
but we find it to be more self-contained and explicit than \cite[Def.~7]{KN}.

\begin{definition}[{\cite[Prop.~17]{KN}}] \label{cFX-def}
Fix $X \in \{\sfB,\sfC,\sfD\}$. Write $\ell=\ell^X$ and suppose $w \in W^X_\infty$. 
The \defn{$K$-Stanley symmetric function} of $w$ for type $X$ is the formal power series
\[
\cFX_w(\bfz) =\sum_{a\in\cH^X(w)}\sum_{b\in\cC^X(a)} 2^{|b| -\gamma(a,b)-o^X(a)} \cdot \beta^{\ell(a) -\ell(w)} \cdot \bfz_b \in \ZZ[\beta]\llbracket\bfz\rrbracket.
\]
\end{definition}

\begin{remark}
This definition corrects a minor error in 
\cite[Prop.~17]{KN}, which
(mistakenly) uses the same notion of compatible sequences for types $\sfB$ and $\sfC$.
Specifically, \cite[Prop.~17]{KN} asserts that the $K$-Stanley symmetric function $\cFB_w(\bfz)$ is equal to the sum
\[\textstyle \sum_{a\in\cH^\sfB(w)}\sum_{b\in\cC^\sfC(a)}  2^{|b| -\gamma(a,b)-o^\sfB(a)} \cdot   \beta^{\ell(a) -\ell(w)} \cdot  \bfz_b.\] However, this expression can fail to be a power series with integer coefficients, as one can check by taking $w=t_0$.
Kirillov and Naruse omit the details for the type $\sfB$ case of their proof of \cite[Prop.~17]{KN},
but filling these in leads to our  formula.
\end{remark}

In type $\sfD$ it is immediate from the definition that 
\be \cFD_w(\bfz) = \cFD_{w^\ast}(\bfz) \quad\text{for all }w \in W^\sfD_\infty\ee
where $\ast$ is the automorphism that exchanges $t_{-1}$ and $t_1$ while fixing $t_i$ for all $i\geq 2$.
For all types, the less obvious symmetry
 $\cFX_w(\bfz) =\cFX_{w^{-1}}(\bfz)$  
holds by \cite[Lem.~7]{KN}.

\begin{example} 
If $\ell(w) = 1$ then $\cH^X(w) = \{ (a), (a,a), (a,a,a), \dots\}$ for some $a \in \{-1,0,1,2,\dots\}$. 
If $a \in \{-1,0,1\}$ and $X \in \{\sfB,\sfD\}$ then all compatible sequences must by strictly increasing, so 
\[
\cFB_{t_0}(\bfz)=\cFD_{t_1}(\bfz)=\cFD_{t_{-1}}(\bfz)= \sum_{n \in \NN} \sum_{i_0<i_1<\dots<i_n} \beta^n  z_{i_0}z_{i_1}\cdots z_{i_n}. 
\]
On the other hand, if $k \in \NN$ is any nonnegative integer then
\[
\ba 
\cFB_{t_{k+1}}(\bfz) 
= \cFC_{t_k}(\bfz) 
=\cFD_{t_{k+2}}(\bfz) 
&=  \sum_{n \in \NN} \sum_{i_0<i_1<\dots<i_n} \beta^n  (2z_{i_0} + \beta z_{i_0}^2) (2z_{i_1} + \beta z_{i_1}^2)\cdots (2z_{i_k} + \beta z_{i_n}^2). 
 \ea
 \]
 \end{example}

We mention an alternate formula for $\cFX_w(\bfz)$ that may appear more natural.
Write $\prec$ for the total order on $\ZZ$ that has $0\prec -1\prec 1 \prec -2 \prec 2 \prec \cdots$.

\begin{definition}
Given an integer sequence $a =(a_1,a_2,\dots,a_k)$, define $\cU(a)$
to be  the set of integer sequences $b=(b_1,b_2,\dots,b_k)$ with the following properties:
\ben
\item[(a)] one has $0 \neq b_1 \preceq b_2 \preceq \dots \preceq b_k$;
\item[(b)] if $b_i=b_{i+1}$ is negative for some $i\in[k-1]$ then $|a_i|>|a_{i+1}|$;
\item[(c)] if $b_i=b_{i+1}$ is positive for some $i\in[k-1]$ then $|a_i| <|a_{i+1}|$.
\een
Following \cite{AHHP}, we refer to elements of $\cU(a)$ as \defn{unimodal factorizations} of $a$.
We also set
\[ \ba
\cU^\sfB(a) & := \{ b \in \cU(a) : b_i\text{ is positive whenever $a_i =0$}\}, \\
\cU^\sfC(a)& := \cU(a), \\
\cU^\sfD(a)  &:= \{ b \in \cU(a) : b_i\text{ is positive whenever $a_i=\pm 1$}\}.
\ea
\]
\end{definition}

For each $b \in \cU^X(a)$ let
 $\bfz_b = z_{|b_1|} z_{|b_2|} \cdots z_{|b_k|}$ where $k=\ell(a)=\ell(b)$.

\begin{proposition}\label{cFX-def2}
Fix $X \in \{\sfB,\sfC,\sfD\}$. Write $\ell=\ell^X$ and suppose $w \in W^X_\infty$. 
Then
\[
\cFX_w(\bfz) =\sum_{a\in\cH^X(w)}\sum_{b\in\cU^X(a)}  \beta^{\ell(a) -\ell(w)} \cdot \bfz_b \in \ZZ[\beta]\llbracket\bfz\rrbracket.
\]
\end{proposition}

\begin{proof}
In type $\sfC$ this result is \cite[Prop.~3.2]{AHHP}. The proof here is similar.

Suppose $b \in \cU^X(a)$, and let $|b| = (|b_1|, |b_2|, |b_3|,\dots)$. The sequence $|b|$ is weakly increasing
and if $|b_{i-1}|=|b_i|=|b_{i+1}|$ then either $b_{i-1}=b_i<0$ (so $|a_{i-1}|>|a_i|$) 
or $b_i=b_{i+1}>0$ (so $|a_{i}|<|a_{i+1}|$).
Hence, if $|a_{i-1}|\leq |a_{i}| \geq |a_{i+1}|$ then we must have $|b_{i-1}|<|b_{i+1}|$.
Moreover, if $X=\sfB$ and $a_i=a_{i+1}=0$,
or if $X=\sfD$ and $|a_i|=|a_{i+1}|=1$,
 then we must have $b_i=|b_i| <|b_{i+1}|< b_{i+1}$.
 
 We conclude that $|b| \in \cC^X(a)$.
 It remains to check that the inverse image of any $b \in\cC^X(a)$ under the map $|\cdot|$
 has size $ 2^{|b| -\gamma(a,b)-o^X(a)}$. This is a straightforward exercise.
\end{proof}

 There are also $K$-Stanley symmetric functions in type $\sfA$, which we briefly review for comparison.
These power series, which will be denoted $G_w(\bfz)$, are often called \defn{stable Grothendieck polynomials} or \defn{symmetric Grothendieck functions} in the literature \cite{Buch2002,BKSTY,BucScr,FK1994}.

Given $w \in S_\infty$, let $\cH(w)$ be the set of finite integer sequences $a=(a_1,a_2,\dots,a_k)$
with
$w = s_{a_1}\circ s_{a_2}\circ \cdots\circ s_{a_k}$.
Then, given $a=(a_1,a_2,\dots,a_k) \in \cH(w)$, let $\cC(a)$
be the set of weakly increasing positive integer sequences $b=(b_1\leq b_2\leq \dots\leq b_k)$
with $b_j < b_{j+1}$ if $a_j \leq a_{j+1}$. Finally, let
\be\label{limit-eq}
  G_w(\bfz) = \sum_{a \in \cH(w)} \sum_{b \in \cC(a)}  \beta^{\ell(a) -\ell(w)}   \bfz_b \in \ZZ[\beta]\llbracket\bfz\rrbracket.
\ee
This object is related to the Grothendieck polynomials from Section~\ref{a-transitions-sect} by the formula \cite[\S2]{BKSTY}
\be\textstyle G_w(\bfz) = \lim_{N\to \infty} \fkG_{1^N\times w}(\bfz)\ee where the limit is in the sense of formal power series 
and $1^N\times w \in S_\infty$ denotes the that fixes $i \in [n]$ and maps $i+n \mapsto w(i)+n$ for $i \in \PP$.

The identity \eqref{limit-eq}, via property (b) in Definition~\ref{fkG-def}, implies that $G_w(\bfz)$ belongs to the ring  $\Sym_\beta \subset \ZZ[\beta]\llbracket\bfz\rrbracket$
of symmetric formal power series (of unbounded degree) with coefficients in $\ZZ[\beta]$.
The symmetric functions $G_w(\bfz)$
 generalize the type A \defn{Stanley symmetric functions} from \cite{Stanley}, which are recovered by setting $\beta=0$.

When $w \in \WA_\infty \cong S_\infty$, there is a close relationship between $\cFX_w(\bfz)$ and $G_w(\bfz)$.
Let $\FM$ be the set of finite integer sequences.
Such a sequence $a =(a_1,a_2,\dots,a_k) \in \FM$ is a \defn{multi-permutation}
if $a_i\neq a_{i+1}$ for all $i \in [k-1]$.
Let $\MP$ be the set of all multi-permutations and define $\mpmap : \FM \to \MP$
by the formula $\mpmap(a) = (a_{i_1}, a_{i_2},\dots,a_{i_r})$ where
\[ 
\{i_1 < i_2< \dots <i_r\} = \{ i \in [\ell(a)] : i=1\text{ or }a_{i-1}\neq a_i\}.
\]
Finally, given any multi-permutation $\pi \in \MP$,
consider the formal power series
\be L^{(\beta)}_\pi := \sum_{\substack{a \in \FM \\ \mpmap(a)=\pi}} \sum_{b \in \cC(a)}  \beta^{\ell(a) -\ell(\pi)}   \bfz_b 
\quand
 K^{(\beta)}_\pi := \sum_{\substack{a \in \FM \\ \mpmap(a)=\pi}} \sum_{b \in \cU(a)}  \beta^{\ell(a) -\ell(\pi)}   \bfz_b \ee
These coincide with the \defn{multi-fundamental quasisymmetric functions} and \defn{multi-peak quasisymmetric functions} studied in \cite{LamPyl,LM,MarbergHopf}.
They recover Gessel's \defn{fundamental quasisymmetric functions} and Stembridge's \defn{peak quasisymmetric functions}
\cite{StembridgeEnriched}
when $\beta=0$.

The power series $L^{(\beta)}_\pi $ only depends of $\ell(\pi)$ and the \defn{descent set} $\Des(\pi) = \{ i : \pi_i > \pi_{i+1}\}$ \cite[\S3.3]{LM}.
Similarly, $K^{(\beta)}_\pi $ only depends on $\ell(\pi)$ and the \defn{peak set} $\Peak(\pi) = \{ i :\pi_{i-1}< \pi_i > \pi_{i+1}\}$ \cite[\S4.3]{LM}.
Both elements belong to the ring of $\QSym_\beta \subset \ZZ[\beta]\llbracket\bfz\rrbracket$
of \defn{quasisymmetric functions} (of unbounded degree) with coefficients in $\ZZ[\beta]$; in fact, every element of $\QSym_\beta$ can be uniquely written as an infinite $\ZZ[\beta]$-linear combination of $L^{(\beta)}$-functions.

By \cite[Cor.~4.22]{LM} there is a (continuous) ring endomorphism
$
\Theta : \QSym_\beta \to \QSym_\beta
$
with
\be\label{Theta-eq} \Theta(L^{(\beta)}_\pi)= K^{(\beta)}_\pi\quad\text{for all $\pi \in \MP$.}
\ee
Below, we define $G_w(\bfz) $ for $w \in W^A_\infty$ by identifying $w$ with an element of $S_\infty$ in the obvious way.

\begin{proposition}
Let $w \in W^A_\infty$ and $X \in \{\sfB,\sfC,\sfD\}$.
Assume that $w(1)=1$ if $X=\sfD$. Then
\[
\cFX_w(\bfz) = \Theta (G_w(\bfz)).
\]
\end{proposition}

\begin{proof}
In these cases the Hecke words for $w$ do not have any letters equal to $0$ in types $\sfB$ and $\sfC$
or to $\pm 1$ in type $\sfD$, so $\cU^X(a) = \cU(a)$ if $a \in \cH^X(w)$ and 
$
\cFX_w(\bfz) = \sum_{\pi \in \cH^X(w) \cap \MP}  \beta^{\ell(\pi) -\ell(w)} K^{(\beta)}_\pi.
$
Since we likewise have
$
G_w(\bfz) = \sum_{\pi \in \cH(w) \cap \MP}  \beta^{\ell(\pi) -\ell(w)} L^{(\beta)}_\pi
$, the result follows from \eqref{Theta-eq}.
\end{proof}

More generally, it is clear for any $w \in W^\sfC_\infty$ that \be
\textstyle \cFC_w(\bfz) = \sum_{\pi \in \cH^\sfC(w) \cap \MP}  \beta^{\ell(\pi) -\ell^\sfC(w)} K^{(\beta)}_\pi.
\ee
Formulating the $K^{(\beta)}$-expansions of $\cFB_w(\bfz) $ and $\cFD_w(\bfz) $ when $w \notin W^A_\infty$
is a more subtle problem, which we will not pursue here. 
Knowing this expansion in type $\sfD$ would generalize \cite[Thm.~11.2]{StembridgeSome}.

\subsection{Supersymmetric functions}\label{ss-sect}

In this section, we explain some tableau formulas 
that give $\cFX_w(\bfz)$ for certain classes of signed permutations.

Recall from Section~\ref{intro2-sect} that $\SS_\beta$ denotes the subring of  power series  
   $f(\bfz)\in\Sym_\beta$ with the   symmetry property \eqref{ksym-eq}.
While not evident from Definition~\ref{cFX-def}, it turns out \cite[Lem.~6]{KN} that
\be \cFX_w(\bfz) \in \SS_\beta\quad\text{for all $X \in \{\sfB,\sfC,\sfD\}$ and $w \in W^X_\infty$}.\ee
Following \cite{IN,KN}, we refer to elements of $\SS_\beta$ as \defn{$K$-supersymmetric functions}.

The ring $\SS_\beta$ was first considered by Ikeda and Naruse in \cite{IN}
since it consists of exactly the formal power series that can be expressed as infinite $\ZZ[\beta]$-linear combinations
of their \defn{$K$-theoretic Schur $P$- and $Q$-functions}, whose definition we now recall.

A \defn{partition} is a weakly decreasing sequence of nonnegative integers $\lambda= (\lambda_1\geq \lambda_2\geq \dots)$ with finite sum $|\lambda| = \sum_i \lambda_i$. A partition is \defn{strict} if $\lambda_i=\lambda_{i+1}$ implies that $\lambda_i=\lambda_{i+1}=0$.
The \defn{shifted diagram} of a strict partition $\lambda$ is the set of pairs 
\be\SD_\lambda = \{(i,i+j-1) \in \PP\times \PP: 1\leq j \leq \lambda_i\}\ee
which we treat as positions in a matrix. We write $\mu \subseteq \lambda$ if $\SD_\mu\subseteq\SD_\lambda$ and define $\SD_{\lambda/\mu} = \SD_\lambda \setminus \SD_\mu$.

 A \defn{set-valued shifted tableau} of  shape $\lambda/\mu$ is a map $T$ from $\SD_{\lambda/\mu}$ to the set of finite, nonempty subsets of the totally ordered alphabet $ \{1'<1<2'<2<\dots\}$.
 If $T$ is such a tableau, then we write $(i,j) \in T$ to indicate that $(i,j)$ is in its domain and 
let $T_{ij}$ denote the set assigned by $T$ to this position. We also define $|T| = \sum_{(i,j) \in T} |T_{ij}|$ and
\be \bfz^T = \prod_{k=1}^\infty z_k^{a_k+b_k}
\quad\text{where $\begin{cases} a_k\text{ is the number of pairs $(i,j)\in T$ with $ k \in T_{ij}$}, \\ b_k\text{ is the number of pairs $(i,j)\in T$ with $ k' \in T_{ij}$}.\end{cases}$}
\ee
A set-valued shifted tableau is \defn{semistandard} if 
\ben
\item[(a)] $\max(T_{ij}) \leq \min(T_{i,j+1})$ with $T_{ij} \cap T_{i,j+1} \subset \{1,2,3,\dots\}$ whenever $(i,j),(i,j+1)\in T$; and
\item[(b)] $\max(T_{ij}) \leq \min(T_{i+1,j})$ with $T_{ij} \cap T_{i+1,j} \subset \{1',2',3',\dots\}$ whenever $(i,j),(i+1,j)\in T$.
\een
Write $\SetShTab^+(\lambda/\mu)$ for the set of  semistandard set-valued shifted tableaux of shape $\lambda/\mu$. Let
\be\SetShTab(\lambda/\mu)= \{ T \in \SetShTab^+(\lambda/\mu): T_{ij} \subset \PP\text{ for all $(i,j)\in T$ with $i=j$}\}.\ee We abbreviate by writing 
\be \SetShTab(\lambda) = \SetShTab(\lambda/\emptyset)
\quand
 \SetShTab^+(\lambda) = \SetShTab^+(\lambda/\emptyset).
 \ee
\begin{example} If $\lambda = (5,3,2)$ and $\mu = (2)$ then 
one element $T \in \SetShTab(\lambda/\mu)$ is given by
\[
\ytableausetup{boxsize = .8cm,aligntableaux=center}
T=\begin{ytableau}
*(lightgray) &*(lightgray) & 2'3'  & 4 & 457
\\
\none & 12 & 3'3 & 5'  \\
\none & \none & 4 & 5'57 
\end{ytableau}
\]
and we have 
 $T_{23} = \{3',3\}$, $|T| = 15$, and $ \bfz^T = z_1 z_2^2 z_3^3 z_4^3 z_5^4 z_7^2$.
\end{example}

For convenience, when $\mu \not\subset\lambda$ we let
$\SetShTab(\lambda/\mu)=\SetShTab^+(\lambda/\mu)=\varnothing$ be the empty set.

\begin{definition}[{\cite[\S9]{IN}}]
If $\mu$ and $\lambda$ are strict partitions, then
 the \defn{$K$-theoretic Schur $P$- and $Q$-functions} indexed by $\lambda/\mu$ are
 the formal power series in $\NN[\beta]\llbracket\bfz\rrbracket$ given by
\[
\GP_{\lambda/\mu}(\bfz) = \sum_{T \in \SetShTab(\lambda/\mu)}\beta^{|T|-|\lambda|} \bfz^T\quand
\GQ_{\lambda/\mu}(\bfz) = \sum_{T \in \SetShTab^+(\lambda)} \beta^{|T|-|\lambda|} \bfz^T 
.\]
\end{definition}

We abbreviate by writing $\GP_{\lambda}(\bfz)=\GP_{\lambda/\emptyset}(\bfz)$ and 
 $\GQ_{\lambda}(\bfz)=\GQ_{\lambda/\emptyset}(\bfz)$. We also set 
 \be
 \GP_0(\bfz) = \GQ_0(\bfz)=1,
 \quad
 \GP_{n}(\bfz) = \GP_{(n)}(\bfz),
 \quand
  \GQ_{n}(\bfz) = \GQ_{(n)}(\bfz)
  \ee
  when $n$ is a positive integer.
If   $\deg(z_i)=1$ and $\deg(\beta)=-1$, then
 $\GP_{\lambda/\mu}(\bfz)$ and $\GQ_{\lambda/\mu}(\bfz)$ are both homogeneous of degree $|\lambda/\mu|=|\lambda|-|\mu|$ 
 whenever $\mu\subset\lambda$. 
 
\begin{example}
If we write $e_k(\bfz) = \sum_{1 \leq i_1 < i_2 <\dots <i_k} z_{i_1} z_{i_2} \cdots z_{i_k}$ for $k \in \PP$ then   
\[ \GP_{1}(\bfz)  = e_1(\bfz) + \beta e_2(\bfz) +\beta^2 e_3(\bfz) +\beta^3 e_4(\bfz)  + \cdots \quand 
\GP_{2}(\bfz) =  \(\GP_{1}(\bfz)\)^2,\]
along with $\GQ_{1}(\bfz)= 2\cdot\GP_{1}(\bfz) + \beta \cdot\GP_{2}(\bfz)$. In general,
$\GQ_{n}(\bfz)= 2\cdot \GP_{n}(\bfz) + \beta\cdot \GP_{n+1}(\bfz)$.
\end{example}

Ikeda and Naruse \cite{IN} have proved that $\GP_\lambda(\bfz)$ and $\GQ_\lambda(\bfz)$ are both elements of $\SS_\beta$ for all
strict partitions $\lambda$.
The same holds for $\GP_{\lambda/\mu}(\bfz)$ and $\GQ_{\lambda/\mu}(\bfz)$ by \cite[Cor.~5.13]{LM}.
The set $\{ \GP_\lambda(\bfz) :\lambda \text{ strict}\}$ is a \defn{pseudobasis} for $\SS_\beta$ in the sense that
each element of $\SS_\beta$ can be uniquely expressed as a formal (meaning possibly infinite) $\ZZ[\beta]$-linear combination of $\GP$-functions. Each $\GQ_\lambda(\bfz)$ is a finite $\ZZ[\beta]$-linear combination of $\GP$-functions \cite{ChiuMarberg},
and it follows that 
each element of $\SS_\beta$ can be uniquely expressed as a formal $\QQ[\beta]$-linear combination of $\GQ$-functions.

Recall from \eqref{gr-eq0} and \eqref{gr-eq} in the introduction that 
a signed permutation $w \in \WB_\infty=\WC_\infty$ is \defn{Grassmannian} if 
$w(1)<w(2)<w(3)<\cdots$, and if $w \in W^X_\infty$ has this property then there are strict partitions
$\lambda_\sfB(w) = \lambda_\sfC(w)$ and $\lambda_\sfD(w) $
such that \cite[Lem.~7]{KN} we have
\[
\cFX_w(\bfz) = \begin{cases} \GP_{\lambda_X(w)}(\bfz)&\text{when }X \in \{\sfB,\sfD\} \\ 
 \GQ_{\lambda_X(w)}(\bfz)&\text{when $X =\sfC$.}
 \end{cases}
 \]
This formula implies that any $\GP_\lambda(\bfz)$ can occur as a $K$-Stanley symmetric function
of type $\sfB$ or $\sfD$, and any $\GQ_\lambda(\bfz)$ can occur as a $K$-Stanley symmetric function
of type $\sfC$. The next result will show that this is also true for $\GP$- and $\GQ$-functions indexed by skew shapes.

Fix two strict partitions $\mu \subseteq \lambda$. 
We define two signed permutations indexed by the skew shape $\lambda/\mu$ by
 the following procedure.
First, let $\TBC_{\lambda/\mu}$ be the shifted tableau of shape $\lambda/\mu$ 
with $j-i$ in each position $(i,j) \in \SD_{\lambda/\mu}$.
Then let $\TD_{\lambda/\mu}$ be the shifted tableau of shape $\lambda/\mu$ 
with $j-i+1$ in position $(i,j) \in \SD_{\lambda/\mu} $ if $i \neq j$ and with $(-1)^i$ in this position if $i=j$.
Finally, define
\be w_\sfB(\lambda/\mu) = w_\sfC(\lambda/\mu) = t_{a_1}t_{a_2}\cdots t_{a_{k}}
\quand
w_\sfD(\lambda/\mu) = t_{b_1}t_{b_2}\cdots t_{b_{k}}\ee
where $k=|\lambda/\mu|$ and 
where $a_1a_2\cdots a_k$ and $b_1b_2\cdots b_k$ are the entries of $\TBC_{\lambda/\mu}$ and $\TD_{\lambda/\mu}$
read row-by-row from left to right, starting with the smallest row index. 

\begin{example}\label{fc-ex}
When $\lambda=(5,3,1)$ and $\mu=(2)$ our two shifted tableaux are
\[
\ytableausetup{boxsize = .6cm,aligntableaux=center}
\TBC_{(5,3,1)/(2)} = 
\begin{ytableau}
*(lightgray) &*(lightgray) & 2 & 3 & 4 \\
 \none & 0 & 1 & 2 \\ 
 \none & \none & 0
\end{ytableau}
\quand
\TD_{(5,3,1)/(2)} = 
\begin{ytableau}
*(lightgray) &*(lightgray) & 3 & 4 & 5 \\
 \none & -1 & 2 & 3 \\ 
 \none & \none & +1
\end{ytableau},
\]
and the corresponding signed permutations are
\be\label{wbd-eq}
\ba w_\sfB(\lambda/\mu) =w_\sfC(\lambda/\mu)  &= t_2t_3t_4 t_0t_1t_2 t_0 = \overline{3}4\overline{1}52
&&\text{and} \\
w_\sfD(\lambda/\mu) &= t_3t_4t_5 t_{-1} t_2t_3 t_1 = 4\overline{2}5\overline{1}63.
\ea
\ee
\end{example}

An element of a Coxeter group is \defn{fully commutative} if each of its reduced words can be transformed to any other
by successively swapping pairs of adjacent commuting generators (thus avoiding the use of any nontrivial braid relations). In the symmetric group $S_\infty$,
the set of fully commutative elements coincides with the class of $321$-avoiding permutations \cite[Cor.~5.8]{StembridgeSome}.

In the terminology of \cite{StembridgeSome},
the signed permutation $w_\sfB(\lambda/\mu)$ is the fully commutative \defn{top element} of $W^\sfB_\infty$ of shape $\lambda/\mu$, and the sequence $a_1a_2\cdots a_k$ 
read off from $\TBC_{\lambda/\mu}$ is its \defn{canonical reduced word} \cite[Cor.~5.6]{StembridgeSome}.
Similarly, $w_\sfD(\lambda/\mu)$ is a fully commutative \defn{$A$-stable} element of $W^\sfD_\infty$ of shape $\lambda/\mu$, and the sequence $b_1b_2\cdots b_k$
read off from $\TD_{\lambda/\mu}$ is its \defn{canonical reduced word} \cite[Cor.~10.2]{StembridgeSome}.

\begin{theorem}
\label{fc-prop}
Suppose $\mu\subseteq \lambda$ are strict partitions. Then
\[
\cFB_{ w_\sfB(\lambda/\mu) }(\bfz) = \GP_{\lambda/\mu}(\bfz),
\quad
\cFC_{ w_\sfC(\lambda/\mu)}(\bfz) = \GQ_{\lambda/\mu}(\bfz),
\quand
\cFD_{ w_\sfD(\lambda/\mu)}(\bfz) = \GP_{\lambda/\mu}(\bfz).
\]
\end{theorem}

\begin{proof}
The type $\sfC$ formula 
is \cite[Cor.~3.12]{AHHP}, and the other formulas are derived by similar arguments.
In detail, let $v = w_\sfB(\lambda/\mu)=w_\sfC(\lambda/\mu)$ and $w=w_\sfD(\lambda/\mu)$.
Then \cite[Thm.~3.11]{AHHP} constructs a bijection 
\[\res : \SetShTab^+(\lambda/\mu)\to \left\{ (a,b) : a \in \cH^\sfC(v)\text{ and }b \in \cU^\sfC(a)\right\}\]
that is weight-preserving in the sense that $|T|=\ell(a)=\ell(b)$ and $\bfz^T = \bfz_{b}$ whenever $\res(T)=(a,b)$.
It is clear from the definition of this map in \cite{AHHP} that it restricts to a bijection
\[ \SetShTab(\lambda/\mu)\to \left\{ (a,b) : a \in \cH^\sfB(v)\text{ and } b \in \cU^\sfB(a)\right\}.\]
The existence of these bijections proves part (a) in view of Proposition~\ref{cFX-def2}.

Given an integer sequence $a=(a_1,a_2,\dots,a_k)$, let $\sigma(a) = (|a_1|-1, |a_2|-1,\dots,|a_k|-1)$.
By \cite[Prop.~9.8 and Cor.~10.2]{StembridgeSome},
the operation $\sigma$ defines a bijection from the set of reduced words for $w\in W^\sfD_\infty$ to the set of reduced words for $v \in W^\sfB_\infty$.
As $v $ and $w$ are both fully commutative, it follows that $\sigma$ extends to a bijection $\cH^\sfD(w) \to \cH^\sfB(v)$.
Since $b \in \cU^\sfD(a)$ for $a \in \cH^\sfD(w)$ if and only if $b \in \cU^\sfB(\sigma(a))$,
the map $(a,b) \mapsto \res^{-1}(\sigma(a), b)$ is a weight-preserving bijection
\[  \left\{ (a,b) : a \in \cH^\sfD(v)\text{ and } b \in \cU^\sfD(a)\right\} \to  \SetShTab(\lambda/\mu)\]
which suffices to prove part (b) by Proposition~\ref{cFX-def2}.
\end{proof}

Tamvakis derives another family of tableau formulas for $\cFX_w(\bfz)$ in \cite[Thms.~2 and 3]{Tamvakis}.
However, it 
does not seem straightforward to obtain the results here as special cases 
of \cite{Tamvakis}; see \cite[\S7.3]{AHHP}.

\subsection{Kirillov--Naruse polynomials}\label{knp-sect}

This section reviews the construction in \cite{KN} of certain elements  $\fkGX_w(\bfx;\bfy;\bfz)$
of the ring
\[
\SS_\beta[\bfx;\bfy] := \SS_\beta[x_1,x_2,x_3,\dots,y_1,y_2,y_3,\dots] \subset \ZZ[\beta]\llbracket \bfx;\bfy;\bfz\rrbracket .
\]
We may directly define these power series using the following formula, 
which involves Grothendieck polynomials for type $\sfA$, $K$-Stanley symmetric functions, and the Demazure product.

\begin{definition}[{\cite[Prop.~2]{KN}}]\label{cont-def}
Fix $X \in \{\sfB,\sfC,\sfD\}$ and let $\ell=\ell^X$. Then, for $w \in W^X_\infty$ define 
\[
\fkGX_w(\bfx;\bfy;\bfz) = \sum_{
(\sigma,u,\tau) 
} \beta^{\ell(\sigma)+\ell(u)+\ell(\tau) - \ell(w)} \cdot \fkG_{\sigma}(\bfy) \cdot \cFX_{u}(\bfz) \cdot \fkG_{\tau}(\bfx) \in \SS_\beta[\bfx;\bfy]
\]
where the sum is over all triples $(\sigma,u,\tau) \in \WA_\infty\times W^X_\infty \times \WA_\infty$
such that 
$\sigma^{-1} \circ u \circ \tau = w$.
\end{definition}

 We refer to $\fkGX_w(\bfx;\bfy;\bfz) $ as a \defn{Kirillov--Naruse polynomial} or \defn{double Grothendieck polynomial}.

\begin{remark}
The cited formula  \cite[Prop.~2]{KN} accidentally omits the power of $\beta$, which should be included to match
\cite[Def.~7]{KN}. Adding this factor
makes $\fkGX_w(\bfx;\bfy;\bfz)$ homogeneous of degree $\ell^X(w)$ under the usual convention that $\deg \beta=-1$ and $\deg x_i= \deg y_i = 1$.
As noted in \cite{KN}, this definition of $\fkGX_w(\bfx;\bfy;\bfz)$ is equivalent to but simpler than an earlier construction in \cite{Kirillov}.
\end{remark}

Setting $\beta=0$ transforms $\fkGX_w(\bfx;\bfy;\bfz) $ to the \defn{double Schubert polynomials} for classical groups introduced
in \cite{IkMiNa}; further setting $\bfy=0$ recovers the \defn{Schubert polynomials} constructed by Billey and Haiman in \cite{BilleyHaiman}.
Notice that we always have $\fkGX_w(0;0;\bfz) = \cFX_w(\bfz)$.

\begin{example}\label{cont-def-ex}
If $X \in \{\sfB,\sfC\}$ and $w= \overline{2}1= t_1t_0\in  \WB_\infty=\WC_\infty$ then the triples $(\sigma,u,\tau)$ indexing 
the sum that defines $\fkGX_w(\bfx;\bfy;\bfz)$ are given by 
$(1, t_1t_0,1)$, $(t_1, t_0,1)$, and $(t_1,t_1t_0,1)$.
Since $\fkG_1(\bfx) = 1$ and $\fkG_{t_1}(\bfx) =\fkG_{21}(\bfx) = x_1$
and since $t_0=\overline{1}$ and $t_1t_0 = \overline{2}1$ are both Grassmannian,
we have 
\[
\fkGB_{\overline 21}(\bfx;\bfy;\bfz) 
=  y_1 \GP_1(\bfz) +(1+ \beta y_1) \GP_2(\bfz)
\quand
\fkGC_{\overline 21}(\bfx;\bfy;\bfz) =  y_1 \GQ_1(\bfz) + (1+\beta y_1) \GQ_2(\bfz).
\]
\end{example}

Kirillov and Naruse show that the formal power series 
$\fkGX_w(\bfx;\bfy;\bfz) $ represent the opposite Schubert classes in the equivariant $K$-theory ring of $G/B$ in types $\sfB$, $\sfC$, and $\sfD$ \cite[\S6]{KN}.
In particular, the structure constants that expand the product
of $\fkGX_v(\bfx;\bfy;\bfz) $ and $\fkGX_w(\bfx;\bfy;\bfz) $ into $\fkGX$-functions are the same as the ones for the opposite Schubert basis of $K_T(G/B)$, in the following sense:

\begin{lemma}[{\cite[Cor.~2]{KN}}]
\label{structure-lem2}

Fix $X \in \{\sfB,\sfC,\sfD\}$
and $u,v \in W^X_\infty$. 
Then there are unique homogeneous coefficients
 $c_{u,v}^w(\bfy) \in \AA$ for each $w\in W^X_\infty$
with
\[\textstyle \fkGX_u(\bfx;\bfy;\bfz)\cdot\fkGX_v(\bfx;\bfy;\bfz) = \sum_{w\in W^X_\infty}c_{u,v}^{w}(\bfy) \cdot \fkGX_w(\bfx;\bfy;\bfz),\]
and for all $n\in \PP$ with $u,v\in W^X_n$ it holds that $\left\{ c_{u,v}^{w}(\bfy) : w \in W^X_n\right\} \subset \AA^{(n)}$ and
\[  \textstyle
 \cKX{u}\cdot \cKX{v} = \sum_{w\in W^X_n}c_{u,v}^{w}(y_i \mapsto 1-x^{\e_i})|_{\beta=-1} \cdot \cKX{w} \in K_T(G/B)  \] 
when $G$ is the rank $n$ classical group of type $X$.
\end{lemma}

This connection will let us derive a \defn{transition formula} for $\fkGX_w(\bfx;\bfy;\bfz) $ 
from the equivariant Chevalley formula in Section~\ref{oppo-sect}.
This application will use the calculations in the following example.

\begin{example}\label{x1BCD-eq}
If $\ell(w) = 1$ and $w \notin  \WA_\infty \cong S_\infty$ then $\fkGX_w(\bfx;\bfy;\bfz) = \cFX_w(\bfz)$,
while if $\ell(w) = 1$ and $w \in \WA_\infty$ then $w=t_k$ for some $i \in \PP$ and it is immediate from Definition~\ref{cont-def} that
\[
1+\beta\fkGX_w(\bfx;\bfy;\bfz) = \( 1+\beta\fkG_{s_k}(\bfy)\) \(1+\beta\cFX_{t_k}(\bfz)\) \(1+\beta\fkG_{s_k}(\bfx)\).
\]
Comparing this equation with Example~\ref{sk-ex} and \eqref{gr-eq} leads to the following identities:
\ben
\item[(a)] In type $X=\sfB$ we have 
$ 1+\beta \fkGB_{t_0}(\bfx;\bfy;\bfz) = 1 + \beta \GP_{1}(\bfz)$
along with
\[\ba 
 1+\beta \fkGB_{t_k}(\bfx;\bfy;\bfz) &= \textstyle (1 + \beta \GQ_{1}(\bfz)) \prod_{i\in[k]} (1+\beta x_i)(1+\beta y_i)
\ea\]
for $k \geq 1$.
Since 
\[ 1 + \beta \GQ_{1}(\bfz) = (1 + \beta \GP_{1}(\bfz))^2 = 1 + 2\beta \GP_{1}(\bfz)+\beta^2 \GP_{2}(\bfz),\]
it follows from Example~\ref{cont-def-ex} that
\[ 1 + \beta \GQ_{1}(\bfz)=
1 + \beta \cdot \tfrac{2+ \beta y_1}{1+\beta y_1}\cdot  \fkGB_{\overline 1}(\bfx;\bfy;\bfz) +\beta^{2} \cdot (1+\beta y_1)^{-1} \cdot \fkGB_{\overline 21}(\bfx;\bfy;\bfz).
\]
Thus, we have
\[\label{b-x1-eq}
1+\beta x_1 =  
\tfrac{1}{1+\beta y_1}\cdot \tfrac{1}{1+\beta \fkGB_{t_{0}}(\bfx;\bfy;\bfz)} \cdot \tfrac{1+\beta \fkGB_{t_1}(\bfx;\bfy;\bfz)}{1+\beta \fkGB_{t_{0}}(\bfx;\bfy;\bfz)}
=
 \tfrac{1+\beta \fkGB_{21}(\bfx;\bfy;\bfz)}{(1+\beta y_1) + \beta  (2+ \beta y_1)  \fkGB_{\overline 1}(\bfx;\bfy;\bfz) +\beta^{2}  \fkGB_{\overline 21}(\bfx;\bfy;\bfz)}
\]
and for $k\geq 2$  it holds that
$
1+\beta x_k = \tfrac{1}{1+\beta y_k} \cdot  \tfrac{1+\beta \fkGB_{t_k}(\bfx;\bfy;\bfz)}{1+\beta \fkGB_{t_{k-1}}(\bfx;\bfy;\bfz)}. 
$

\item[(b)] In type $X=\sfC$ we have for all $k \geq 0$ that
\[\textstyle 1+\beta \fkGC_{t_k}(\bfx;\bfy;\bfz) = (1 + \beta \GQ_{1}(\bfz)) \prod_{i\in[k]} (1+\beta x_i)(1+\beta y_i).\]
Hence if $k\geq 1$ then
$
1+\beta x_k = 
 \tfrac{1}{1+\beta y_k} \cdot \tfrac{1+\beta \fkGC_{t_k}(\bfx;\bfy;\bfz)}{1+\beta \fkGC_{t_{k-1}}(\bfx;\bfy;\bfz)}.
 $
 
\item[(c)] In type $X=\sfD$ we have
\[\ba
 1+\beta \fkGD_{t_{-1}}(\bfx;\bfy;\bfz) &= 1 + \beta \GP_{1}(\bfz), 
 \\
  1+\beta \fkGD_{t_{1}}(\bfx;\bfy;\bfz) &= \(1 + \beta \GP_{1}(\bfz)\)(1+\beta x_1)(1+\beta y_1), 
 \\
1+\beta \fkGD_{t_{k}}(\bfx;\bfy;\bfz)  &= \textstyle \(1+\beta \GQ_{1}(\bfz)\)\prod_{i\in[k]} (1+\beta x_i)(1+\beta y_i),
\ea
\]
for all $k\geq 2$.
Hence
\[
1+\beta x_1 = 
 \tfrac{1}{1+\beta y_1} \cdot \tfrac{1+\beta \fkGD_{t_1}(\bfx;\bfy;\bfz)}{1+\beta \fkGD_{t_{-1}}(\bfx;\bfy;\bfz)}
\]
and
\[
1+\beta x_2=
 \tfrac{1}{1+\beta y_2}\cdot \tfrac{1}{1+\beta \fkGD_{t_{-1}}(\bfx;\bfy;\bfz)} \cdot \tfrac{1+\beta \fkGD_{t_2}(\bfx;\bfy;\bfz)}{1+\beta \fkGD_{t_{1}}(\bfx;\bfy;\bfz)}
\]
and for $k\geq 3$ we have
\[1+\beta x_k =  \tfrac{1}{1+\beta y_k}\cdot\tfrac{1+\beta \fkGD_{t_k}(\bfx;\bfy;\bfz)}{1+\beta \fkGD_{t_{k-1}}(\bfx;\bfy;\bfz)}.\]

\een
\end{example}

For Grassmannian signed permutations, there are explicit tableau formulas \cite[Prop.~5]{KN} for $\fkGX_w(\bfx;\bfy;\bfz)$
refining the identities for $\cFX_w(\bfz)$ in \eqref{gr-eq}; see also \cite{HIMN}.
Tamvakis describes some related tableau formulas for a certain class of \defn{skew} signed permutations in \cite[\S3 and 4]{Tamvakis}, generalizing results of Matsumura \cite{Matsumura} in type $\sfA$.
 
 Outside of these cases, it is not particularly easy to calculate the power series $\fkGX_w(\bfx;\bfy;\bfz)$.
 The results in the next section will provide an inductive solution to this problem.

\subsection{Transitions for classical types}

Fix $X \in \{\sfB,\sfC,\sfD\}$ for the rest of this section.
Define $\bVV^X $ 
to be the $\AA$-module consisting of all formal $\AA$-linear combinations 
$
\sum_{w \in W^X_\infty} f_w(\bfy) \cdot \fkG_w^X(\bfx;\bfy;\bfz)
$
such for each $k \in \NN$  the set
\[\left\{ w \in W^X_\infty : \ell^X(w) = k\text{ and }f_w(\bfy)\neq 0\right\}\]
is finite.
We say that a $\ZZ[\beta]$- or $\AA$-linear map $\Phi : \bVV^X \to\bVV^X$
is \defn{continuous} if  
\[
\Phi\Bigl(\textstyle \sum_{w} f_w(\bfy) \cdot \fkG_w^X(\bfx;\bfy;\bfz)\Bigr) = \sum_{w} \Phi\Bigl(f_w(\bfy) \cdot \fkG_w^X(\bfx;\bfy;\bfz)\Bigr).
\]
We may view $\bVV^X$ as an $\AA$-submodule of $\ZZ[\beta]\llbracket \bfx;\bfy;\bfz\rrbracket$ 
since $\AA \subset \ZZ[\beta]\llbracket \bfx;\bfy\rrbracket$
and since if we regard $\beta$ as having degree zero and all other variables as having degree one,
then the lowest degree monomials appearing in $\fkG_w^X$  all have degree at least $\ell^X(w)$.

For any parameter $a$ let $\ominus a = \frac{-a}{1+\beta a}$. Notice that $\ominus(\ominus a) = a$.
For $i \in \PP$ define 
\be \textstyle y_{-i} =\ominus y_i =  \frac{-y_i}{1+\beta y_i}
\quad\text{so that}
\quad
\frac{1}{1+\beta y_{-i}} = 1+\beta y_i.
\ee
Using this convention, we let $w \in \WB_\infty=\WC_\infty$ act on $f(\bfy) \in \AA$
by the formula \eqref{ast-eq}, so that
\be
w\ast f(\bfy)  = f(y_{w(1)},y_{w(2)},y_{w(3)},\dots).
\ee
We now define several continuous $\ZZ[\beta]$-linear operators on $\bVV^X$,
to be denoted $\bftX_{ij}$, $\bfuX_{ij}$, $\bfvX_k$, $\bfoX_k$, and $\bfMX_k$.
Fix integers $ i < j$  and $k$ and $j,k\in\PP$.
Recall from \eqref{tij-eq} that $t_{ij}\in\WB_\infty=\WC_\infty$
 is the identity element when $i=-j$ and  a reflection otherwise.
We have $t_{ij} \in \WA_\infty$ if and only if $0<i<j$ and $t_{ij} \in \WD_\infty$ if and only if $i\neq 0$.
Let 
$\bftX_{ij}$ be the continuous $\AA$-linear operator sending
\be
\bftX_{ij}: 
\fkGX_w(\bfx; \bfy;\bfz)\mapsto
\begin{cases}
\fkGX_{wt_{ij}}(\bfx; \bfy;\bfz) &\text{if $t_{ij} \in W^X_\infty$ and $\ell^X(w) +1 = \ell^X(wt_{ij})$}
\\
0&\text{otherwise} 
\end{cases}
\ee
for each $w \in W^X_\infty$. Notice that we have 
\[
\bftD_{0j}=0\quand \bftB_{-j,j}=\bftC_{-j,j}=\bftD_{-j,j}=0.
\]
Then let 
$\bfuX_{ij}$ be the continuous $\ZZ[\beta]$-linear operator sending
\be
\bfuX_{ij}: 
f(\bfy)\cdot \fkGX_w(\bfx; \bfy;\bfz)\mapsto
\(wt_{ij}w^{-1}\ast f(\bfy)\)\cdot \bftX_{ij}\fkGX_{w}(\bfx; \bfy;\bfz)
\ee
for each $f(\bfy) \in \AA$ and $w \in W^X_\infty$. 
We likewise have 
\[
\bfuD_{0j}=0\quand \bfuB_{-j,j}=\bfuC_{-j,j}=\bfuD_{-j,j}=0.
\]
Next 
 let $\bfv_{k} :\bVV^X\to\bVV^X$ be the continuous $\AA$-linear map with
\be
\bfvX_{k}:  \fkGX_w(\bfx; \bfy;\bfz)\mapsto
\tfrac{1}{1+\beta y_{w(k)}}\cdot \fkGX_w(\bfx; \bfy; \bfz).
\ee
Now define $\bfoB_{k} : \bVV^\sfB\to\bVV^\sfB$ to be the continuous $\ZZ[\beta]$-linear operator sending
 \be
\bfoB_{k} : f(\bfy)\cdot \fkGB_w(\bfx; \bfy;\bfz)\mapsto
\begin{cases}
 f(0) \cdot \fkGB_{w(-k,k)}(\bfx;\bfy;\bfz) &\text{if }\ellB(w(-k,k))=\ellB(w) + 1\\
0 & \text{otherwise}
\end{cases}
\ee
for each $f(\bfy) \in \AA$ and $w \in W^\sfB_\infty$. Also let $\bfoC_k = \bfoD_k = 0$.
We finally set
\be
\bfM^X_{k} =  
\prod_{k<l<\infty}^{[<]}(1 +\beta \bftX_{kl}) \cdot (1-\beta \bfoX_k)\cdot \prod_{-\infty<j<k}^{[<]}(1 -\beta \bfuX_{jk}) \cdot \bfvX_{k}
\ee
This infinite product is a well-defined continuous $\ZZ[\beta]$-operator $ \bVV^X\to\bVV^X$  since 
$1 +\beta \bftX_{kl}$ and $1 +\beta \bfuX_{jk}$ act as the identity on  any fixed term $f(\bfy) \cdot \fkGX_w(\bfx;\bfy;\bfz)$ if $|j|$ and $l$ are sufficiently large. 
Applying $\bfMX_k$ to $\fkGX_w(\bfx;\bfy;\bfz)$ 
can yield an infinite sum of Kirillov--Naruse polynomials, however.

\begin{theorem}\label{bcd-thm1}
Fix $X \in \{\sfB,\sfC,\sfD\}$.
If $u \in W^X_\infty$ and $k \in \PP$ then
\[
(1+\beta x_k)  \fkG_u^X(\bfx;\bfy;\bfz)
=
\bfM^X_k \fkG_u^X(\bfx;\bfy;\bfz)
.\]
\end{theorem}

\begin{proof}
After comparing the formulas for $[\cL_{\e_k}]$ and $1+\beta x_k$ in Proposition~\ref{cL-prop} and Example~\ref{x1BCD-eq},
this follows from Theorem~\ref{lp-cor2} 
using 
Lemma~\ref{structure-lem2} 
by the same argument as in the proof of Theorem~\ref{len-cor1}.
In fact, the argument here is simpler since the representation ring $\ZZ[T]$ of the rank $n$ torus in type $X$
is equal to the Laurent polynomial ring
 $\ZZ[x^{\pm \e_1},\dots,x^{\pm \e_n}]$ rather than to a quotient.
\end{proof}

\begin{example}
When $X=\sfB$, $k=1$, and $w \in \{1, \overline{1}, \overline{2}1\}$,
one can use Theorem~\ref{bcd-thm1} to check that
\[\ba
 (1+\beta x_1) \fkGB_1  &=
(1+\beta y_{1})^{-1} \cdot \(\fkGB_1 +\beta\cdot \fkGB_{21}\) 
 -\beta\cdot (2+\beta y_{1}) \cdot \(\fkGB_{\overline{1}} +\beta\cdot \fkGB_{2\overline{1}}\) 
 \\&\quad +
  \beta^2\cdot (1+\beta y_{2}) \cdot \(\fkGB_{\overline{2}1} + \beta\cdot \fkGB_{1\overline{2}}\)
 - \beta^3 \cdot \Psi,
\\[-10pt]\\
(1+\beta x_1) \fkGB_{\overline 1}
&= (1+\beta y_{1}) \cdot \(\fkGB_{\overline{1}} +\beta\cdot \fkGB_{2\overline{1}}\)- \beta \cdot (1+\beta y_{2}) \cdot \(\fkGB_{\overline{2}1} + \beta\cdot \fkGB_{1\overline{2}}\)+\beta^2 \cdot \Psi,
\\[-10pt]\\
(1+\beta x_1) \fkGB_{\overline 21}
&=(1+\beta y_{2}) \cdot \(\fkGB_{\overline{2}1} + \beta\cdot \fkGB_{1\overline{2}}\) - \beta \cdot \Psi,
\ea
\]
where
\[
\Psi := \textstyle\sum_{n=3}^\infty (-\beta)^{n-3} \cdot (1+\beta y_{n}) \cdot \(\fkGB_{\overline{n}12\cdots(n-1)} + \beta\cdot \fkGB_{1\overline{n}2\cdots(n-1)}\).
\] As $\fkGB_1=1$, one can combine these identities to deduce that
\[
(1+\beta x_1)\cdot \((1+\beta y_1) + \beta\cdot  (2+ \beta y_1)\cdot  \fkGB_{\overline 1} +\beta^{2}  \cdot\fkGB_{\overline 21} \) = 1+\beta \cdot \fkGB_{21}, \]
recovering the formula in \eqref{b-x1-eq}.
\end{example}

\begin{example} For each type $X \in \{\sfB,\sfC,\sfD\}$, we can use Theorem~\ref{bcd-thm1} to compute that
\[\ba
(1+\beta x_3)  \fkGX_{\overline{6} \overline{1}3\overline{4}\overline{2}5} &=
                        (1 + \beta y_{3})^{-1} \cdot \fkGX_{\overline{6} \overline{1} 3 \overline{4} \overline{2} 5}
&&        + \beta \cdot (1 + \beta y_{3})^{-1} \cdot \fkGX_{\overline{6} \overline{1} 5 \overline{4} \overline{2} 3}
\\&\quad        -\beta\cdot (1 + \beta y_{2})^{-1} \cdot \fkGX_{\overline{6} \overline{1} 2 \overline{4} \overline{3} 5}
&&        -\beta^2\cdot (1 + \beta y_{2})^{-1} \cdot \fkGX_{\overline{6} \overline{1} 5 \overline{4} \overline{3} 2}
\\&\quad        -\beta\cdot (1 + \beta y_{1})^{-1} \cdot \fkGX_{\overline{6} \overline{3} 1 \overline{4} \overline{2} 5}
&&        -\beta^2\cdot (1 + \beta y_{1})^{-1} \cdot \fkGX_{\overline{6} \overline{3} 5 \overline{4} \overline{2} 1}
\\&\quad        -\beta \cdot(1 + \beta y_{1}) \cdot \fkGX_{\overline{6} 3 \overline{1} \overline{4} \overline{2} 5}
&&        -\beta^2\cdot (1 + \beta y_{1}) \cdot \fkGX_{\overline{6} 3 5 \overline{4} \overline{2} \overline{1}}
\\&\quad        + \beta^2\cdot (1 + \beta y_{3}) \cdot \fkGX_{\overline{6} 1 \overline{3} \overline{4} \overline{2} 5}
&&        + \beta^3\cdot (1 + \beta y_{3}) \cdot \fkGX_{\overline{6} 1 \overline{2} \overline{4} \overline{3} 5}
\\&\quad        -\beta^3\cdot (1 + \beta y_{5}) \cdot \fkGX_{\overline{6} 1 \overline{5} \overline{4} \overline{2} 3}
&&        -\beta^4\cdot (1 + \beta y_{5}) \cdot \fkGX_{\overline{6} 1 \overline{4} \overline{5} \overline{2} 3}
\\&\quad        + \beta^2 \cdot  (1 + \beta y_{6}) \cdot \fkGX_{\overline{1} 3 \overline{6} \overline{4} \overline{2} 5}
&&        + \beta^3 \cdot (1 + \beta y_{6}) \cdot \fkGX_{\overline{1} 3 \overline{4} \overline{6} \overline{2} 5}
\\&\quad        + \beta^2 \cdot \Delta^X
&&        - \beta^3 \cdot\Sigma^X
\ea
\]
where we set 
\[\Delta^X := 
\begin{cases} \fkGB_{\overline{6} \overline{3} \overline{1} \overline{4} \overline{2} 5}
  +
   \beta\cdot \fkGB_{\overline{6} \overline{3} 5 \overline{4} \overline{2} \overline{1}} &\text{if }X=\sfB \\
   0&\text{if }X \in \{\sfC,\sfD\}\end{cases} 
\]
and
\[\textstyle
 \Sigma^X := 
 \sum_{n=7}^\infty (-\beta)^{n-7} \cdot (1+\beta y_n)\cdot \(
 \fkGX_{u(n)}
 +
 \beta \cdot  \fkGX_{v(n)}\)
\]
for the signed permutations  
$u(n) := \overline{1}3\overline{n}\overline{4}\overline{2}56\cdots (n-1)
$ and $ v(n) := \overline{1}3\overline{4}\overline{n}\overline{2}56\cdots (n-1).$
In type $\sfC$ this computation generalizes \cite[Eq.~(1.1)]{Billey}.
\end{example}

Only in type $\sfB$, we additionally define 
\be
\bfnB_k = \bftB_{0k} \cdot \bfvB_k
\ee
 to be the continuous $\AA$-linear map $\bVV^\sfB\to\bVV^\sfB$
sending
\[
\bfnB_{k} : \fkG^\sfB_w(\bfx; \bfy; \bfz)\mapsto
\begin{cases}
\tfrac{1}{1+\beta y_{w(k)}}\fkG^\sfB_{w(-k,k)}(\bfx; \bfy; \bfz) &\text{if }\ell^\sfB(w) +1 = \ell^\sfB(w(-k,k))
\\
0&\text{otherwise} .
\end{cases}
\]
Notice that $\ell^\sfB(w) +1 = \ell^\sfB(w(-k,k))$ can only hold if $w(k)>0$.
Now set
\be
\ba
\bfR^\sfB_{k} &=  
\prod_{k>j>-\infty}^{[>]} (1+\beta \bftB_{jk}) \cdot (1+ \beta\bfnB_{k}),
\\
\bfR^\sfC_{k}  &=\prod_{k>j>-\infty }^{[>]} (1+\beta \bftC_{jk}), 
\\
\bfR^\sfD_{k} &= \prod_{\substack{k>j>-\infty \\ j \neq 0} }^{[>]}(1+ \beta\bftD_{ik}).
\ea
\ee
Recall that $\bftX_{-k,k}=0$ in all types.
Alternatively, for $X\in\{\sfB,\sfC,\sfD\}$ we could define 
\be
\bfR^X_{k} =  
\prod_{k>j>-\infty}^{[>]} (1+\beta \bft^X_{jk}) \cdot (1+ \beta\bfnX_{k})
\ee
under the convention that $\bftD_{0k}= \bfnC_{k}=\bfnD_{k}=0$.
As with $\bfMX_k$, this infinite product 
is a well-defined continuous $\ZZ[\beta]$-linear operator $\bVV^X\to\bVV^X$. Define
\[ \VV^X = \AA\spanning\left\{ \fkGX_w(\bfx;\bfy;\bfz)  : w \in W^X_\infty\right\}\subset \bVV^X.\]
Unlike $\bfMX_k$, the following holds:

\begin{proposition}
The operator $\bfR^X_{k}$ restricts to an $\AA$-linear map $\VV^X\to\VV^X$.
\end{proposition}

\begin{proof}
For each $w \in W^X_\infty$, one has 
$\bfR^X_{k} \fkGX_w(\bfx;\bfy;\bfz) = \prod_{k>j>-N}^{[>]} (1+\beta \bft^X_{jk}) \cdot (1+ \beta\bfnX_{k})\cdot \fkGX_w(\bfx;\bfy;\bfz) $
for all sufficiently large $N >0$ by Lemmas~\ref{len-lem1} and \ref{len-lem3}, and the latter expression is clearly in
$\VV^X$.
\end{proof}

In the following lemmas, let $\UU^X = \ZZ[\beta]\spanning\left\{ \fkG_w^X(\bfx;\bfy;\bfz) : w \in W^X_\infty\right\}\subset \VV^X$.

\begin{lemma}\label{weig-lem2}
Fix $X \in \{\sfB,\sfC,\sfD\}$ and $k \in \PP$. Suppose $F \in\UU^X $. Then 
\be\label{pi-pi-eq} \( \prod_{-\infty<j<k}^{[<]}(1 -\beta \bfuX_{jk}) \cdot \bfvX_{k}
\cdot \prod_{k>j>-\infty}^{[>]} (1+\beta \bft^X_{jk})\) F =   \bfvX_{k} F.\ee
\end{lemma}

\begin{proof}
As in the proof of Lemma~\ref{weig-lem}, it is straightforward to check that
\be\label{pi-pi-eq2}
 ( (1-\beta \bfuX_{jk})\cdot \bfv_{k} \cdot (1+\beta \bftX_{jk}))\fkGX_w(\bfx;\bfy;\bfz) =  
 \bfvX_k \fkGX_w(\bfx;\bfy;\bfz)
 \ee
 for any fixed $w \in W^X_\infty$ and $j\in \ZZ$ with $j<k$. 
 Since $\bftX_{jk}$, $\bfuX_{jk}$, and $\bfvX_k$ are all $\ZZ[\beta]$-linear operators,
 and since $1+\beta \bftX_{jk}$ acts as the identity map on $F$ if $|j|$ is sufficiently large,
we deduce from \eqref{pi-pi-eq2} by induction that the right hand side of \eqref{pi-pi-eq}
 is equal to 
 \be\label{pi-pi-eq3}\textstyle
 \( \prod_{-\infty<j<k-N}^{[<]}(1 -\beta \bfuX_{jk}) \cdot \bfvX_{k}
\cdot \prod_{k-N>j>-\infty}^{[>]} (1+\beta \bft^X_{jk})\) F
\ee
for all $N \in \NN$. 
If  $N$ is sufficiently large then every operator $1+\beta \bft^X_{jk}$ with $j<k-N$ 
 fixes $F$ and every operator $1 -\beta \bfuX_{jk}$ with $j<k-N$ fixes $\bfvX_k F$,
and then
 \eqref{pi-pi-eq3} gives $\bfvX_kF$ as needed.
\end{proof}

\begin{lemma}\label{weig-lem3}
Fix $X \in \{\sfB,\sfC,\sfD\}$ and $k \in \PP$. Suppose $F \in\UU^X $. Then 
\be\label{weig-lem3-eq}
(1+\beta x_k) \bfR^X_k F = \(\prod_{k<l<\infty}^{[<]}(1 +\beta \bftX_{kl})\cdot \bfvX_{k}\)F.
\ee
\end{lemma}

\begin{proof}
First assume $X \in \{\sfC,\sfD\}$. Then 
 $\bfRX_k$ preserves $\UU^X$, so Theorem~\ref{bcd-thm1} implies that 
\[
(1+\beta x_k) \bfRX_k F = \bfMX_k \bfRX_k F
=
 \(\prod_{k<l<\infty}^{[<]}(1 +\beta \bftX_{kl}) \cdot \prod_{-\infty<j<k}^{[<]}(1 -\beta \bfuX_{jk}) \cdot \bfvX_{k}
\cdot \prod_{k>j>-\infty}^{[>]} (1+\beta \bft^X_{jk})\) F
\]
since $\bfMX_k$ is a $\ZZ[\beta]$-linear map, and the desired identity follows from Lemma~\ref{weig-lem2}.

We must be more careful to handle type $X=\sfB$,
as the operator $\bfRB_k$ does not preserve $\UU^\sfB$.
Define 
$
\bfQB_k = \prod_{k>j>-\infty}^{[>]} (1+\beta \bft^X_{jk})
$ so that 
$
 \bfRB_k = \bfQB_k \cdot (1+\beta \bfnB_k).
$
Notice that $\bfQB_k$ does preserve $\UU^\sfB$.

Now suppose $F = \fkGB_w(\bfx;\bfy;\bfz)$ for some $w \in W^\sfB_\infty$.
By linearity, it suffices to prove the lemma in this special case.
If $\ell(w(-k,k))\neq \ell(w)+1$,
then $1+\beta \bfnB_k$ fixes $F$ so
 Theorem~\ref{bcd-thm1} implies that
\[
\ba
(1+\beta x_k) \bfRB_k F &
=(1+\beta x_k) \bfQB_k F 
= \bfMB_k \bfQB_k F
\\&=
 \(\prod_{k<l<\infty}^{[<]}(1 +\beta \bftB_{kl}) \cdot (1-\beta \bfoB_k)\cdot \prod_{-\infty<j<k}^{[<]}(1 -\beta \bfuB_{jk}) \cdot \bfvB_{k}
\cdot \prod_{k>j>-\infty}^{[>]} (1+\beta \bftB_{jk})\) F.
\ea
\]
Using Lemma~\ref{weig-lem2}, the last expression simplifies to 
$
  \(\prod_{k<l<\infty}^{[<]}(1 +\beta \bftB_{kl}) \cdot (1-\beta \bfoB_k)\cdot \bfvB_{k}\)F
  $,
  which is equal to the left hand side of \eqref{weig-lem3-eq} since 
$(1-\beta \bfoB_k)\cdot \bfvB_{k}F = \bfvB_{k}F$ if $\ell(w(-k,k))\neq \ell(w)+1$.

Now instead suppose that $\ell(w(-k,k))= \ell(w)+1$.
Let $G =  \fkGB_{w(-k,k)}(\bfx;\bfy;\bfz)$. Then since $\bfQB_k$ is $\AA$-linear, Theorem~\ref{bcd-thm1} implies that
\be\label{almost-eq}
(1+\beta x_k) \bfRB_k F 
=(1+\beta x_k) \bfQB_k \(F + \tfrac{\beta}{1+\beta y_{w(k)}} G\) 
 = \bfMB_k \bfQB_k F + \tfrac{\beta}{1+\beta y_{w(k)}} \bfMB_k\bfQB_k G.
\ee
In this setup we have
\[
\ba
 (1-\beta \bfoB_k)   \bfvB_{k} F &= \bfvB_k F - \beta \bfoB_k\(\tfrac{1}{1+\beta y_{w(k)}} \fkGB_{w}(\bfx;\bfy;\bfz)\) = \bfvB_kF - \beta G\quand \\
 (1-\beta \bfoB_k)   \bfvB_{k} G &= \tfrac{1}{1+\beta y_{w(-k)}} G= \tfrac{1}{1+\beta y_{-w(k)}} G = (1+\beta y_{w(k)}) G,
 \ea
 \]
 so
 it follows from Lemma~\ref{weig-lem2} that
\[
\ba
\bfMB_k \bfQB_k F 
&=\textstyle
 \(\prod_{k<l<\infty}^{[<]}(1 +\beta \bftB_{kl}) \cdot (1-\beta \bfoB_k)\cdot  \bfvB_{k}\) F
=
  \(\prod_{k<l<\infty}^{[<]}(1 +\beta \bftB_{kl})\)\(\bfvB_k F - \beta G\)
  \ea
\]
and
\[
\ba
 \bfMB_k \bfQB_k G 
&=
\textstyle
\(\prod_{k<l<\infty}^{[<]}(1 +\beta \bftB_{kl}) \cdot (1-\beta \bfoB_k)\cdot  \bfvB_{k}\) G
=
 \(\prod_{k<l<\infty}^{[<]}(1 +\beta \bftB_{kl})\)  \((1+\beta y_{w(k)}) G\).
 \ea
\]
As each operator $1+\beta \bftB_{kl}$ is $\AA$-linear,
  substituting these identities into \eqref{almost-eq}
gives \eqref{weig-lem3-eq}.
 \end{proof}

 Recall that a \defn{descent} of a signed permutation $w$ is a positive integer $i$ with $w(i)>w(i+1)$.

\begin{theorem}\label{bcd-thm2}
Fix $X \in \{\sfB,\sfC,\sfD\}$. Suppose $w \in W^X_\infty$ has at least one descent.
Let $a \in [1,\infty)$ be maximal with $w(a) > w(a+1)$, define $b \in [a+1,\infty)$ to be maximal
with $w(a)>w(b)$, and set $v = wt_{ab}$ and $c=v(a) = w(b)$. Then
\[
\fkG^X_w(\bfx; \bfy; \bfz) =\beta^{-1}\Bigl(   (1+\beta y_{c})(1+\beta x_a) \bfRX_a \fkGX_v(\bfx; \bfy;\bfz) -\fkGX_v(\bfx;\bfy;\bfz)\Bigr) 
\]
where for positive integers $c<0$ we set $y_{c} := \ominus y_{-c} = \frac{-y_{-c}}{1+\beta y_{-c}}$.
\end{theorem}

\begin{proof}
The structure of the proof is the same as for Theorem~\ref{a-thm2}.
The maximality of $a$ and $b$ implies that if $a<i<\infty$ then $\ell^X(v(a,i)) = \ell^X(v)+1$ if and only if $i=b$,
and that $\ell^X(w(a,i))\neq \ell^X(w)+1$ for all $a<i<b$.
Hence, by Lemma~\ref{weig-lem3} we have 
\[\ba
(1+\beta x_a) \bfRX_a \fkGX_v(\bfx; \bfy;\bfz) &= \textstyle\(\prod_{a<i<\infty}^{[<]}(1 +\beta \bft_{ai})\cdot \bfv_a\) \fkGX_v(\bfx; \bfy;\bfz)
\\& = \tfrac{1}{1+\beta y_c} \( \fkGX_v(\bfx;\bfy;\bfz) + \beta \fkG_w(\bfx;\bfy;\bfz)\).
\ea
\]
Substituting this identity into the right hand side of the desired formula gives $\fkGX_w(\bfx;\bfy;\bfx)$.
\end{proof}

Theorem~\ref{bcd-thm2} is a $K$-theoretic extension of \cite[Prop.~6.12]{IkMiNa}, which is itself an equivariant generalization of \cite[Thm.~4]{Billey}.
We can  recursively use this theorem to expand $\fkG^X_w(\bfx; \bfy; \bfz)$ as an
 $\AA$-linear combination of Kirillov--Naruse polynomials
indexed by Grassmannian signed permutations (which have the explicit formulas given in \cite[Prop.~5]{KN}).
The proof of Theorem~\ref{k-stanley-positivity-thm} in the next section will show that this recurrence always terminates in a finite number of iterations.
 
\subsection{Applications to K-Stanley positivity}\label{k-positivity-sect}

Fix integers $j<k$ with $k>0$. Given a signed permutation $w \in W^X_\infty$, we let 
\[
T^X_{jk} \cFX_w(\bfz) = \begin{cases}
\cFX_{wt_{jk}}(\bfz) &\text{if $t_{ij} \in W^X_\infty$ and $\ell^X(w)+1=\ell^X(wt_{jk})$} \\
0&\text{otherwise}
\end{cases}
\]
and extend by $\ZZ[\beta]$-linearity.
Unlike $\bft^X_{jk}$, our  definition of $T^X_{jk}$ in this context only makes sense as a formal symbolic operator since
the $K$-Stanley symmetric functions $\cFX_w(\bfz)$ are not linearly independent.
We compose these operators to define
\be
\ba
R^\sfB_{k} &=  
\prod_{k>j>-\infty}^{[>]} (1+\beta T^\sfB_{jk}) \cdot (1+ \beta T^\sfB_{0k}),
\\
R^\sfC_{k}  &=\prod_{k>j>-\infty }^{[>]} (1+\beta T^\sfC_{jk}),
\\
R^\sfD_{k} &= \prod_{\substack{k>j>-\infty \\ j \neq 0} }^{[>]}(1+ \beta T^\sfD_{ik}).
\ea
\ee
We   interpret $R^X_k$ as a formal symbolic operator on the module of finite $\ZZ[\beta]$-linear combinations of $K$-Stanley symmetric functions. 

\begin{corollary}\label{bcd-thm3}
Fix $X \in \{\sfB,\sfC,\sfD\}$. Suppose $w \in W^X_\infty$ has at least one descent.
Let $a \in [1,\infty)$ be maximal with $w(a) > w(a+1)$, define $b \in [a+1,\infty)$ to be maximal
with $w(a)>w(b)$. Then
\[
\cFX_w(\bfz) = \beta^{-1} (R^X_{a}- 1) \cFX_{wt_{ab}}(\bfz) \in \NN[\beta]\spanning\left\{ \cFX_u(\bfz) : u \in W^X_\infty\right\}.
\]
\end{corollary}

\begin{proof}
This follows immediately from Theorem~\ref{bcd-thm2} on setting $\bfx=\bfy=0$.
\end{proof}

\begin{example} Applying Corollary~\ref{bcd-thm3} to $w =\overline{3}4\overline{1}52 \in \WB_\infty=\WC_\infty$  gives 
\[
\cFX_{\overline{3}4\overline{1}52}(\bfz) = 
\cFX_{\overline{3} 4 2 \overline{1}}(\bfz)
 + \cFX_{\overline{3} 4 \overline{2} 1}(\bfz)
  + \beta \cdot \cFX_{\overline{3} 4 \overline{2} \overline{1}}(\bfz)
   + \beta \cdot \cFX_{\overline{3} 4 1 \overline{2}}(\bfz)
    + \beta^2 \cdot \cFX_{\overline{3} 4 \overline{1} \overline{2}}(\bfz)
\]
for both types $X \in \{\sfB,\sfC\}$. We can  recursively use Corollary~\ref{bcd-thm3} to expand each term
on the right, until we reach a linear combination of Kirillov--Naruse polynomials
indexed by signed permutations with no descents. This 
gives the distinct formulas
 \[\ba
\cFB_{\overline{3}4\overline{1}52}(\bfz)&=
4 \cdot \cFB_{\overline{4} \overline{2} \overline{1} 3}(\bfz)
 + 2 \cdot \cFB_{\overline{4} \overline{3} 1 2}(\bfz)
  + 2 \cdot \cFB_{\overline{5} \overline{2} 1 3 4}(\bfz)
\\&\quad 
 + 5\beta \cdot \cFB_{\overline{4} \overline{3} \overline{1} 2}(\bfz)
 + 5\beta \cdot \cFB_{\overline{5} \overline{2} \overline{1} 3 4}(\bfz)
 + 3\beta \cdot \cFB_{\overline{5} \overline{3} 1 2 4}(\bfz)
+ 6\beta^2 \cdot \cFB_{\overline{5} \overline{3} \overline{1} 2 4}(\bfz)
 \ea
\]
and
 \[
 \ba
 \cFC_{\overline{3}4\overline{1}52}(\bfz) & 
 =
 2 \cdot \cFC_{\overline{4} \overline{2} \overline{1} 3}(\bfz)
 + 2 \cdot \cFC_{\overline{4} \overline{3} 1 2}(\bfz)
  + 2 \cdot \cFC_{\overline{5} \overline{2} 1 3 4}(\bfz)
\\&\quad 
 + 3\beta \cdot \cFC_{\overline{4} \overline{3} \overline{1} 2}(\bfz)
 + 3\beta \cdot \cFC_{\overline{5} \overline{2} \overline{1} 3 4}(\bfz)
 + 3\beta \cdot \cFC_{\overline{5} \overline{3} 1 2 4}(\bfz)
 + 4\beta^2 \cdot \cFC_{\overline{5} \overline{3} \overline{1} 2 4}(\bfz).
 \ea
 \]
 Recall from Example~\ref{fc-ex}
 that $\overline{3}4\overline{1}52 = w_{B}(\lambda/\mu)=w_{C}(\lambda/\mu)$ for
$\lambda=(5,3,1)$ and $\mu=(2)$.
Therefore, via \eqref{gr-eq} and Theorem~\ref{fc-prop}, the preceding identities give a $\NN[\beta]$-linear $\GP$-expansion
of $\GP_{(5,3,1)/(2)}(\bfz)$ and a $\NN[\beta]$-linear $\GQ$-expansion of $\GQ_{(5,3,1)/(2)}(\bfz)$.
 \end{example}

Recall from \eqref{LD-eq} our definition of the descent set $\Des(w)$ and the least descent $\LD(w) \in \NN$
for $w \in W^\sfB_\infty$.
Following \cite[\S4]{Billey}. let $\prec_\LD$ be the partial order on $W^\sfB_\infty$ that has 
$v \prec_\LD w$ if and only if 
\be
\LD(v) < \LD(w)\quord 0<\LD(v) = \LD(w)\text{ and }v(\LD(v)) < w(\LD(w)).
\ee
 The minimal elements under this partial order and the signed permutations that are \defn{Grassmannian}
 in the sense of having no descents.
The following technical lemma is a slight generalization of the main claim in the proof of \cite[Thm.~4]{Billey}.

\begin{lemma}\label{tech-lem2}
Fix $X \in \{\sfB,\sfC,\sfD\}$ and $n \in \PP$. Choose $v \in W^X_n$
and $j_1,j_2 , \dots , j_l, k\in\ZZ $ with $0< k = \LD(v)$ such that $t_{j_ak} \in W^X_\infty$ and $\ell^X(vt_{j_1k}t_{j_2k}\cdots t_{j_ak}) =\ell^X(v)+a$ for all $a \in [l]$.
Assume
\[j_1<j_2<\dots<j_l<k
\quord \text{$j_1= 0$ and $j_2<\dots<j_l<k$ when $X=\sfB$},\]
where we fix $l\geq1$.
Then  the element $u := vt_{j_1k}t_{j_2k}\cdots t_{j_ak} \in W^X_\infty$ satisfies either
\ben
\item[(a)] $u \in W^X_n$ and $u\prec_\LD v$, or
\item[(b)] $u \in W^X_{n+1}$ and $\LD(u) < \LD(v)$.
\een
\end{lemma}

\begin{proof}
The first three parts of Lemma~\ref{tech-lem} imply that $u \in W^X_{n+1}$ along with $u(k) < v(k)$ and 
$\Des(u) \subseteq \Des(v)\cup\{k-1\}$, so 
$u\prec_\LD v$.
If $u \notin W^X_n$ then the last two parts of Lemma~\ref{tech-lem}
imply that we must have $vt_{j_1k}t_{j_2k}\cdots t_{j_{a}k} \in W^X_n$ for all $a \in [l-1]$ and $\LD(u) \leq k-1 <\LD(v)$.
\end{proof}

This leads to our main application of Theorem~\ref{bcd-thm3}.
 Consider the additive semigroups
\[ \GPP:=\NN[\beta]\spanning\left\{\GP_\lambda(\bfz) : \lambda\text{ strict}\right\}
\quand
\GQP:=\NN[\beta]\spanning\left\{\GQ_\lambda(\bfz) : \lambda\text{ strict}\right\}
\]
Neither of these sets is contained in the other; for a characterization of which $\GQ$-functions are in $\GPP$, see \cite[Cor.~1.2]{ChiuMarberg}.
The following theorem was conjectured in \cite[Rem.~3]{KN} and \cite[Conj.~7.2]{AHHP}.

\begin{theorem}\label{k-stanley-positivity-thm}
Fix $X \in \{\sfB,\sfC,\sfD\}$, $n \in \PP$, and $w \in W^X_n$ 
with $k=\LD(w)$. Then 
\be\label{ksp-eq} \cF^X_w(\bfz)  \in \NN\spanning\Bigl\{ \beta^{\ell^X(u)-\ell^X(w)}\cF^X_u(\bfz) : u \in W^X_{n+k}\text{ is Grassmannian}\Bigr\}.\ee
In particular, it holds that $\cFX_w(\bfz)  \in \GPP$ if $X \in \{\sfB,\sfD\}$ and $\cFX_w(\bfz)  \in \GQP$ if $X=\sfC$.
\end{theorem}

For the type $\sfA$ versions of Corollary~\ref{bcd-thm3} and Theorem~\ref{k-stanley-positivity-thm}, see \cite[Thm.~4]{LascouxTransitions}.

\begin{proof}
It suffices by \eqref{gr-eq} just to prove \eqref{ksp-eq}.
The desired containment is trivial if $k=0$ so assume $k>0$.
Consider the order $\prec_\LD$ restricted to the set $U  := \bigcup_{i=0}^k \{ u \in W^X_{n+i} : \LD(u) \leq k - i\}$.
Since every minimal element in this restricted partial order is Grassmannian, we may assume by induction that \eqref{ksp-eq} holds when $w$ is replaced by 
any $u \in U $ with $u\prec_\LD w$.

 Let $l \in [k+1,\infty)$ be maximal with $w(l)>w(k)$ and set $v=wt_{kl}$. Then $\LD(v) = \LD(w) = k$, $v \prec_\LD w$, 
and $\ell^X(w) = \ell^X(v)+1$.  
Theorem~\ref{bcd-thm3} therefore expands $\cFX_w(\bfz)$ as a sum of terms  $\beta^{\ell^X(u)-\ell^X(w)} \cFX_u(\bfz)$
where $u$ is defined relative to $v$  as in Lemma~\ref{tech-lem2}, meaning   that 
$u\prec_\LD v\prec_\LD w $ and $u \in U $.
Hence, the property \eqref{ksp-eq} holds by induction.
\end{proof}

\begin{remark}
All named functions in this work are homogeneous when we set $\deg (\beta)=-1$ and $\deg (x_i) = \deg (y_i )= \deg (z_i)  =1$. 
For example, both $\fkGX_w(\bfx;\bfy;\bfz)$ and $\cFX_w(\bfz)$ are homogenous of degree $\ell^X(w)$.
As $\GP_\lambda(\bfz)$ and $\GQ_\lambda(\bfz)$ are homogeneous of degree $|\lambda|$,
asserting that a homogenous power series $f \in \ZZ[\beta]\llbracket \bfz \rrbracket$ belongs to $\GPP$ or $\GQP$
is equivalent to claiming that we can write
\be\label{fab-eq}\textstyle f = \sum_\lambda a_\lambda\cdot \beta^{|\lambda|-\deg(f)}\cdot \GP_\lambda(\bfz)
\quord
f = \sum_\lambda a_\lambda\cdot \beta^{|\lambda|-\deg(f)}\cdot \GQ_\lambda(\bfz)\ee
for nonnegative integers $a_\lambda \in \NN$ that are nonzero for only finitely many strict partitions $\lambda$.
This is consistent with \eqref{ksp-eq} since $\ell^X(w) = |\lambda_X(w)|$ when $w \in W^X_\infty$ is Grassmannian.

It is an open problem to find a combinatorial formula 
for the integers $a_\lambda$ when 
$f = \fkGX_w(\bfx;\bfy;\bfz)$. Recent work of Arroyo--Hamaker--Hawkes--Pan \cite[Conj.~1.3]{AHHP} gives a conjectural solution to this problem in type $\sfC$; compare with \cite[Thm.~1]{BKSTY} in type $\sfA$.
\end{remark}

We conclude this discussion with a few consequences of 
Theorem~\ref{k-stanley-positivity-thm}.
Since we can only have $\sigma^{-1} \circ u \circ \tau  \in W^X_n$ for $u \in W^X_n$ and $\sigma,\tau \in W^A_\infty$
if $\sigma$ and $\tau$ are in the image of $S_n \hookrightarrow W^A_\infty$,
we deduce the next corollary from  \eqref{nn-eq} and Definition~\ref{cont-def}.
This generalizes \cite[Thm.~8.13]{IkMiNa}.

 \begin{corollary}
 If $v\in \WB_n=\WC_n$ and $w \in \WD_n$ then 
 \[\ba
 \fkGB_v(\bfx;\bfy;\bfz) &\in R\spanning\left\{\GP_\lambda(\bfz) :\lambda\text{ strict}\right\} \\
  \fkGC_v(\bfx;\bfy;\bfz) &\in R\spanning\left\{\GQ_\lambda(\bfz) :\lambda\text{ strict}\right\} \\
   \fkGD_w(\bfx;\bfy;\bfz) &\in R\spanning\left\{\GP_\lambda(\bfz) :\lambda\text{ strict}\right\}
 \ea
\quad\text{for $R:=\NN[\beta][x_1,x_2,\dots,x_{n-1},y_1,y_2,\dots,y_{n-1}]$.}\]
 \end{corollary}
 
Combining Theorems~\ref{fc-prop} and \ref{k-stanley-positivity-thm} implies the following shifted analogue of \cite[Thm.~6.9]{Buch2002}. This result was predicted as \cite[Conj.~5.14]{LM}, and is equivalent to Corollary~\ref{coprod-cor} in the introduction.

\begin{corollary}\label{skew-cor}
If $\lambda \subseteq \nu$ are strict partitions
then $\GP_{\nu/\lambda}(\bfz)  \in \GPP$ and $\GQ_{\nu/\lambda}(\bfz)  \in \GQP$.
\end{corollary}

\begin{remark}
As explained in the remark after \cite[Problem~4.37]{MarbergHopf},
the preceding corollary also confirms \cite[Conjs.~4.35 and 4.36]{MarbergHopf}. These conjectures assert that the coproduct operation $\Delta$ for the \defn{(linearly compact) Hopf algebra} $\QSym_\beta$ (see \cite[\S2]{LM}) restricts to well-defined maps 
\[ \Delta : \GPP\to \GPP\otimes_{\ZZ[\beta]} \GPP
\quand 
 \Delta : \GQP\to \GQP\otimes_{\ZZ[\beta]} \GQP.
 \]
 The analogous property for the product operation $\nabla$ on $\QSym_\beta$ holds by \cite[Thm.~1.6]{LM2}.
\end{remark}

Suppose $\lambda$ is a partition with $n$ nonzero parts. Let $\delta = (n,n-1,\dots,3,2,1)$
and define 
\be
\GS_\lambda(\bfz) = \GP_{(\lambda+\delta)/\delta}(\bfz) = \GQ_{(\lambda+\delta)/\delta}(\bfz)
\quad
\text{where }\lambda+\delta = (\lambda_1+n,\lambda_2+n-1,\dots,\lambda_n+1).\ee
These functions are of interest 
since the map $\Theta$ from \eqref{Theta-eq} sends 
$ \Theta(G_\lambda(\bfz)) = \GS_\lambda(\bfz)$
for all partitions $\lambda$
by \cite[\S4.6]{LM}; see also \cite[Cor.~5.17]{LM}. The following is a special case of Corollary~\ref{skew-cor}:

\begin{corollary}
If $\lambda$ is any partition then $\GS_{\lambda}(\bfz) \in \GPP\cap \GQP$.
\end{corollary}


\begin{thebibliography}{99}
 
 

\bibitem{AGM} Dave Anderson, Stephen Griffeth, and Ezra Miller,
Positivity and Kleiman transversality in equivariant $K$-theory of homogeneous spaces,
\emph{J. Eur. Math. Soc.} 13 (2011), no. 1, pp. 57--84.

\bibitem{AHHP} Joshua Arroyo, Zachary Hamaker, Graham Hawkes, and Jianping Pan,
Type C $K$-Stanley symmetric functions and Kra\'skiewicz-Hecke insertion, preprint (2025), {\tt arXiv:2503.16641}.



\bibitem{Billey} Sara Billey, Transition equations for isotropic flag manifolds, \emph{Discrete Math.} \textbf{193} (1998), 69--84.

\bibitem{BilleyHaiman} Sara Billey and Mark Haiman,
Schubert polynomials for the classical groups,
\emph{J. Amer. Math. Soc.} \textbf{8} (1995), 443--482.

\bibitem{BL} Sara Billey and Tao Kai Lam  Lam, Vexillary elements in the hyperoctahedral group, \emph{J. Alg. Combin.} \textbf{8}
(1998), no. 2, 139--152.

\bibitem{Buch2002} Anders Skovsted Buch, A Littlewood-Richardson rule for the $K$-theory of Grassmannians, \emph{Acta Math.} \textbf{189} (2002), no. 1, 37--78.

\bibitem{BKSTY} Anders Skovsted Buch, Andrew Kresch, Mark Shimozono, Harry Tamvakis, and Alexander Yong,
Stable Grothendieck polynomials and $K$-theoretic factor sequences, 
\emph{Math. Ann.} \textbf{340} (2) (2008), 359--382.

\bibitem{BucScr} Valentin Buciumas and Travis Scrimshaw, 
Double Grothendieck polynomials and colored lattice models,
\emph{Int. Math. Res. Not.}(2020).


\bibitem{ChiuMarberg} Yu-Cheng Chiu and Eric Marberg,
Expanding $K$-theoretic Schur $Q$-functions,
\emph{Algebr. Comb.} \textbf{6} (2023), no. 6, 1419--1445.

\bibitem{CTY} Edward Clifford, Hugh Thomas, and Alexander Yong, $K$-theoretic Schubert calculus for OG$(n, 2n + 1)$ and jeu de taquin for shifted increasing tableaux, \emph{J. Reine Angew. Math.} \textbf{690} (2014) 51--63.


\bibitem{FK1994}
Sergey Fomin and Anatol N. Kirillov, Grothendieck polynomials and the Yang-Baxter equation, Proceedings of the
Sixth Conference in Formal Power Series and Algebraic Combinatorics, DIMACS (1994), 183--190.

\bibitem{FK1996} Sergey Fomin and Anatol Kirillov,
Combinatorial $B_n$-analogues of Schubert polynomials,
Transactions of the American Mathematical Society \textbf{348} (1996), no. 9, 3591--3620.

\bibitem{FL} William Fulton and Alain Lascoux,
A Pieri formula in the Grothendieck ring of a flag bundle,
\emph{Duke Math. J.} \textbf{76} (1994), 711--729.


\bibitem{GK} William Graham and Shrawan Kumar, 
On Positivity in $T$-Equivariant $K$-Theory of Flag Varieties,
\emph{IMRN} (2008), rnn093.

\bibitem{GR} Stephen Griffeth and Arun Ram,
Affine Hecke algebras and the Schubert calculus, 
\emph{European J. Combin.} \textbf{25} (2004), no. 8, 1263--1283.


\bibitem{Hudson} 
Thomas Hudson,
A Thom-Porteous formula for connective K-theory using algebraic cobordism,
\emph{Journal of K-theory} \textbf{14} (2014), no. 2, 343--369.

\bibitem{HIMN} Thomas Hudson, Takeshi Ikeda, Tomoo Matsumura, and Hiroshi Naruse,
Double Grothendieck polynomials for symplectic and odd orthogonal Grassmannians,
\emph{J. Algebra} \textbf{546} (2020), 294--314.

\bibitem{IN}
Takeshi Ikeda and Hiroshi Naruse,
$K$-theoretic analogues of factorial Schur $P$- and $Q$-functions,
\emph{Adv. Math.} \textbf{243} (2013), 22--66.

\bibitem{IkMiNa} Takeshi Ikeda, Leonardo C. Mihalcea and Hiroshi Naruse, 
Double Schubert polynomials for the classical groups, 
\emph{Adv. Math.} \textbf{226} (2011), 840--886.

\bibitem{Kirillov}  Anatol N. Kirillov,
On double Schubert and Grothendieck polynomials for classical groups,
preprint (2015), {\tt arXiv:1504.01469}.

\bibitem{KN} Anatol N. Kirillov and Hiroshi Naruse,
Construction of double Grothendieck polynomials of classical types using IdCoxeter algebras,
\emph{Tokyo J. Math.} \textbf{39} (2017), no. 3, 695--728.


\bibitem{KnutsonMiller1}
Allen Knutson and Ezra Miller, Subword complexes in Coxeter groups,
\emph{Adv. Math.} \textbf{184} (2004), 161--176.


\bibitem{KMY} Allen Knutson, Ezra Miller, and Alexander Yong, Gr\"obner geometry of vertex decompositions and of
flagged tableaux, \emph{J. Reine Angew. Math.} \textbf{630} (2009), 1--31.
	
\bibitem{KV} Axel Kohnert and S\'ebastien Veigneau, 
Using Schubert basis to compute with multivariate polynomials,
\emph{Adv. Appl. Math.} \textbf{19} (1997), no. 1, 45--60.

\bibitem{KK} Bertram Kostant and Shrawan Kumar,
$T$-equivariant $K$-theory of generalized flag varieties,
\emph{J. Diff. Geom.} \textbf{32} (1990), 549--603.



\bibitem{TKLamThesis} Tao Kai Lam,  
$B$ and $D$ analogues of stable Schubert polynomials and related insertion algorithms, 
PhD thesis, MIT, 1995.

\bibitem{TKLam} Tao Kai Lam,  $B_n$ Stanley symmetric functions, \emph{Discrete Math.} \textbf{157} (1996), 241--270.

\bibitem{LLS} Thomas Lam, Seung Jin Lee, and Mark Shimozono,
Back stable $K$-theory Schubert calculus, 
\emph{IMRN} Volume 2023, Issue 24, 21381--21466.

\bibitem{LamPyl} Thomas Lam and Pavlo Pylyavskyy,
Combinatorial Hopf algebras and $K$-homology of Grassmannians,
\emph{IMRN} (2007), rnm125.


%
\bibitem{LascouxTransitions}
Alain Lascoux, 
Transition on Grothendieck polynomials, \emph{Physics and combinatorics}, 2000
(Nagoya), World Sci. Publishing, River Edge, NJ, 2001, 164--179.


\bibitem{Lascoux} Alain Lascoux, Chern and Yang through ice, Preprint, 2002.

\bibitem{LS1983} Alain Lascoux and Marcel-Paul Sch\"utzenberger, Symmetry and flag manifolds, in: Invariant Theory,
\emph{Lect. Notes in Math.} \textbf{996} (1983), 118--144.

\bibitem{LS1985} Alain Lascoux and Marcel-Paul Sch\"utzenberger, Schubert polynomials and the
Littlewood-Richardson rule, \emph{Lett. Math. Phys.} \textbf{10}  (1985), no. 2, 111--124.


\bibitem{Lenart} 
Cristian Lenart,
A $K$-theory version of Monk's formula and some related multiplication formulas,
\emph{J. Pure Appl. Alg.} \textbf{179} (2003), 137--158.

\bibitem{LP} Cristian Lenart and Alexander Postnikov, 
Affine Weyl groups in $K$-theory and representation theory,
\emph{IMRN} (2007), rnm038.


\bibitem{LM} Joel B. Lewis and Eric Marberg,
Enriched set-valued $P$-partitions and shifted stable Grothendieck polynomials,
\emph{Math. Z.} \textbf{299} (2021), 1929--1972.

\bibitem{LM2} Joel B. Lewis and Eric Marberg,
Combinatorial formulas for shifted dual stable Grothendieck polynomials,
\emph{Forum Math. Sigma} \textbf{12} (2024), Paper e22.

%
\bibitem{Manivel} Laurent Manivel, \emph{Symmetric Functions, Schubert Polynomials, and Degeneracy Loci},
American Mathematical Society, 2001.
%
%

\bibitem{MarbergHopf} Eric Marberg,
Shifted combinatorial Hopf algebras from $K$-theory
\emph{Algebraic Combinatorics} \textbf{7} (2024), no. 4, 1123--1156.

\bibitem{MarScr} Eric Marberg and Travis Scrimshaw,
Key and Lascoux polynomials for symmetric orbit closures,
preprint (2023), arXiv:2302.04226.



\bibitem{Matsumura} Tomoo Matsumura,
A tableau formula of double Grothendieck polynomials for 321-avoiding permutations, 
\emph{Ann. Comb.} \textbf{24} (2020), 55--67.

\bibitem{MNS} Leonardo C. Mihalcea, Hiroshi Naruse, and Changjian Su,
Left Demazure--Lusztig Operators on Equivariant (Quantum) Cohomology and $K$-Theory,
\emph{IMRN} Volume 2022, Issue 16, 12096--12147.


 \bibitem{PechenikYong} Oliver Pechenik and Alexander Yong,
  Genomic tableaux,
  \emph{J. Algebr. Comb.} 45 (2017), 649--685.

\bibitem{Stanley} Richard P. Stanley, On the number of reduced decompositions of elements of Coxeter groups,
\emph{European J. Combin.} \textbf{5} (1984), 359--372.


\bibitem{StembridgeEnriched}  John R. Stembridge,
Enriched $P$-partitions,
\emph{Trans. Amer. Math. Soc.} \textbf{349} (1997), no. 2, 763--788.

\bibitem{StembridgeSome}  John R. Stembridge,
Some combinatorial aspects of reduced words in finite Coxeter groups,
\emph{Trans. Amer. Math. Soc.}, \textbf{349} (1997), no. 4, 1285--1332.

\bibitem{Tamvakis} Harry Tamvakis, 
Tableau formulas for skew Grothendieck polynomials,
\emph{J. Math. Soc. Jpn.} \textbf{76} (2024), no. 1, 147--172.

\bibitem{Weigandt} Anna Weigandt, Bumpless pipe dreams and alternating sign matrices, \emph{J. Comb. Theory, Ser. A} \textbf{182} (2021), 105470.

\bibitem{Willems} Matthieu Willems,
A Chevalley formula in equivariant $K$-theory,
\emph{J. Algebra} \textbf{308} (2007), 764--779.


\end{thebibliography}
\end{document}